\mag=\magstep1
\documentclass{article}
\usepackage{amsmath,amssymb,amsthm}
\usepackage{amscd}
\makeatletter
\@addtoreset{equation}{section}

\makeatother
\makeatletter
\long\def\@makecaption#1#2{
  \vskip\abovecaptionskip
  \sbox\@tempboxa{#1~ #2}%
  \ifdim \wd\@tempboxa >\hsize
    #1~ #2\par
  \else
    \global \@minipagefalse
    \hb@xt@\hsize{\hfil\box\@tempboxa\hfil}%
  \fi
  \vskip\belowcaptionskip}
\makeatother
\makeatletter
\renewcommand\section{\@startsection {section}{1}{\z@}%
                                   {-3.5ex \@plus -1ex \@minus -.2ex}%
                                   {2.3ex \@plus.2ex}%
                                   {\normalfont\large\bfseries}}
\renewcommand\subsection{\@startsection{subsection}{2}{\z@}%
                                     {-3.25ex\@plus -1ex \@minus -.2ex}%
                                     {1.5ex \@plus .2ex}%
                                     {\normalfont\normalsize\bfseries}}
\renewcommand\subsubsection{\@startsection{subsubsection}{3}{\z@}%
                                     {-3.25ex\@plus -1ex \@minus -.2ex}%
                                     {1.5ex \@plus .2ex}%
                                     {\normalfont\normalsize\bfseries}}
\makeatother
\newcommand{\OP}[1]{\operatorname{#1}}
\newcommand{\BF}[2]{\langle #1,#2 \rangle} 
\newcommand{\LF}[1]{\langle #1 \rangle}
\newcommand{\MC}[1]{\mathcal{#1}}
\newcommand{\MF}[1]{\mathfrak{#1}}
\newcommand{\Z}{\mathbb{Z}}
\newcommand{\Q}{\mathbb{Q}}
\newcommand{\R}{\mathbb{R}}
\newcommand{\C}{\mathbb{C}}
\renewcommand{\P}{\mathbb{P}}
\newcommand{\T}{\mathsf{t}}
\newcommand{\X}{\mathcal{X}}
\newcommand{\RL}{\rule{0ex}{1.8ex}}
\def\rank{\operatorname{rank}}
\newcommand{\disc}{\operatorname{disc}}
\newcommand{\ord}{\operatorname{ord}}
\newcommand{\ORDa}{\OP{Cl}^+(T^*)}
\newcommand{\ORDb}{\overline{\OP{Cl}^+}(T^*)}
\newcommand{\GAMa}{\Gamma^1(\ORDb)[2]}
\newcommand{\GAMb}{\Gamma^1(\ORDb)}
\newcommand{\GAMc}{\Gamma'}
\newcommand{\GAMd}{\Gamma^+(\ORDb)}
\newcommand{\DOMb}{\MC{D}(\ORDb,I)}
\newcommand{\HALF}{\frac{1}{2}}
\newcommand{\XXXa}{\theta^4_{0000}}
\newcommand{\XXXb}{\theta_{000\HALF}^4}
\newcommand{\XXXc}{\theta_{00\HALF0}^4}
\newcommand{\XXXd}{\theta_{00\HALF\HALF}^4}
\newcommand{\XXXe}{\theta_{0\HALF00}^4}
\newcommand{\SYMG}{{\mathfrak{S}_5}}
\newcommand{\PRIM}{\text{\rm prim}}

\newcommand{\CONGL}{\mathop{\longrightarrow}^\sim}
\newcommand{\PHIG}{\Phi_{\rm grp}}
\newcommand{\PHID}{\Phi_{\rm dom}}
\newcommand{\PER}{\mathit{per}}
\DeclareMathSizes{10}{10}{7}{7}
\begin{document}
\theoremstyle{plain}
\newtheorem{thm}{Theorem}[section]
\newtheorem*{thmn}{Theorem}
\newtheorem{thmdash}[thm]{Theorem$'$}
\newtheorem{lem}[thm]{Lemma}
\newtheorem{cor}[thm]{Corollary}
\newtheorem{prop}[thm]{Proposition}
\newtheorem{clm}[thm]{Claim}
\newtheorem{ex}[thm]{Example}
\theoremstyle{definition}
\newtheorem{defn}[thm]{Definition}
\newtheorem{rem}[thm]{Remark}
\newenvironment{itemize2}
 {%
  \begin{itemize}%
  \setlength{\itemsep}{0pt}%
  \vspace{-1ex}%
 }%
 {%
  \vspace{-0.5ex}%
  \end{itemize}%
 }
\title{Period Map of a Certain %
$K3$ Family with an $\mathfrak{S}_5$-Action}
\author{Kenji Hashimoto\footnote{
This work was supported by Grant-in-Aid for JSPS Fellows No.20-56181.}}
\date{\today}

\maketitle
\setcounter{section}{-1}
%
%

\begin{abstract}
In this paper, we study the period map of a certain one-parameter
 family of quartic
 $K3$ surfaces
 with an $\MF{S}_5$-action.
We construct automorphic forms on the period domain as the pull-backs
 of theta constants of genus $2$ by a modular embedding.
Using these automorphic forms,
 we give an explicit presentation of the inverse period map.
\end{abstract}

\section{Introduction}
A compact complex surface $X$ is called a $K3$ surface
if it is simply connected
 and has a nowhere vanishing holomorphic $2$-form
$\omega_X$.
A nonsingular quartic surface in $\P^3$ is
 a typical example of $K3$ surface.
Period maps play an important role in the study on $K3$ surfaces.
In this paper, we study the period map
of the family $\{\MC{X}_\T \}_{\T=(t_0:t_1)\in\P^1}$
of quartic $K3$ surfaces defined by
\begin{equation} \label{defaffine}
\MC{X}_\T:x_1+\cdots+x_5
=t_0(x_1^4+\cdots+x_5^4)+t_1(x_1^2+\cdots+x_5^2)^2=0
\end{equation}
in $\P^4$
 with homogeneous coordinates
 $(x_1:\cdots:x_5)$.
We give an explicit presentation of the inverse
 of the period map using
automorphic forms (see below).
The symmetric group $\mathfrak{S}_5$ of degree $5$ acts
on the family $\{ \MC{X}_\T \}$ by the permutation of the variables
$x_1,\ldots,x_5$.

This family $\{ \X_\T \}$ is an analogue
 of the Hasse family of
 elliptic curves for $K3$ surfaces.
For the Hasse family, the inverse of the period map
 is described explicitly
 using theta constants as follows
 (see \cite{krazer70}).
Let $(X_1:X_2:X_3)$ be
 homogeneous coordinates
 of $\P^2$.
The symmetric group $\MF{S}_3$
 of degree $3$ acts on $\P^2$
 by the permutation of
 the variables $X_1,X_2,X_3$.
Let $G$ be the subgroup of
 $\OP{Aut}(\P^2)$ generated by
 $\MF{S}_3$ and the transformation
 $(X_1:X_2:X_3)\mapsto
 (X_1:\zeta X_2:\zeta^2 X_3)$,
 where $\zeta=\exp(2\pi i/3)$.
Then, a smooth cubic curve
 in $\P^2$ stable under the action of $G$ is
 of the form
\begin{equation}
 C_\mu:X_1^3+X_2^3+X_3^3-3\mu X_1X_2X_3=0,\quad
 \mu\in\C\smallsetminus
 \{ 1,\zeta,\zeta^2 \}.
\end{equation}
Let $\omega_\mu$ be
 a nowhere vanishing holomorphic $1$-form on $C_\mu$.
We choose $1$-cycles $\alpha_\mu,\beta_\mu$
 on $C_\mu$
 which give a symplectic basis of $H_1(C_\mu,\Z)$.
Then, $\tau_\mu
 :=\int_{\beta_\mu}\omega_\mu
 / \int_{\alpha_\mu}\omega_\mu$
 is an element in the upper half-plane $\MF{H}$.
The $\Gamma(3)$-orbit
 $[\tau_\mu]:=\Gamma(3)\cdot\tau_\mu$
 of $\tau_\mu$
 depends only on $\mu$,
 not on the choice of $\alpha_\mu,\beta_\mu$,
 where $\Gamma(3)$ is the principal
 congruence subgroup of level $3$.
By attaching $[\tau_\mu]$ to $\mu$,
 we have the period map
\begin{equation}
 \C\smallsetminus \{ 1,\zeta,\zeta^2 \} \ni\mu
 \longmapsto[\tau_\mu]\in\Gamma(3)\backslash\MF{H}.
\end{equation}
In fact, the period map is bijective
 and the inverse period map
 is given by theta constants as follows.
The theta constant $\theta_m$
 for a character $m=(m',m'')\in\Q^2$ is
 defined by
\begin{equation}
 \theta_m(\tau)
 =\sum_{n\in\Z}
 \exp \left( 2\pi i \left(
  \frac{1}{2}(n+m')\cdot\tau\cdot
  (n+m')+
  (n+m')\cdot
  m''\right)
 \right)
\end{equation}
for $\tau\in\MF{H}$.
Then, we have
\begin{equation} \label{formula_a}
 \mu=
  \frac
   {
    \theta_{0,0}(\tau_\mu)^3
    +2\,\theta_{0,\frac{1}{3}}(\tau_\mu)^3
   }
   {
    \theta_{\frac{1}{3},0}(\tau_\mu)^3
    +\zeta^2\,\theta_{\frac{1}{3},\frac{1}{3}}(\tau_\mu)^3
    +\zeta\,\theta_{\frac{1}{3},-\frac{1}{3}}(\tau_\mu)^3
   }.
\end{equation}
In this paper, we find
 an analogous formula of
 (\ref{formula_a}) for the family
 $\{ \X_\T \}$.

We study this particular family
 $\{ \X_\T \}$
 from the view point of
 finite symplectic actions on $K3$ surfaces
 after Nikulin \cite{nikulin79fin}.
An automorphism $g$ of a $K3$ surface $X$ is said to be
symplectic if $g^\ast \omega_X=\omega_X$.
The restriction of the action of $\MF{S}_5$
 to $\MF{A}_5$
on each $\X_\T$ is symplectic.
We review the history of studies
 on symplectic actions of finite
 groups on $K3$ surfaces.

Nikulin \cite{nikulin79fin} classified
 symplectic actions of finite abelian groups
 on $K3$ surfaces.
After Nikulin's result,
Mukai \cite{mukai88} classified finite groups
 which can be realized as symplectic actions
 on a $K3$ surface.
He listed all maximal finite groups
 with this property.
Another proofs of this classification are given
 by Xiao \cite{xiao96}
 and Kond\=o \cite{kondo98}.
Kond\=o's method is to embed the coinvariant part of the
$K3$ lattice into a Niemeier lattice, i.e.\ definite
even unimodular lattice of rank $24$,
and yo describe symplectic actions
 as automorphisms of the Niemeier lattice.
Kond\=o also used this method to
study symplectic actions of the Mathieu group $M_{20}$
and showed that there exists a unique $K3$
surface with a non-symplectic action of
 $M_{20}.C_4$,
 i.e.\ an extension of $M_{20}$
 by the cyclic group $C_4$ of order $4$.
For similar results for $L_2(7),\mathfrak{A}_6,F_{384}$,
we refer to \cite{oguisozhang02,keumoguisozhang05,oguiso05}.
These four groups appear in the Mukai's list of maximal finite
symplectic actions.
For a maximal finite group $G$
 which appears in the Mukai's list,
 the moduli of algebraic $K3$ surfaces with a symplectic
action of $G$ is $0$-dimensional (countably infinite).
If an algebraic $K3$ surface $X$ admits a maximal finite symplectic action,
the Picard number of $X$ is equal to $20$
and such $X$ is determined uniquely by its transcendental lattice
with the orientation.

On the other hand,
 the symplectic action of $\MF{A}_5$ on
 each $\X_\T$ is not maximal.
If a $K3$ surface $X$ admits
 a symplectic action of $\mathfrak{A}_5$,
the Picard number of $X$ is either $19$ or $20$.
As a consequence,
the moduli of algebraic $K3$ surfaces
 with a symplectic
$\mathfrak{A}_5$-action is $1$-dimensional,
 and
 the family $\{ \X_\T \}$ is a maximal family
 of quartic $K3$ surfaces
 with an $\mathfrak{S}_5$-action
 (see Remark~\ref{REMmaximality}).

We recall the period map for polarized $K3$ surfaces.
A pair $(X,l)$ of a $K3$ surface $X$
and an ample class $l\in H^2(X,\Z)$
is called a polarized $K3$ surface of degree $d$
if $l$ is primitive and $\BF{l}{l}=d$.
Let $\Lambda=E_8[-1]^{\oplus 2}\oplus U^{\oplus 3}$ be the
$K3$ lattice.
We fix
 a primitive vector $\delta\in\Lambda$
 with $\BF{\delta}{\delta}=d$.
Then we have $\Lambda_\PRIM:=(\delta)^\bot_\Lambda
\cong E_8[-1]^{\oplus 2}\oplus U^{\oplus 2}\oplus\LF{-d}$.
A marking of a polarized $K3$ surface $(X,l)$
 is an isomorphism
 $\alpha:H^2(X,\Z)\rightarrow \Lambda$ of lattices
 such that $\alpha(l)=\delta$.
The triple $(X,l,\alpha)$ is called a marked polarized $K3$ surface.
For a lattice $L$ of signature $(2,\rho)$,
 we define the symmetric domain $\Omega_L$ by
\begin{equation}
 \Omega_{L}
 =\{\C\omega\in\P(L\otimes\C) \bigm|
 \BF{\omega}{\omega}=0,\BF{\omega}{\bar{\omega}}>0
 \}.
\end{equation}
This is a bounded symmetric domain of type IV.
The period
 of a marked polarized $K3$ surface
 $(X,l,\alpha)$ is defined as
 the point
 $(\alpha\otimes\C)(H^{2,0}(X))\in\Omega_{\Lambda_\PRIM}$.
The period map
\begin{equation}
 (X,l,\alpha)\mapsto (\alpha\otimes\C)(H^{2,0}(X))
\end{equation}
 gives a one-to-one correspondence
 between
 the set of isomorphism classes of marked polarized $K3$ surfaces
 of degree $d$
 and $\Omega_{\Lambda_\PRIM}$
 by Torelli theorem and the surjectivity of the period map.
By forgetting markings,
 this correspondence induces a one-to-one correspondence
 between the set of isomorphism classes of polarized $K3$ surfaces
 of degree $d$ and
 $\OP{Aut}(\Lambda_\PRIM,\delta)\backslash\Omega_{\Lambda_\PRIM}$,
 where
 $\OP{Aut}(\Lambda_\PRIM,\delta)$
 is the subgroup of $\OP{O}(\Lambda)$ consisting of
 elements fixing $\delta$.
Since the period domain $\Omega_{\Lambda_\PRIM}$
 is a bounded symmetric domain of type IV,
the inverse period map can be written by automorphic forms.
However, since the dimension of $\Omega_{\Lambda_\PRIM}$ is $19$,
it will be difficult
 to write down these automorphic forms concretely.
On the other hand,
 the Grammian matrix of the transcendental lattice $T$
 of a generic fiber of our family $\{ \X_\T \}$
 is equal to
\begin{equation}
 \begin{pmatrix} 4 & 1 & 0 \\ 1 & 4 & 0 \\
 0 & 0 & -20 \end{pmatrix},
\end{equation}
 and the image of the period map for our family is
 a $1$-dimensional subdomain $\Omega_T$
 of $\Omega_{\Lambda_\PRIM}$.
Let $\omega_\T$ be a nowhere vanishing holomorphic $2$-form
 on $\X_\T$.
By attaching
 $[\C\omega_\T]
 :=\OP{O}(T)\cdot\C\omega_\T\in\OP{O}(T)\backslash\Omega_T$,
 we have a map
\begin{equation}
\P^1\ni\T\longmapsto[\C\omega_\T]\in\OP{O}(T)\backslash\Omega_T,
\end{equation}
 which is called the period map
 for the $K3$ family $\{ \X_\T \}_{\T\in\P^1}$.
The inverse period map for our family can be written by
 one-variable automorphic forms
 on a connected component $\Omega_T^\circ$ of $\Omega_T$.
Since the Fuchsian group $\OP{O}(T)$ has no cusp,
 we use a modular embedding
 $i:\Omega_T^\circ\rightarrow\MF{H}_2$
 to the Siegel upper half-space $\MF{H}_2$ of genus $2$
 to construct automorphic forms on $\Omega_T^\circ$.
The pull-backs of the fourth powers of theta constants
 of genus $2$ (see (\ref{theta}))
\begin{equation}
 \begin{array}{c}
 X_1=i^*\theta_{0000}^4,~
 X_2=i^*\theta_{000\HALF}^4,~
 X_3=i^*\theta_{00\HALF0}^4, \vspace{1ex} \\
 X_4=i^*\theta_{00\HALF\HALF}^4,~
 X_5=i^*\theta_{0\HALF00}^4
 \end{array}
\end{equation}
 become automorphic forms on $\Omega_T^\circ$
 of a finite index subgroup of $\OP{O}(T)$.
The main theorem of this paper is the following.

\begin{thmn}
Using the above automorphic forms $X_1,\ldots,X_5$,
 the parameter $\T=(t_0:t_1)\in\P^1$ is expressed by
\begin{equation}
 \frac{t_1}{t_0}
 =-\frac{1}{5}
 \left(
 \frac
 {s(\omega_\T)-1}
 {s(\omega_\T)+2}
 \right)^2-\frac{1}{4},
\end{equation}
where
\begin{gather*}
 r_1=\frac
 {(X_2-X_3+2X_5)(2X_1-X_2+X_3-2X_4-2X_5)}
 {(X_2+X_3)(2X_1-X_2-X_3+2X_4)}, \\
 r_2=\frac
 {(X_2+X_3)^2}
 {(2X_1-X_2-X_3+2X_4)^2}, \\
 s=\frac{324 r_2 - 389 r_1^2
 + 458 r_1 - 325}{r_1(340 r_1 - 20)}.
\end{gather*}
\end{thmn}

This paper proceeds as follows.
We recall basic facts on $K3$ surfaces
 and lattices in Section~\ref{SECTnotations}.
In Section~\ref{SECTs5},
 we study symplectic actions of $\mathfrak{S}_5$
 on $K3$ surfaces
 using Kond\=o's method.
To study actions of $\MF{S}_5$ on the second cohomology groups of $K3$
surfaces,
we embed the coinvariant lattice $H(X)_{\MF{S}_5}$
 into a Niemeier lattice.
The goal of this section is
 to classify the transcendental lattices of
 $K3$ surfaces
 with certain $\MF{S}_5\times C_2$-actions.
In Section~\ref{SECTfamily},
 we determine singular fibers
 of the family $\{ \X_\T \}$.
We construct a smooth birational model after local base change.
We remark that the fiberwise symplectic action of $\MF{A}_5$
 extends to this model.
In Section~\ref{SECTtrans},
 we determine
 the transcendental lattice $T$
 of a generic fiber of the family.
For this purpose,
 we study a special fiber of the family,
 which has an extra symmetry.
The minimal desingularization of
 this fiber has
 a symplectic action of $\mathfrak{S}_5$
 extending
 that of $\mathfrak{A}_5$.
We apply the result in Section~\ref{SECTs5}
 to determine the transcendental lattice of
 this fiber.
By the modification in Section~\ref{SECTfamily},
 we determine $T$
 from the information on this fiber.
In Section~\ref{SECTquaternion},
 we recall basic facts on quaternion algebras
 and study some relations
 between $\OP{O}(T)$ 
 and the unit groups of orders
 in a quaternion algebra $\OP{Cl}^+(T^\vee[12])\otimes\Q$,
 which is the even part of
 a Clifford algebra $\OP{Cl}(T^\vee[12])\otimes\Q$.
It follows that the space $\OP{O}(T)\backslash\Omega_T$
 is compact.
To obtain automorphic forms on $\Omega^\circ_T$,
 we construct a modular embedding
 $i:\Omega^\circ_T\rightarrow\MF{H}_2$
 by Kuga--Satake construction
 in Section~\ref{SECTmodular}.
The pull-backs of the fourth powers of theta constants
 on $\mathfrak{H}_2$
 by the modular embedding $i$
 become automorphic forms on $\Omega_T^\circ$.
In Section~\ref{SECTperiod},
we give an explicit formula of the inverse period map by using these
 automorphic forms.
We find relations of
 these automorphic forms to obtain a uniformizing parameter of the
 Shimura curve $\OP{O}(T)\backslash\Omega_T$.
To get enough relations of automorphic forms, we compute special values
 of them at several CM points, which is discussed in Appendix.

After finishing this work, we noticed a relevant paper \cite{smith06},
which treats the Picard-Fuchs equation for this $K3$ family.

\subsubsection*{Acknowledgement}
The author would like to express his thanks to the supervisor Professor Tomohide
Terasoma
for suggesting this problem.
The author also would like to express his thanks to Professor Keiji
Matsumoto and Professor Keiji Oguiso
for comments and encouragements.

\section{Some basic facts on $K3$ surfaces and lattice theory}
\label{SECTnotations}
In this paper, a free $\Z$-module $L$
of finite rank
equipped with an integral symmetric bilinear form
$\langle~,~\rangle$
is simply called a lattice.
An automorphism of $L$ is defined
 as a $\Z$-automorphism of $L$
 preserving $\BF{~}{~}$.
We denote by $\OP{O}(L)$ the group of automorphisms
of $L$.
If $\BF{x}{x}\equiv 0 \bmod 2$ for all $x\in L$,
a lattice $L$ is said to be even.
For an integral symmetric $n\times n$ matrix $M$,
the lattice $(\Z^n,\BF{~}{~})$ defined by
$\BF{x}{y}={}^t x M y$ is also denoted by $M$.
The lattice of rank $1$ generated by $x$
with $\BF{x}{x}=a$ is denoted by $\LF{a}$.
The discriminant $\disc(L)$ of $L$ is
the determinant of the Grammian matrix of $L$.
If $\disc(L)=\pm 1$, a lattice $L$ is said to be unimodular.
We denote by $L[\lambda]$ the lattice $L$ equipped with
$\lambda$ times the bilinear form $\BF{~}{~}$,
i.e.\ $(L,\lambda\BF{~}{~})$.
A sublattice $K$ of $L$ is said to be primitive
 if $L/K$ is torsion-free.
We assume that an action of a group $G$ on $L$ preserves
the bilinear form unless otherwise stated.
For a lattice $L$ with an action of $G$,
 we define the invariant part $L^G$
 and the coinvariant part $L_G$
 of $L$ by
\begin{equation}
L^G=\{x\in L\bigm|g\cdot x=x ~ (\forall g\in G) \},\quad
L_G=(L^G)^\bot_L.
\end{equation}

An even unimodular lattice of signature $(3,19)$
is isomorphic to
$\Lambda:=E_8[-1]^{\oplus 2}\oplus U^{\oplus 3}$,
which is called the $K3$ lattice.
Here $E_8$ (resp.\ $U$) is the root lattice of type $E_8$
 (resp.\ the hyperbolic lattice),
 which is a unique even unimodular lattice
 of signature $(8,0)$ (resp.\ (1,1)) up to isomorphism.
A compact complex surface $X$ is called a $K3$ surface
if it is simply connected
and has a nowhere vanishing holomorphic $2$-form $\omega_X$.
(See~\cite{bhpv}.)
The second cohomology group $H(X):=H^2(X,\Z)$
 of a $K3$ surface $X$
 with its intersection form
 is isomorphic to the $K3$ lattice $\Lambda$.
Since any $K3$ surface $X$ is a K\"ahler manifold,
we have the Hodge decomposition
$H(X)\otimes\C=H^{2,0}(X)\oplus H^{1,1}(X)\oplus H^{0,2}(X)$.
The space $H^{2,0}(X)\cong\C$ is generated by
the class $[\omega_X]$ of $\omega_X$.
The algebraic lattice $S(X)$
 and the transcendental lattice $T(X)$ of
$X$ are defined by
\begin{equation}
S(X)=\{ x\in H(X) \bigm| \langle x,[\omega_X] \rangle=0 \},
\quad T(X)=(S(X))^\bot_{H(X)}.
\end{equation}
Here we extend the bilinear form $\BF{~}{~}$ on $H(X)$ to
that on $H(X)\otimes\C$ linearly.
The open subset $\MC{C}_X\subset H^2(X,\R)$
of classes of K\"ahler $(1,1)$-forms 
is called the K\"ahler cone of $X$.
The study of automorphims of $K3$ surfaces is reduced to
lattice theory using the following theorems.

\begin{thm}[Torelli theorem for $K3$ surfaces,
cf.~\cite{bhpv}] \label{THMtorelli}
Let $X$ and $X'$ be $K3$ surfaces.
Let $\varphi:H(X)\rightarrow H(X')$
be an isomorphism of lattices
satisfying the following conditions:
\begin{itemize2}
\item[(i)]
$\varphi\otimes\C:H^2(X,\C)\rightarrow H^2(X',\C)$ preserves
the Hodge structure;
\item[(ii)]
$(\varphi\otimes\R)(\MC{C}_X)=\MC{C}_{X'}.$
\end{itemize2}
Then, there exists a unique isomorphism $f:X'\rightarrow X$
such that $f^*=\varphi$.
Moreover, under the condition (i),
the condition (ii) is equivalent to
\begin{itemize2}
\item[(ii)$'$]
there exists a vector $\kappa\in\MC{C}_X$
 such that $(\varphi\otimes\R)(\kappa)\in\MC{C}_{X'}$.
\end{itemize2}
\end{thm}

\begin{thm}
[The surjectivity of the period map for $K3$ surfaces,
cf.~\cite{bhpv}] \label{THMsurjperiod}
\rule{0ex}{0ex}\\
Assume that vectors $\omega\in\Lambda\otimes\C$
and $\kappa\in\Lambda\otimes\R$ satisfy
the following conditions:
\begin{itemize2}
\item[(i)]
$\BF{\omega}{\omega}=0$, $\BF{\omega}{\bar{\omega}}>0$;
\item[(ii)]
$\BF{\kappa}{\omega}=0$, $\BF{\kappa}{\kappa}>0$;
\item[(iii)]
 $\BF{\kappa}{x}\neq 0$
 for all $x\in(\omega)^\bot_\Lambda$
 such that $\BF{x}{x}=-2$.
\end{itemize2}
Then, there exist a $K3$ surface $X$ and an isomorphism
$\alpha:H(X)\rightarrow\Lambda$
 between lattices such that
$\C\omega=(\alpha\otimes\C)(H^{2,0}(X))$
and $\kappa\in(\alpha\otimes\R)(\MC{C}_X)$.
\end{thm}

If a group $G$ acts on a $K3$ surface $X$,
the action of $G$ induces a left action on $H(X)$ by
\begin{equation} \label{defaction}
g\cdot x=(g^{-1})^\ast x,\quad g\in G,\, x\in H(X).
\end{equation}
An automorphism $g$ of $X$ is said to be
symplectic
if $g^\ast\omega_X=\omega_X$.
To study symplectic actions of finite groups on $K3$ surfaces,
we use the following theorems.

\begin{thm}[\cite{nikulin79fin}] \label{THMsymplectic}
Let $G$ be a finite group and $X$ a $K3$ surface.
If $G$ acts on $X$ symplectically,
we have the following:
\begin{itemize2}
\item[(i)]
$(H(X))_G$ is negative definite;
\item[(ii)]
$x^2\neq -2$ for all $x\in (H(X))_G$;
\item[(iii)]
$T(X)\subset H(X)^G$.
\end{itemize2}
Conversely, if $G$ acts on $H(X)$
satisfying the conditions (i)--(iii),
there exists an element $\alpha\in\OP{O}(H(X))$
such that the action of $\alpha \circ G \circ \alpha^{-1}$
(i.e.\ the conjugation by $\alpha$)
on $H(X)$ is induced by a symplectic action of $G$ on $X$.
\end{thm}

\begin{thm}[\cite{nikulin79fin,mukai88}] \label{THMfixedpts}
Let $g$ be a symplectic automorphism
of a $K3$ surface $X$ of finite order
and $X^g$ the set of fixed points under the action of $g$.
Then, we have $\ord(g) \leq 8$ and the cardinality $\sharp X^g$ of $X^g$
depends only on $\ord(g)$ as follows:
\begin{center}
\begin{tabular}{c|ccccccc}
$\ord(g)$ & $2$ & $3$ & $4$ & $5$ & $6$ & $7$ & $8$ \\ \hline
$\sharp X^g$\rule{0ex}{2.2ex} & $8$ & $6$ & $4$ & $4$ & $2$ & $3$ & $2$
\end{tabular}\lower1.5ex\hbox{.}
\end{center}
Moreover, $\sharp X^g$ is equal to the trace of $g$ on
$H^\ast(X,\Q)=
H^0(X,\Q)\oplus H^2(X,\Q)\oplus H^4(X,\Q)
\cong\Q^{24}$.
(This is a consequence of Lefschetz fixed point formula.)
\end{thm}

As a consequence
of Theorems~\ref{THMtorelli}--\ref{THMfixedpts},
we have the
following (cf.\ \cite{nikulin79fin,mukai88}).

\begin{lem} \label{LEMtrace}
Let $g$ be an element in $O(\Lambda)$.
Assume that the coinvariant lattice $\Lambda_{\langle g \rangle}$
is negative definite
and contains no vector $x$ such that $\BF{x}{x}=-2$.
Then, we have $\ord(g)\leq 8$
and $\OP{tr}(g)+2=24,8,6,4,4,2,3,2\,$
if $\,\ord(g)=1,2,3,4,5,6,7,8$, respectively.
\end{lem}

We recall some basic properties on discriminant forms of lattices
for the sake of reader's convenience.
See \cite{nikulin79int} for details.
Let $L$ be a non-degenerate even lattice.
We denote by $L^\vee=\OP{Hom}(L,\Z)$ the dual of $L$,
which we consider as a submodule of $L\otimes\Q$ in a natural way.
We extend the bilinear form $\BF{~}{~}$ on $L$
to that on $L^\vee$ linearly.
The discriminant group $A(L)$
 and the discriminant form $q(L)$ of $L$
are defined by
\begin{equation}
A(L)=L^\vee/L,\quad
q(L):A(L)\rightarrow\Q/2\Z;~
x \bmod L \mapsto \BF{x}{x} \bmod 2\Z,
\end{equation}
respectively.
Let $l(A(L))$ denote the minimum number of generators of $A(L)$.
We have $\left| A(L) \right|=\left| \disc(L) \right|$ and
$q(L)$ is a non-degenerate quadratic form on a finite abelian group $A(L)$.
For a prime number $p$,
 let $A(L)_p$ and $q(L)_p$ denote
the $p$-components of $A(L)$ and $q(L)$, respectively.
The quadratic form $q(L)_p$ is considered as
$\Q_p/\Z_p$-valued
(resp.\ $\Q_2/2\Z_2$-valued)
for odd $p$ (resp.\ $p=2$).
Let $M$ be a module such that $L\subset M\subset L^\vee$.
We say $M$ (equipped with $\BF{~}{~}$) is an overlattice of $L$
if $\BF{~}{~}$ takes integer values on $M$.
Any lattice which includes $L$ as a sublattice of finite index
is considered as an overlattice of $L$.
For $\varphi\in \OP{O}(L)$,
we denote by $\bar{\varphi}$ the image of $\varphi$
under the natural homomorphism $\OP{O}(L)\rightarrow \OP{O}(A(L),q(L))$.

\begin{prop}[\cite{nikulin79int}] \label{PROPoverlattice}
Let $L$ be a non-degenerate even lattice.
For a module $M$ such that $L\subset M\subset L^\vee$,
$M$ is an even overlattice of $L$ if and only if the image of $M$ in $A(L)$
is an isotropic subgroup,
i.e., the restriction of $q(L)$ to $M/L$ is zero.
In other words, there exists a natural one-to-one correspondence between
the set of even overlattices of $L$
 and the set of isotropic subgroups of $A(L)$.
\end{prop}

\begin{prop} \label{PROPextension}
Let $\Gamma$ be a non-degenerate lattice
and $L$ be a non-degenerate sublattice of $\Gamma$
(which may be not primitive).
Then, for any $\varphi\in \OP{O}(L)$
whose action on $L^\vee/L$ is trivial,
there exists a unique extension $\Phi\in \OP{O}(\Gamma)$ of $\varphi$
such that $\Phi$ is the identity on $(L)^\bot_\Gamma$.
\end{prop}

\begin{prop}[\cite{nikulin79int}] \label{PROPembedding}
Let $K$ and $L$ be non-degenerate even lattices.
Then, there exists
a primitive embedding of $K$ into an even unimodular lattice $\Gamma$
such that $(\Gamma)^\bot_K\cong L$,
if and only if $(A(K),q(K))\cong (A(L),-q(L))$.
More precisely, any such $\Gamma$ is of the form $\Gamma_\gamma
\subset K^\vee \oplus L^\vee$
for some isomorphism $\gamma:(A(K),q(K))
{\displaystyle\mathop{\rightarrow}^{\sim}}(A(L),-q(L))$,
where $\Gamma_\gamma$ is the lattice corresponding to the isotropic subgroup
$\{x\oplus \gamma(x) \in A(K)\oplus A(L)\bigm|x\in A(K)\}\subset A(K)\oplus A(L)$.
For $\gamma,\gamma':(A(K),q(K))
{\displaystyle\mathop{\rightarrow}^{\sim}}(A(L),-q(L))$,
$\varphi\in \OP{O}(K)$ and $\psi\in \OP{O}(L)$,
the isomorphism $\varphi\oplus\psi$ is extended to
an isomorphism
$\Gamma_\gamma{\displaystyle\mathop{\rightarrow}^{\sim}}\Gamma_{\gamma'}$
if and only if
$\gamma'\circ\bar{\varphi}=\bar{\psi}\circ\gamma$ in $\OP{O}(A(L),q(L))$.
\end{prop}

\begin{thm}[\cite{nikulin79int}] \label{THMunimodularemb}
There exists a primitive embedding
of a non-degenerate even lattice $L$ with signature $(t_{(+)},t_{(-)})$
into an even unimodular lattice with signature $(l_{(+)},l_{(-)})$
if the following conditions are satisfied:
\begin{itemize2}
\item[(i)]
$l_{(+)}-l_{(-)}\equiv 0 \bmod 8$;
\item[(ii)]
$l_{(+)}-t_{(+)}\geq 0,~ l_{(-)}-t_{(-)}\geq 0,~
l_{(+)}+l_{(-)}-t_{(+)}-t_{(-)}>l(A(L)).$
\end{itemize2}
\end{thm}

We need the following classification theorem for quadratic forms
on finite abelian groups.

\begin{thm}[\cite{nikulin79int}] \label{THMfinquad}
Let $p$ be an odd prime number.
Let $A$ and $Q:A\rightarrow\Q_p/\Z_p$ be an abelian $p$-group and
a non-degenerate quadratic form on $A$, respectively.
Then, $Q$ is decomposed into the following form:
$$Q\cong{\textstyle\bigoplus}_{k\geq 1}
\left( q(p^k)^{\oplus n_k} \oplus q'(p^k)^{\oplus m_k} \right),
\quad m_k=0 \rm{\ or\ } 1,$$
where $q(p^k)=\langle 1/p^k \rangle$
and $q'(p^k)=\langle \theta/p^k \rangle$
with non-square $\theta\in\Z_p^\times$
are quadratic forms on $\Z/p^k\Z$.
The numbers $n_k,m_k$ are independent of the choice of decompositions.
\end{thm}

\section{Symplectic actions of $\SYMG$ on $K3$ surfaces}
\label{SECTs5}
In this section,
 we study symplectic actions of $\SYMG$
 on $K3$ surfaces,
 using Kond\=o's method~\cite{kondo98,kondo99}.
The goal is to classify $K3$ surfaces with
 a faithful action of $\SYMG\times C_2$ such that
 the action of $\SYMG$ is symplectic
 (see Theorem~\ref{THMs5x2}).
The result will be applied to determine
 the transcendental lattices of
 special fibers of the family
 $\{ \X_\T \} _{\T\in\P^1}$ of $K3$ surfaces
 defined in the following section.

Let $X$ be a $K3$ surface with
 a faithful symplectic action of a finite group $G$.
By the Mukai's classification of finite groups
 of symplectic actions
 on $K3$ surfaces~\cite{mukai88},
 $G$ is contained in one of the eleven maximal groups:
 $G=M_{20}$, $L_2(7)$, $\MF{A}_6$,
 $F_{384}$, $\SYMG$, etc.
Kond\=o~\cite{kondo98} gave a simplified proof of
 the Mukai's classification in a lattice theoretic way.
The method used in \cite{kondo98} is to embed the coinvariant
 lattice $H(X)_G$
 into a negative definite Niemeier lattice
 (i.e.\ a negative definite even unimodular lattice
 of rank $24$)
 which is not the Leech lattice
 and to describe $G$ as a subgroup of $\OP{O}(N)$.
Using this method,
 Kond\=o also studied a maximal case $G=M_{20}$,
 where $M_{20}$ is the Mathieu group of degree $20$.
He showed that there exists a unique $K3$ surface with
 a faithful action of $M_{20}.C_4$
 (i.e.\ an extension of $M_{20}$ by $C_4$)
 up to isomorphism.
(Since $M_{20}$ is a perfect group,
 the action of $M_{20}$
 is automatically symplectic.)
For similar results
 in the cases
 $G=L_2(7)$, $\MF{A}_6$, $F_{384}$,
 see \cite{oguisozhang02,keumoguisozhang05,oguiso05}.
In this section,
 we consider the case $G=\SYMG$
 and show a similar result.

By Theorem~\ref{THMsymplectic},
 the study on symplectic actions
 of $\SYMG$ on $K3$ surfaces $X$
 is reduced to
 the study on certain actions
 of $\SYMG$
 on the $K3$ lattice $\Lambda$.
Moreover, using Kond\=o's method,
 this is reduced to
 the study on
 certain Niemeier lattices $N$
 with actions of $\SYMG$.
In Subection~\ref{SUBSECTsympA},
 we study such Niemeier lattices $N$
 with actions of $\SYMG$.
We show that the coinvariant lattice
 $N_\SYMG$ is unique up to isomorphism
 in Subsection~\ref{SUBSECTsympB}.
We denote it by $S_{(\SYMG)}$.
We classify actions of $\SYMG$
 on $\Lambda$ induced by
 symplectic actions of $\SYMG$
 on $K3$ surfaces
 in Subsection~\ref{SUBSECTsympC}.
Finally,
 in Subsection~\ref{SUBSECTsympD},
 we classify $K3$ surfaces
 with faithful actions
 of $\SYMG\times C_2$ such that
 the actions of $\SYMG$ are symplectic
 (Theorem~\ref{THMs5x2}).

\subsection{Certain actions of $\SYMG$ on  Niemeier lattices}
\label{SUBSECTsympA}
In this subsection,
 we study certain actions of $\SYMG$ on
 Niemeier lattices $N$
 which are not the Leech lattice.
We will see later that
 the coinvariant lattice $N_\SYMG$ (denoted by $S_{(\SYMG)}$)
 is isomorphic to the coinvariant
 lattice $H(X)_\SYMG$
 for any $K3$ surface $X$ with
 a symplectic action of $\SYMG$
 (see Proposition~\ref{PROPcharacterization}).

We recall some facts on Niemeier lattices.
See \cite{SP} for details.
A definite even unimodular lattice of rank $24$ is
 called a Niemeier lattice.
Let $N$ be a negative definite Niemeier lattice.
A root of $N$ is defined
 as a vector $x\in N$ such that $\BF{x}{x}=-2$.
The isomorphism class of $N$ depends only on
 the $A$-$D$-$E$-type of the root system of $N$.
If $N$ contains no root,
 the lattice $N$ is called the Leech lattice.
Assume that $N$ is not the Leech lattice.
Then, the root system of $N$ is of rank $24$.
There exists a fundamental root system $R$ of $N$,
 i.e.\ a subset $R\subset N$ such that
 $R$ is a basis of $N\otimes\R$ and
 any root is of the form
 $\pm\sum_{r\in R} a_r r$
 with $ a_r \in\Z_{\geq 0}$.
Let $W(N)$ be the Weyl group of $N$,
 i.e.\ the subgroup of $\OP{O}(N)$
 generated by reflections
 by roots of $N$.
The Weyl chamber
 $\{ x\in N\otimes\R\bigm|
 \BF{x}{r}\leq 0~(\forall r\in R) \}$
 corresponding to $R$
 is a fundamental domain of
 the action of $W(N)$ on $N\otimes\R$.
The action of $W(N)$
 on the set of fundamental root systems
 is free and transitive.

In the rest of the subsection,
 we prove the following Lemma~\ref{LEMdiscroot}.
This is a preparation
 for Proposition~\ref{PROPunique6d4},
 where we will see that
 the coinvariant lattice $N_\SYMG$
 in Lemma~\ref{LEMdiscroot} is
 unique up to isomorphism.
In fact, any $N_{\SYMG}$ is reduced to the case (III)
 in Lemma~\ref{LEMdiscroot}
 (see the proof of Proposition~\ref{PROPunique6d4}).

\begin{lem} \label{LEMdiscroot}
Let $N$ be a negative definite Niemeier lattice
 which is not the Leech lattice.
Assume that $\SYMG$ acts on $N$
 with the character $\chi$ satisfying
 the following conditions:
\begin{itemize2}
\item[(i)]
 $\BF{x}{x}\neq -2$ for all $x\in N_\SYMG$;
\item[(ii)]
 $\chi(g)=24,8,6,4,4,2\,$
 if $\,\ord(g)=1,2,3,4,5,6,$ respectively.
\end{itemize2}
Then, we have 
\begin{equation} \label{rankdisc}
 \rank N_\SYMG=19\text{,}\quad
 \operatorname{disc}(N_\SYMG)=-300\text{.}
\end{equation}
Moreover, there exists a fundamental root system $R$
 stable under the action of $\SYMG$.
The $A$\text{\rm -}$D$\text{\rm -}$E$\text{\rm -}type
 of $R$ and
 the type of the orbit decomposition
 (i.e.\ the numbers of elements in the orbits)
 of $R$ under the action of $\SYMG$
 are one of the following cases:
\begin{equation} \label{adetype}
 \begin{array}{c|c|c}
 \text{\rm case}
  & A\text{\rm -}D\text{\rm -}E\text{\rm -type}
  & \text{\rm orbit decomposition} \\
  \hline
 \text{(I)}   & 24A_1 & [1,2,5,6,10] \\
 \text{(II)}  & 24A_1 & [1,1,2,5,15] \\
 \text{(III)} &  6D_4 & [1,1,2,5,15]
 \end{array}.
\end{equation}
Let $R_1,\ldots,R_5$ be the orbits
 under the action of $\SYMG$,
 where we assume that $\sharp R_1\leq\cdots\leq\sharp R_5$.
We set $\nu^{}_S=\sum_{r\in S} r \,$
 for $S\subset R$.
In each case in (\ref{adetype}), we have the following.

The case (I).
The invariant lattice $N^\SYMG$ is generated by
 the following vectors:
\begin{equation} \label{genI_}
 \begin{array}{c}
 \displaystyle
 \nu^{}_{R_1},\ldots,\nu^{}_{R_5},~
 \frac{1}{2}(\nu^{}_{R_1}+\nu^{}_{R_2}+\nu^{}_{R_3}),~
 \frac{1}{2}(\nu^{}_{R_2}+\nu^{}_{R_4}),~ \vspace{1ex} \\
 \displaystyle
 \frac{1}{2}(\nu^{}_{R_1}+\nu^{}_{R_2}+
 \nu^{}_{R_3}+\nu^{}_{R_4}+\nu^{}_{R_5}).
 \end{array}
\end{equation}

The case (II).
The invariant lattice $N^\SYMG$ is generated by
 the following vectors:
\begin{equation} \label{genII_}
 \nu^{}_{R_1},\ldots,\nu^{}_{R_5},~
 \frac{1}{2}
 (\nu^{}_{R_2}+\nu^{}_{R_3}+\nu^{}_{R_4}),~
 \frac{1}{2}(\nu^{}_{R_1}+\nu^{}_{R_2}
 +\nu^{}_{R_3}+\nu^{}_{R_4}+\nu^{}_{R_5}).
\end{equation}

The case (III).
There exists a numbering
 $R=\{
 r^{(j)}_k
 \} ^{1\leq j\leq 6}_{0\leq k\leq 3}$
 satisfying the following conditions:
\begin{itemize2}
\item[(a)]
 $\BF{r^{(j)}_k}{r^{(j')}_{k'}}=0$ if $j\neq j'$,
 $\BF{r^{(j)}_0}{r^{(j)}_k}=1$ if $k\in \{ 1,2,3 \}$,
 and $\BF{r^{(j)}_1}{r^{(j)}_2}
 =\BF{r^{(j)}_1}{r^{(j)}_3}
 =\BF{r^{(j)}_2}{r^{(j)}_3}=0$;
\item[(b)]
 the orbit decomposition of $R$
 under the action of $\SYMG$ is
\begin{equation} \label{orbitIII}
 \{ r^{(1)}_{0} \}
 \sqcup \{ r^{(1)}_{1} \}
 \sqcup \{ r^{(1)}_{2},r^{(1)}_{3} \}
 \sqcup \{ r^{(j)}_{0} \}^{2\leq j \leq 6}
 \sqcup \{ r^{(j)}_{k}
 \}^{2\leq j \leq 6}_{1\leq k \leq 3} \text{;}
\end{equation}
\item[(c)]
 the lattice $N$ is generated by $R$ and
 the following vectors:
\setbox0=\hbox{$v^{(3)}_3$}
\begin{equation} \label{codeIII}
 \begin{array}
  {c@{~}c@{~}c@{~}c@{~}c@{~}c@{~}%
   c@{~}c@{~}c@{~}c@{~}c@{~}c}
 v^{(1)}_1 &+& v^{(2)}_1 &+& v^{(3)}_1 &+&
 v^{(4)}_1 &+& v^{(5)}_1 &+& v^{(6)}_1, \\[5pt]
 v^{(1)}_2 &+& v^{(2)}_2 &+& v^{(3)}_2 &+&
 v^{(4)}_2 &+& v^{(5)}_2 &+& v^{(6)}_2, \\[5pt]
           & &           & & v^{(3)}_2 &+&
 v^{(4)}_3 &+& v^{(5)}_3 &+& v^{(6)}_2, \\[5pt]
           & & v^{(2)}_2 &+& v^{(3)}_3 &+&
 v^{(4)}_3 &+& v^{(5)}_2 & & \hspace{\wd0}, \\[5pt]
           & & v^{(2)}_3 &+& v^{(3)}_3 &+&
 v^{(4)}_2 & &           &+& v^{(6)}_2, \\[5pt]
           & & v^{(2)}_3 &+& v^{(3)}_2 & &
           &+& v^{(5)}_2 &+& v^{(6)}_3,
 \end{array}
\end{equation}
 where we set
$$v^{(j)}_k
 =\frac{1}{2}
 (r^{(j)}_1+r^{(j)}_2+r^{(j)}_3-r^{(j)}_k)
 \text{;}$$
 \item[(d)]
 the coinvariant lattice $N^\SYMG$ is generated by
 the following vectors:
\begin{equation} \label{genIII_}
 \nu^{}_{R_1},\ldots,\nu^{}_{R_5}
 \text{.}
\end{equation}
\end{itemize2}
\end{lem}

\begin{rem} \label{REMexistence_s5}
We will see later that
 there exists an action of $\SYMG$
 on a Niemeier lattice
 satisfying the assumptions in
 Lemma~\ref{LEMdiscroot}
 (see the proof of
 Proposition~\ref{PROPinvariantlattice}).
\end{rem}

We give the character table of $\SYMG$
 in Table~\ref{TBLcharacters} (cf.\ \cite{atlas}).
In the table, $\chi_1,\chi'_1,\ldots$ are characters of
 irreducible representations,
 and $nA,nB$ are conjugacy classes of order $n$.
The third row in the table indicates
 the number of elements in each conjugacy class.
The irreducible decomposition of $\chi$ is
\begin{equation} \label{character}
 \chi=5\chi_1+\chi'_1+2\chi_4+\chi_5+\chi'_5.
\end{equation}
In particular, we have $\rank N_\SYMG=19$.
By the condition (i),
 there exists a vector
 $w\in N^\SYMG\otimes\R$ such that
 $\BF{w}{x}\neq 0$ for any root $x$ of $N$.
Let $R$ be the fundamental root system of $N$
 corresponding
 to the Weyl chamber containing $w$.
Then, the action of $\SYMG$ preserves $R$.
We fix such $R$.
It follows by (\ref{character}) that
 $R$ is decomposed into $5$ orbits $R_1,\ldots,R_5$
 under the action of $\SYMG$,
 where we assume that
 $\sharp R_1\leq\cdots\leq \sharp R_5$.

\begin{table}[tbh]
\caption{}
\label{TBLcharacters}
\begin{center}
\begin{tabular}{c|cccccccc}
\hline
type      &$(1)$&$(2)^2$&$(3)$&$(5)$&$(2)$&$(4)$&$(2)(3)$ \\
 \hline 
notation  &$1A $&$2A $&$3A $&$5A $&$2B $&$4A $&$6A $ \\
 \hline
$\sharp$  &$ 1 $&$15 $&$20 $&$24 $&$10 $&$30 $&$20 $ \\
 \hline
$\chi_1$  &$ 1 $&$ 1 $&$ 1 $&$ 1 $&$ 1 $&$ 1 $&$ 1 $ \\
$\chi'_1$ &$ 1 $&$ 1 $&$ 1 $&$ 1 $&$-1 $&$-1 $&$-1 $ \\
$\chi_4$  &$ 4 $&$ 0 $&$ 1 $&$-1 $&$ 2 $&$ 0 $&$-1 $ \\
$\chi'_4$ &$ 4 $&$ 0 $&$ 1 $&$-1 $&$-2 $&$ 0 $&$ 1 $ \\
$\chi_5$  &$ 5 $&$ 1 $&$-1 $&$ 0 $&$ 1 $&$-1 $&$ 1 $ \\
$\chi'_5$ &$ 5 $&$ 1 $&$-1 $&$ 0 $&$-1 $&$ 1 $&$-1 $ \\
$\chi_6$  &$ 6 $&$-2 $&$ 0 $&$ 1 $&$ 0 $&$ 0 $&$ 0 $ \\
 \hline
\end{tabular}
\end{center}
\end{table}

\begin{lem} \label{LEMorbit}
Let $\psi_i$ be the character
 of the permutation representation
 corresponding to $R_i$
 ($i=1,\ldots,5$).
Then, we have
 $$\psi_i=\chi_1+a\chi'_1+b\chi_4+c\chi_5+d\chi'_5,$$
 where $(a,b,c,d)$ is equal to one of the following:
\begin{equation} \label{orbitabcd}
 (0,0,0,0),(1,0,0,0),(0,1,0,0),(0,0,0,1),
 (0,1,1,0),(0,1,1,1). \tag{$\star$}
\end{equation}
\end{lem}
\begin{proof}
By (\ref{character}), we find that
 $\psi_i$ is of the form
 $\chi_1+a\chi'_1+b\chi_4+c\chi_5+d\chi'_5$
 for some non-negative integers $a,b,c,d$
 with $a,c,d\leq 1$ and $b\leq 2$.
(Since $R_i$ is a orbit under the action of $\SYMG$,
 $\chi_1$ appears exactly once in $\psi_i$.)
Since $\psi_i$ is the character
 of a permutation representation,
 we have
\begin{align*}
 \psi_i(3A)&=1-a+b-c-d\geq 0, \\
 \psi_i(4A)&=1-a-c+d\geq 0, \\
 \psi_i(6A)&=1-a-b+c-d\geq 0.
\end{align*}
Assume that $(a,b,c,d)$ is not one of (\ref{orbitabcd}).
Then the above inequalities implies
 that $(a,b,c,d)=(0,0,1,0)$, i.e., $\psi_i=\chi_1+\chi_5$.
Since $\psi_i(6A)=2$ and $(6A)^2=3A$,
 any element in $3A$ fixes
 at least two elements in $R_i$.
This is a contradiction to $\psi_i(3A)=0$.
\end{proof}

\begin{lem} \label{LEMorbittype}
The $A$\text{\rm -}$D$\text{\rm -}$E$\text{\rm -}type
 of $R$ and
 the type of the orbit decomposition of $R$
 under the action of $\SYMG$ are
 one of the cases in (\ref{adetype}).
\end{lem}

\begin{proof}
By (\ref{character}) and Lemma~\ref{LEMorbit}, we find that
the decomposition of $\chi$
corresponding to the orbit decomposition of $R$
must be either
$$\chi=(\chi_1)+(\chi_1+\chi'_1)+(\chi_1+\chi_4)
+(\chi_1+\chi'_5)+(\chi_1+\chi_4+\chi_5)$$
or
$$\chi=(\chi_1)+(\chi_1)+(\chi_1+\chi'_1)+(\chi_1+\chi_4)
+(\chi_1+\chi_4+\chi_5+\chi'_5).$$
Therefore,
 the type of the orbit decomposition of $R$ is either
 $[1,2,5,6,10]$ or $[1,1,2,5,15]$.

Since ${\SYMG}$ preserves the Weyl chamber
 corresponding to $R$,
 which is a fundamental domain
 of the Weyl group $W(N)$,
 the group $\SYMG$ is considered
 as a subgroup of $G(N):=\OP{O}(N)/W(N)$.
This implies that $R=24A_1,12A_2,6A_4$ or $6D_4$,
 because $\left| {\SYMG} \right|=120$
 does not divide $\left| G(N) \right|$ in other cases
 (see Table~16.1 of \cite{SP}).
Assume that $R=12A_2$ or $6A_4$.
Then, for each odd $k$,
the number of the $R_i$ with $\sharp R_i=k$ is even.
Since the type of the orbit decomposition of $R$
 is either
 $[1,2,5,6,10]$ or $[1,1,2,5,15]$,
 this is a contradiction.
Therefore these cases are excluded.
Hence we have $R=24A_1$ or $6D_4$.
We consider the case $R=6D_4$.
Since $\SYMG$ fixes at least one element in $R$,
 $\SYMG$ preserves at least one component
 ($=D_4$) of $R$.
Hence we have $\sharp(\bigcup_{i\in I}R_i)=4$
for some $I\subset\{1,\ldots,5\}$.
This implies that
 the type of the orbit decomposition of $R$
 is $[1,1,2,5,15]$.
\end{proof}

\begin{lem} \label{LEMdisc300}
$\disc(N_\SYMG)=-300$.
\end{lem}
\begin{proof}
Since we have $\disc(N_\SYMG)=\disc(N^\SYMG)$,
 it is sufficient to show that
 $\disc(N^\SYMG)=-300$ in each case
 in (\ref{adetype}).

The case (I).
The orbits $R_1,\ldots,R_5$ consist
 of $1,2,5,6,10$ elements, respectively.
We set
 $\MC{G}
 =\{ S\subset R \bigm|
 \nu^{}_S/2 \in N \}$.
Then, $\MC{G}$ is identified with
 $N/\langle R \rangle$,
 where
 $\langle R \rangle=\bigoplus_{r\in R}\Z r$
 is
 the sublattice of $N$ generated by $R$.
Also, $\MC{G}$ is identified
 with the (extended) binary Golay code
 (see \cite{SP} for details).
In particular,
 $\MC{G}_8:=\{ O\in \MC{G}\bigm| \sharp O=8 \}$ is
 the Steiner system $S(5,8,24)$, i.e.,
 for any subset
 $S\subset R$ with $\sharp S=5$,
 there is one and only one $O\in\MC{G}_8$
 with $S\subset O$ (Steiner property).
Therefore, there exists unique $O'\in\MC{G}_8$
 with $R_3\subset O'$.
Since $R_3$ is stable under the action of $\SYMG$,
 so is $O'$.
Hence we have $O'=R_1\sqcup R_2 \sqcup R_3$.
Next, we shall show
 that $R_2\sqcup R_4\in\MC{G}_8$.
By Steiner property,
 there exists $O''\in\MC{G}_8$ such that
 $\sharp(O''\cap R_4)\geq 5$.
Assume that $R_4\not\subset O''$.
Then Steiner property implies that
$H:=\{g\in \SYMG \bigm| g\cdot O''=O''\}$
is a subgroup of $\SYMG$ of index 6.
We fix $\alpha,\beta\in \SYMG$
such that $\alpha\in H$, $\operatorname{ord}(\alpha)=5$
and $\operatorname{ord}(\beta)=3$.
Since we have $\chi(\alpha)=4$,
 $\alpha$ fixes four elements in $R$.
Hence, the orbit decomposition of $O''$
under the action of $\langle \alpha \rangle$
is of type $[1,1,1,5]$.
Since $\alpha$ fixes $R_2$ point-wise,
$O''\cap R_2$ contains at least one element fixed by $\alpha$.
In particular, $O''\cap R_2\neq\varnothing$.
Therefore, we have
$$\sharp((\beta\cdot O'')\cap O'')\geq
 \sharp((\beta\cdot O'')\cap O''\cap R_2)
 +\sharp((\beta\cdot O'')\cap O'' \cap R_4)
 \geq 1+4=5.$$
By Steiner property, we have $\beta\in H$.
This is a contradiction to $[\SYMG:H]=6$,
i.e., $\left| H \right|=20$.
Hence we have $R_4\subset O''$.
Since $R_4$ is stable under the action of $\SYMG$,
 so is $O''$.
This implies that $O''=R_2\sqcup R_4$.
Recall that
 any $S\in \MC{G}$ consists of
 $0,8,12,16,24$ elements (cf.\ \cite{SP}).
As a consequence,
 for the union $S$ of some of $R_i$'s,
 $S$ is an element in $\MC{G}$
 if and only if $\sharp S=0,8,12,16,24$.
Hence, we find that
\begin{equation*} \label{genI}
 \nu^{}_{R_1},\ldots,\nu^{}_{R_5},~
 \frac{1}{2}(\nu^{}_{R_1}+\nu^{}_{R_2}+\nu^{}_{R_3}),~
 \frac{1}{2}(\nu^{}_{R_2}+\nu^{}_{R_4}),~
 \frac{1}{2}(\nu^{}_{R_1}+\nu^{}_{R_2}+
 \nu^{}_{R_3}+\nu^{}_{R_4}+\nu^{}_{R_5})
\end{equation*}
 generate $N^\SYMG$.
Therefore,
 $\operatorname{disc}(N^\SYMG)=
 (-2)^5\cdot 1\cdot 2\cdot 5\cdot 6\cdot 10/2^6
 =-300$.

The case (II).
The orbits $R_1,\ldots,R_5$ consist
of $1,1,2,5,15$ elements, respectively.
Moreover, by the same argument as above,
 we may assume that
 $R_2\sqcup R_3 \sqcup R_4\in\MC{G}$
 after interchanging $R_1$ and $R_2$
 if necessary.
Then, $N^\SYMG$ is generated by
\begin{equation*} \label{genII}
 \nu^{}_{R_1},\ldots,\nu^{}_{R_5},~
 \frac{1}{2}
 (\nu^{}_{R_2}+\nu^{}_{R_3}+\nu^{}_{R_4}),~
 \frac{1}{2}(\nu^{}_{R_1}+\nu^{}_{R_2}
 +\nu^{}_{R_3}+\nu^{}_{R_4}+\nu^{}_{R_5}).
\end{equation*}
Hence, we have
 $\operatorname{disc}(N^{\SYMG})=
 (-2)^5\cdot 1\cdot 1\cdot 2\cdot 5\cdot 15/2^4
 =-300$.

The case (III).
By Table~16.1 of \cite{SP},
 the assertions (a) and (c) holds
 for some numbering $R=\{r^{(j)}_k\}$.
Let $G(N,R)$ (resp.\ $G_1(N,R)$) be
 the subgroup of $\OP{O}(N)$ of elements
 preserving $R$
 (resp.\ each component of $R$).
We set $G_2(N,R)=G(N,R)/G_1(N,R)$.
Again, by Table~16.1 of \cite{SP},
 we find that $G_1(N,R)\cong\Z/3\Z$
 and $G_2(N,R)\cong\MF{S}_6$.
A generator of $G_1(N,R)$ acts
on each component of $R$ non-trivially.
The group $G_2(N,R)$ is
the group of all permutations of the set of components of $R$.
Since $\SYMG$ fixes at least one vertex of $R$,
$\SYMG$ preserves at least one component of $R$.
We may assume that $\SYMG$ preserves the component
$\{r^{(1)}_{k}\}_{0\leq k \leq3}$,
after replacing ${\SYMG}\subset \OP{O}(N)$
with $\alpha {\SYMG} \alpha^{-1}$
for suitable $\alpha\in G(N,R)$ if necessary.
Since the type
 of the orbit decomposition of $R$
 is $[1,1,2,5,15]$,
 that of
 $\{r^{(1)}_{k}\}$
 is $[1,1,2]$.
We may assume that $\SYMG$ fixes $r^{(1)}_{1}$,
 after replacing $\SYMG$
 with $\beta {\SYMG} \beta^{-1}$
 for suitable $\beta\in G_1(N,R)$ if necessary.
Then, the orbit decomposition of $R$ is as follows:
\begin{equation}
 \{ r^{(1)}_{0} \}
 \sqcup \{ r^{(1)}_{1} \}
 \sqcup \{ r^{(1)}_{2},r^{(1)}_{3} \}
 \sqcup \{ r^{(j)}_{0} \}^{2\leq j \leq 6}
 \sqcup \{ r^{(j)}_{k}
 \}^{2\leq j \leq 6}_{1\leq k \leq 3}.
\end{equation}
Since $N$ is generated by $R$ and
 the vectors in (\ref{codeIII}),
 we can check that
 $\nu^{}_{R_1},\ldots,\nu^{}_{R_5}$
 form a basis of $N^\SYMG$.
The Grammian matrix of $N^\SYMG$ for this basis is
\begin{equation} \label{6d4gram}
 \begin{pmatrix}
 -2 & 1 & 2 & 0 & 0 \\
 1  &-2 & 0 & 0 & 0 \\
 2  & 0 &-4 & 0 & 0 \\
 0  & 0 & 0 &-30&15 \\
 0  & 0 & 0 &15 &-10
 \end{pmatrix}.
\end{equation}
Hence we have $\operatorname{disc}(N^\SYMG)=-300$.
\end{proof}

By the proof of Lemma~\ref{LEMdisc300},
 we can check the latter part
 of Lemma~\ref{LEMdiscroot}.
Now we have proved the lemma.

\subsection{Properties on the coinvariant lattice $S_{(\SYMG)}$}
\label{SUBSECTsympB}
In this subsection,
 we study the lattice $S_{(\SYMG)}$,
 which is defined as the coinvariant lattice
 $N_\SYMG$ in Lemma~\ref{LEMdiscroot}
 (see Definition~\ref{DEFNs5}).
In the following subsection,
 we will find that $S_{(\SYMG)}$ is isomorphic to
 the coinvariant lattice $H(X)_\SYMG$ for
 any $K3$ surface $X$
 with a faithful symplectic action of $\SYMG$.
Using the result in the previous subsection,
 we prove the uniqueness of
 $N_\SYMG$ in Lemma~\ref{LEMdiscroot}.

\begin{prop} \label{PROPunique6d4}
Under the assumptions
 in Lemma~\ref{LEMdiscroot},
 the coinvariant lattice $N_\SYMG$
 is unique up to isomorphism.
Moreover,
 The discriminant group and the discriminant form
 of $N_\SYMG$ are as follows:
\begin{equation} \label{disc_form_coinv}
 \begin{array}{c}
 A(N_\SYMG) \cong
  \Z/4\Z \oplus \Z/3\Z
  \oplus \Z/5\Z \oplus \Z/5\Z \\[5pt]
 q(N_\SYMG) \cong
  \LF{3/4} \oplus \LF{1/3}
  \oplus \LF{1/5} \oplus \LF{2/5}.
 \end{array}
\end{equation}
\end{prop}

\begin{defn} \label{DEFNs5}
We denote by $S_{(\SYMG)}$
 the coinvariant lattice $N_\SYMG$ determined
 in Proposition~\ref{PROPunique6d4}.
\end{defn}

\begin{rem}
The lattice $S_{(\SYMG)}$ actually exists
 (see the proof
 of Proposition~\ref{PROPinvariantlattice}).
We denote by $S_{(\SYMG)}$
 only the lattice structure of $N_\SYMG$.
However, the action of $\SYMG$ is
 characterized
 by Lemma~\ref{LEMsurjectivityS5} below
 (up to isomorphism).
\end{rem}

\begin{proof}
[Proof of Proposition~\ref{PROPunique6d4}]
Let $M$ be the negative definite
 root lattice of type $A_2$
 and $A$ the discriminant group of
 $N_\SYMG\oplus M$.
We apply Theorem~\ref{THMunimodularemb}
 to embed $N_\SYMG\oplus M$
 into another negative definite
 Niemeier lattice.
By Lemma~\ref{LEMdiscroot},
 we have $|A|=300\cdot 3=2^2\cdot 3^2\cdot 5^2$.
Hence we have $l(A)\leq 2$.
Since $\rank(N_\SYMG\oplus M)+l(A)\leq 23<24$,
 it follows by Theorem~\ref{THMunimodularemb}
 that
 $N_\SYMG\oplus M$ can be embedded primitively
 into a negative definite Niemeier lattice $N'$.
By Proposition~\ref{PROPembedding},
 the action of $\SYMG$ on $A(N_\SYMG)$ is trivial.
By Proposition~\ref{PROPextension},
 the action of $\SYMG$ on $N_\SYMG$ can be
 extended to that on $N'$ such that
 ${N'}_\SYMG=N_\SYMG$.
Since $N'$ contains the root lattice $M$
 of type $A_2$,
 we find that $N'$ is not the Leech lattice.
Hence $N'$ also satisfies the assumptions in
 Lemma~\ref{LEMdiscroot}.
Since the root system of type $24A_1$ can not
 contain that of type $A_2$,
 it follows that $N'$ is of type (III)
 in Lemma~\ref{LEMdiscroot}.
Therefore, ${N'}_\SYMG=N_\SYMG$
 is unique up to isomorphism.
By the Grammian matrix (\ref{6d4gram})
 of $(N')^\SYMG$,
 we can check the latter part of the assertion
 by a direct computation
 because we have
 $(A({N'}_\SYMG),q({N'}_\SYMG))\cong
 (A((N')^\SYMG),-q((N')^\SYMG))$ by
 Proposition~\ref{PROPembedding}.
\end{proof}

\begin{rem} \label{REMoccur}
In the cases (I)--(III)
 in Lemma~\ref{LEMdiscroot},
 we can check by a direct computation that
 the discriminant forms of
 $N^\SYMG$ are isomorphic to each other.
On the other hand, the number of vectors
 $x\in N^\SYMG$ with $x^2=-2$
 in each case
 is as follows:
\begin{center}
 \begin{tabular}{c|ccc}
 case & (I) & (II) & (III) \\ \hline
 $\sharp$ & $2$ & $4$ & $12$
 \end{tabular}\lower1.5ex\hbox{.}
\end{center}
Hence $N^\SYMG$ for the cases (I)--(III)
 are not isomorphic to each other.
By Propositions~\ref{PROPextension},
 \ref{PROPembedding}
 and Lemma~\ref{LEMdiscroot},
 the cases (I)--(III) can actually occur.
\end{rem}

We will consider
 primitive embeddings of $S_{(\SYMG)}$
 into the $K3$ lattice $\Lambda$
 in the following subsection.
By Proposition~\ref{PROPembedding}
 and the following lemma,
 we find that any automorphism
 of $(S_{(\SYMG)})^\bot_\Lambda$
 can be extended to that of $\Lambda$
 (see Lemma~\ref{LEMext_lambda}).

\begin{lem} \label{LEMsurjectivityS5}
We have the natural short exact sequence
\begin{equation}
 1 \longrightarrow \SYMG
 \mathop{\longrightarrow}^{\alpha}
 \OP{O}(S_{(\SYMG)})
 \mathop{\longrightarrow}^{\beta}
 \OP{O}(A(S_{(\SYMG)}),q(S_{(\SYMG)}))
 \longrightarrow 1.
\end{equation}
Moreover, we have the natural decomposition
\begin{equation}
 \OP{O}(S_{(\SYMG)})
 = \SYMG \times
 \OP{O}(A(S_{(\SYMG)}),q(S_{(\SYMG)})).
\end{equation}
\end{lem}
\begin{proof}
First, we shall show that
 $\OP{Im}(\alpha)=\OP{Ker}(\beta)$.
We use the identification $S_{(\SYMG)}=N_\SYMG$,
 where $N$ is of type (III) in
 Lemma~\ref{LEMdiscroot}
 (cf.\ Remark~\ref{REMoccur}).
We use the same notation as
 in the proof of Lemma~\ref{LEMorbittype}
 in the case (III).
Let $G$ be the kernel of $\beta$.
We have $\SYMG\subset G$.
By Proposition~\ref{PROPextension},
 the action of $G$ on $S_{(\SYMG)}$
 can be extended to that on $N$
 such that $N^{G}=N^\SYMG$.
Hence $G$ preserves $R$ and
 the orbit decomposition of $R$
 under the action of $G$ coincides with
 that under the action of $\SYMG$.
The groups $\SYMG$ and $G$ are considered
 as subgroups of $G_2(N,R)\cong\mathfrak{S}_6$,
 which is the symmetry group
 of the set of components of $R$.
Since $\SYMG$ and $G$ preserve
 one component of $R$,
 we have $G=\SYMG$.

Next, we shall prove the surjectivity of $\beta$.
We use the identification
 $S_{(\SYMG)}=N_\SYMG$,
 where $N$ is of type (I)
 in Lemma~\ref{LEMdiscroot}
 (cf.\ Remark~\ref{REMoccur}).
We use the same notation as in the proof of
 Lemma~\ref{LEMdisc300} in the case (I).
We set $A=A(N^\SYMG)$ and $q=q(N^\SYMG)$.
Let $H$ be
 the subgroup of the Weyl group $W(N)$ of $N$
 which consists of elements preserving $N_\SYMG$.
By Proposition~\ref{PROPembedding},
 it is sufficient to show that
 the natural map
 $H\rightarrow\OP{O}(A,q)$
 is surjective.
For any subset $S\subset R$,
 there exists a unique element
 $\varphi^{}_S\in W(N)$
 such that $\varphi(r)=r$ if $r\not\in S$ and
 $\varphi(r)=-r$ if $r\in S$.
(This gives the natural one-to-one correspondence
 between $2^R$ and $W(N)$.)
The element $\varphi^{}_S\in W(N)$
 is an element in $H$
 if and only if $S$ is the union of
 some of $R_i$'s.
Let $e_1,\ldots,e_4$ be the images in $A$ of
$$\frac{1}{4}
 (2\nu^{}_{R_1}+\nu^{}_{R_2}+\nu^{}_{R_4}+\nu^{}_{R_5}),~
 \frac{1}{3}\nu^{}_{R_4},~
 \frac{1}{5}\nu^{}_{R_5},~
 \frac{1}{5}\nu^{}_{R_3}
 \in (N^\SYMG)^\vee,$$
respectively.
We can see that $e_1,\ldots,e_4$ generate $A$.
Under these generators, we have
$A\cong\Z/2\oplus\Z/3
 \oplus\Z/5\oplus\Z/5$.
Also, we find that $e_i\bot e_j$ for $i\neq j$ and
\begin{gather*}
e_1^2\equiv 5/4 \bmod 2\Z;\quad
e_2^2\equiv 2/3,~
e_3^2\equiv 1/5,~
e_4^2\equiv 2/5 \bmod \Z.
\end{gather*}
We can check by a direct computation that
$\OP{O}(A,q)$ consists of all $\gamma$
such that $\gamma(e_i)=\pm e_i$ for all $i$.
In particular,
 we have $\OP{O}(A,q)\cong (\Z/2\Z)^4$.
We can check the surjectivity of
 $H\rightarrow\OP{O}(A,q)$ by a direct computation.
For example, we consider
 the case that $\gamma\in \OP{O}(A,q)$ sends
 $e_1,e_2,e_3,e_4$
 to $-e_1,e_2,e_3,e_4$, respectively.
We put $\varphi=\varphi^{}_{R_2}$.
Then, it is clear that $\bar{\varphi}(e_i)=e_i$
 for $i=2,3,4$.
Since
 $\nu^{}_{R_1},(\nu^{}_{R_4}+\nu^{}_{R_5})/2\in N$,
 we have
\begin{align*}
 \bar{\varphi}(e_1)
 \equiv \frac{1}{4}
 (2\nu^{}_{R_1}-\nu^{}_{R_2}+\nu^{}_{R_4}+\nu^{}_{R_5})
 \equiv \frac{1}{4}
 (-2\nu^{}_{R_1}-\nu^{}_{R_2}-\nu^{}_{R_4}-\nu^{}_{R_5})
 \equiv -e_1.
\end{align*}
Hence we have $\bar{\varphi}=\gamma$.
Other cases are similar.

The latter part of the lemma follows from
 the fact that
 the outer automorphism group and
 the center of $\SYMG$ are trivial.
\end{proof}

We use the following lemma in Subsection~\ref{SUBSECTtransB}
 (see Lemma~\ref{LEMsurjectivity_T}).

\begin{lem} \label{LEMsurjectivityA5}
Let $K=(S_{(\SYMG)})_{\mathfrak{A}_5}$ be
the coinvariant part of $S_{(\SYMG)}$ under
the action of $\mathfrak{A}_5\subset \SYMG$.
Then, we have
 $\rank K=18$
 and $\operatorname{disc}(K)=300$.
Moreover, we have the natural
 short exact sequence
\begin{equation}
 1 \longrightarrow \MF{A}_5
 \mathop{\longrightarrow}^{\alpha}
 \OP{O}(K)
 \mathop{\longrightarrow}^{\beta}
 \OP{O}(A(K),q(K))
 \longrightarrow 1.
\end{equation}
\end{lem}
\begin{proof}
We can prove the lemma by the same argument
 as in the proof
 of Lemma~\ref{LEMdiscroot}
 and \ref{LEMsurjectivityS5}.
We use the identification
 $S_{(\SYMG)}=N_\SYMG$,
 where $N$ is of type (III)
 in Lemma~\ref{LEMdiscroot}.
By (\ref{character}),
 we have $\rank K=18$.
The type of the orbit decomposition
 $R=R'_1\sqcup\cdots\sqcup R'_6$
 under the action of $\mathfrak{A}_5$
 is $[1,1,1,1,5,15]$.
We can check that
 $\nu^{}_{R'_1},\ldots,\nu^{}_{R'_6}$ generate
 $N^{\MF{A}_5}$.
By a direct computation, we have
 $\disc(K)=\disc(N^{\MF{A}_5})=300$ and
\begin{gather*}
A(K) \cong (\Z/2\Z)^2 \oplus
 \Z/3\Z \oplus (\Z/5\Z)^2, \\
q(K) \cong
 \begin{pmatrix}
   1 & 1/2 \\
 1/2 &   1
 \end{pmatrix}
 \oplus \LF{2/3}
 \oplus \LF{1/5} \oplus \LF{2/5}.
\end{gather*}
Hence we have
 $\OP{O}(A(K),q(K))
 \cong \mathfrak{S}_3\times(\Z/2\Z)^3$.
We use $N$ of type (III) to prove
 that $\OP{Im}(\alpha)=\OP{Ker}(\beta)$
 and use $N$ of types (I)--(III) to prove
 that $\beta$ is surjective,
 as in the proof of
 Lemma~\ref{LEMsurjectivityS5}.
We omit the details of the computation.
\end{proof}

\subsection{Certain actions of $\SYMG$ on the $K3$ lattice}
\label{SUBSECTsympC}
In this subsection,
 we classify actions of $\SYMG$
 on the $K3$ lattice $\Lambda$ induced
 by faithful symplectic actions
 on $K3$ surfaces
 (cf.\ Theorem~\ref{THMsymplectic}).
By the following proposition,
 the coinvariant lattice $\Lambda_\SYMG$
 is isomorphic to the lattice
 $S_{(\SYMG)}$ (see Definition~\ref{DEFNs5}).

\begin{prop} \label{PROPinvariantlattice}
There exist exactly two isomorphism classes
 of faithful actions of $\SYMG$ on $\Lambda$
 satisfying the following conditions:
\begin{itemize2}
\item[(i)]
 $\Lambda_\SYMG$ is negative definite;
\item[(ii)]
 $\BF{x}{x}\neq -2$ for all $x\in \Lambda_{\SYMG}$.
\end{itemize2}
Such an action depends only on
 the isomorphism class of $\Lambda^\SYMG$.
Moreover, we have the following:
\begin{itemize2}
\item[(a)]
 $\Lambda_{\SYMG}$ is isomorphic to
 $S_{(\SYMG)}$;
\item[(b)]
 $\Lambda^{\SYMG}$ is isomorphic
 to either of the following lattices:
\begin{equation} \label{invariantlattice}
 K_1=
 \begin{pmatrix}
 4 & 1 & 0 \\
 1 & 4 & 0 \\
 0 & 0 & 20
 \end{pmatrix},\quad
 K_2=
 \begin{pmatrix}
 4 & 2 & 2 \\
 2 & 6 & 1 \\
 2 & 1 & 16
 \end{pmatrix}.
\end{equation}
\end{itemize2}
\end{prop}
\begin{proof}
There exists a $K3$ surface with
 a symplectic action of $\SYMG$
 (see \cite{mukai88,kondo98}).
By Theorem~\ref{THMsymplectic},
 there exists an action of $\SYMG$
 on $\Lambda$ satisfying
 the conditions (i) and (ii).

We shall check the assertions (a) and (b).
We can embed $\Lambda_\SYMG$ primitively
 into a negative definite Niemeier lattice $N$
 which is not the Leech lattice \cite{kondo98}.
By Proposition~\ref{PROPembedding},
 the action of $\SYMG$ on $A(\Lambda_\SYMG)$
 is trivial.
By Proposition~\ref{PROPextension},
 the action of $\SYMG$
 on $\Lambda_\SYMG$ can be extended
 to that on $N$
 such that $N_\SYMG=\Lambda_\SYMG$.
Now we apply Proposition~\ref{PROPunique6d4}.
In fact, the assumptions
 in Lemma~\ref{LEMdiscroot} are satisfied
 by Lemma~\ref{LEMtrace}.
Hence we have $\Lambda_\SYMG\cong S_{(\SYMG)}$.
On the other hand,
 by Proposition~\ref{PROPunique6d4}
 and
 $(A(\RL\Lambda^\SYMG),q(\RL\Lambda^\SYMG))
 \cong
 (A(\RL\Lambda_\SYMG),-q(\RL\Lambda_\SYMG))$,
 we have
\begin{align*}
 A(N^\SYMG) &\cong
  \Z/4\Z \oplus \Z/3\Z
  \oplus \Z/5\Z \oplus \Z/5\Z \\
 q(N^\SYMG) &\cong
  \LF{5/4} \oplus \LF{2/3}
  \oplus \LF{1/5} \oplus \LF{2/5}.
\end{align*}
From the table of primitive definite even
 ternary quadratic forms
 \cite{brandtintrau58,schiemann},
 we can check that
 the positive definite even lattices $K$
 of rank $3$
 such that
 $(A(K),q(K)) \cong
 (A(\RL\Lambda^\SYMG),q(\RL\Lambda^\SYMG))$
 are $K_1$ and $K_2$.
Hence we have
 $\Lambda^\SYMG\cong K_i$ ($i=1,2$).

For $i=1,2$,
 we shall show the existence and uniqueness
 of an action of $\SYMG$ on $\Lambda$
 satisfying the conditions (i), (ii)
 and $\Lambda^\SYMG\cong K_i$.
By Proposition~\ref{PROPembedding}
 and Lemma~\ref{LEMsurjectivityS5},
 there exists a unique primitive embedding
 of $S_{(\SYMG)}$ into $\Lambda$
 such that
 $(S_{(\SYMG)})^\bot_\Lambda\cong K_i$
 up to isomorphism.
By Proposition~\ref{PROPextension},
 the action of $\SYMG$ on $S_{(\SYMG)}$
 can be extended to that on $\Lambda$
 such that $\Lambda_\SYMG=S_{(\SYMG)}$.
The existence of an action follows.
Note that $\SYMG$ is characterized
 as the kernel of the natural map
 $\OP{O}(S_{(\SYMG)})\rightarrow
 \OP{O}(A(S_{(\SYMG)}),q(S_{(\SYMG)}))$
 by Lemma~\ref{LEMsurjectivityS5}.
Since the outer automorphism group
 of $\SYMG$ is trivial,
 the uniqueness of an action
 follows.
\end{proof}

To study non-symplectic extensions
 of $\SYMG$ in $\OP{Aut}(X)$
 in the above setting,
 we need the following lemma
 (see Theorem~\ref{THMs5x2}).

\begin{lem} \label{LEMext_lambda}
Assume that $\SYMG$ acts on $\Lambda$
 satisfying the conditions (i) and (ii)
 in Proposition~\ref{PROPinvariantlattice}.
Then, any automorphism of $\Lambda^\SYMG$
 can be extended to
 that of $\Lambda$ which commutes
 with the action of $\SYMG$
 in a unique way.
\end{lem}
\begin{proof}
Let $\varphi$ be an automorphism of $\Lambda^\SYMG$
 and $\gamma:(A(\Lambda^\SYMG),q(\Lambda^\SYMG))
 \rightarrow (A(\Lambda_\SYMG),-q(\Lambda_\SYMG))$
 a canonical isomorphism.
By Lemma~\ref{LEMsurjectivityS5} and
 Proposition~\ref{PROPinvariantlattice},
 there exists a unique element
 $\psi\in\OP{O}(\Lambda_\SYMG)$
 such that
 $\gamma\circ\bar{\varphi}=\bar{\psi}\circ\gamma$
 and $\psi$ commutes with the action
 of $\SYMG$ on $\Lambda_\SYMG$.
By Proposition~\ref{PROPembedding},
 $\varphi\oplus\psi$ can be extended
 to an automorphism of $\Lambda$
 which commutes with the action of $\SYMG$.
\end{proof}

\subsection{Special $K3$ surfaces with a symplectic action %
 of $\SYMG$}
\label{SUBSECTsympD}
In this subsection,
 we classify $K3$ surfaces
 with a faithful action of $\SYMG\times C_2$
 such that the action of $\SYMG$ is symplectic,
 using the result on
 certain actions of $\SYMG$
 on the $K3$ lattice $\Lambda$
 in the previous section.
Let $X$ be a $K3$ surface
 with a faithful symplectic action of $\SYMG$.
Since $\SYMG$ appears in the Mukai's list
 of maximal finite
 symplectic actions on $K3$ surfaces
 \cite{mukai88},
 such $X$ exists and the action of $\SYMG$ is
 a maximal finite symplectic action.
As a preparation, we characterize a $K3$ surface
 which admits a faithful symplectic action
 of $\SYMG$,
 using the lattice $S_{(\SYMG)}$
 (see Definition~\ref{DEFNs5}).

\begin{prop} \label{PROPcharacterization}
A $K3$ surface $X$ admits
 a faithful symplectic action of $\SYMG$
 if and only if
 the lattice $S_{(\SYMG)}$
 can be embedded primitively into
 the algebraic lattice $S(X)$ of $X$.
For a faithful symplectic action
 of $\SYMG$ on $X$,
 the coinvariant lattice $H(X)_\SYMG$ is
 isomorphic to $S_{(\SYMG)}$
 and
 the invariant lattice $H(X)^\SYMG$ is
 isomorphic to $K_i$ ($i=1,2$)
 in (\ref{invariantlattice}).
\end{prop}
\begin{proof}
By Proposition~\ref{PROPinvariantlattice},
 it is sufficient to show
 that if $S_{(\SYMG)}$ can be embedded
 primitively into $S(X)$,
 $X$ admits a faithful symplectic action.
In fact, by Proposition~\ref{PROPextension},
 the action of $\SYMG$ on $S_{(\SYMG)}$
 can be extended to that on $H(X)$
 such that $H(X)_\SYMG=S_{(\SYMG)}$
 (cf.\ Lemma~\ref{LEMsurjectivityS5}).
By Theorem~\ref{THMsymplectic},
 $X$ admits a faithful symplectic action
 of $\SYMG$.
\end{proof}

In the above setting,
 we consider a finite subgroup
 $\widetilde{G}$
 of $\operatorname{Aut}(X)$
 such that $\SYMG\subset\widetilde{G}$.
By the maximality of the action of $\SYMG$,
 it follows that
 $\SYMG$ is a normal subgroup of $\widetilde{G}$
 and $\widetilde{G}/\SYMG$ acts
 on $H^{2,0}(X)\cong\C$ faithfully.
The index $[\widetilde{G}:\SYMG]$
 is equal to $1$ or $2$ \cite{zhang06}.
We are interested
 in the case $[\widetilde{G}:\SYMG]=2$.
In this case, we have
 $\widetilde{G}\cong \SYMG\times C_2$,
 where $C_2$ acts on $H^{2,0}(X)$ faithfully.

\begin{thm} \label{THMs5x2}
Let $X$ be a $K3$ surface
 with a faithful action
 of a finite group $\widetilde{G}$.
Assume that $\widetilde{G}$ contains $\SYMG$ properly
 and the action of $\SYMG$ is symplectic.
Then, we have $\widetilde{G}\cong\SYMG\times C_2$,
 where the action of $C_2$ on $H^{2,0}(X)$ is faithful.
There exist exactly six such $X$
 with actions of $\widetilde{G}$
 up to isomorphism.
Moreover, such $X$ corresponds one-to-one to
 each case in Table~\ref{TBLtranslattice}.
\end{thm}

\begin{table}[tbh]
\caption{}
\label{TBLtranslattice}
\begin{center}
\begin{tabular}{c|c|c|c}
  \hline
  $H(X)^\SYMG$ & $T(X)$
      & $H(X)^{\widetilde{G}}$ \rule[-1.2ex]{0ex}{4ex}
      & $[H(X)^\SYMG:T(X)\oplus H(X)^{\widetilde{G}}]$ \\ \hline

  & $\begin{pmatrix} 10 &  0 \\  0 & 20 \end{pmatrix}$
      & $\langle  6 \rangle$ & $2$\rule[-2.8ex]{0ex}{7ex} \\
  \smash{$\begin{pmatrix}
            4 & 1 & 0 \\
            1 & 4 & 0 \\
            0 & 0 & 20
  \end{pmatrix}$}
  & $\begin{pmatrix}  6 &  0 \\  0 & 20 \end{pmatrix}$
      & $\langle 10 \rangle$ & $2$\rule[-2.8ex]{0ex}{7ex} \\
  & $\begin{pmatrix}  4 &  1 \\  1 &  4 \end{pmatrix}$
      & $\langle 20 \rangle$ & $1$\rule[-2.8ex]{0ex}{7ex} \\ \hline

  & $\begin{pmatrix} 20 & 10 \\ 10 & 20 \end{pmatrix}$
      & $\langle  4 \rangle$ & $2$\rule[-2.8ex]{0ex}{7ex} \\
  \smash{$\begin{pmatrix}
            4 & 2 & 2 \\
            2 & 6 & 1 \\
            2 & 1 & 16
  \end{pmatrix}$}
  & $\begin{pmatrix}  4 &  2 \\  2 & 16 \end{pmatrix}$
      & $\langle 20 \rangle$ & $2$\rule[-2.8ex]{0ex}{7ex} \\
  & $\begin{pmatrix}  4 &  2 \\  2 &  6 \end{pmatrix}$
      & $\langle 60 \rangle$ & $2$\rule[-2.8ex]{0ex}{7ex} \\ \hline
\end{tabular} \\
\end{center}
\end{table}

\begin{proof}
Let $X$ be such a $K3$ surface.
By Proposition~\ref{PROPcharacterization},
 we have $H(X)_\SYMG\cong S_{(\SYMG)}$.
By Theorem~\ref{THMsymplectic}
 and Lemma~\ref{LEMsurjectivityS5},
 we find that the subgroup of $\widetilde{G}$
 of symplectic actions coincides with $\SYMG$.
In particular, the action of $\widetilde{G}/\SYMG$
 on $H^{2,0}(X)\cong\C$ is faithful.
By Proposition~\ref{PROPcharacterization},
 the invariant lattice $H(X)^\SYMG$ is
 isomorphic to $K_i$ ($i=1,2$).
We can check that
 $\OP{O}(K_1)\cong\OP{O}(K_2)\cong(\Z/2\Z)^3$.
Hence we have $\widetilde{G}/\SYMG\cong \Z/2\Z$.
Since the outer automorphism group
 and the center of $\SYMG$ are trivial,
 we have a canonical decomposition
 $\widetilde{G}=\SYMG\times C_2$,
 where the action of $C_2$
 on $H(X)^\SYMG$ is faithful.
Let $\rho\in\OP{O}(H(X)^\SYMG)$ be the restriction
 of the action of the generator of $C_2$.
In general, a $K3$ surface
 which admits a non-symplectic
 action of a finite group is algebraic \cite{nikulin79fin}.
Therefore, there exists
 a primitive ample class $l\in H(X)^{\widetilde{G}}$.
By Theorem~\ref{THMsymplectic},
 we find that $(l)^\bot_{H(X)^\SYMG}=T(X)$ and
 $H(X)^{\widetilde{G}}=\Z l$.
Since we have $(\rho\otimes\C)(\omega_X)=-\omega_X$,
 the eigenvalues of $\rho$ are $1,-1,-1$.
In particular, it follows that
 $\rho\in\OP{SO}(H(X)^\SYMG)$.
We fix the $K3$ lattice $\Lambda_i$ with
 an action of $\SYMG$
 satisfying the conditions (i) and (ii)
 in Proposition~\ref{PROPinvariantlattice}
 and the condition $\Lambda_i^\SYMG\cong K_i$,
 which exists and is unique up to isomorphism
 by Proposition~\ref{PROPinvariantlattice}.
By the uniqueness of $\Lambda_i$,
 there exists an isomorphism
 $\alpha:H(X)\rightarrow\Lambda_i$
 compatible with the action of $\SYMG$
 (cf.\ Theorem~\ref{THMsymplectic}).

We consider the case
 where $H(X)^\SYMG\cong K_1\cong \Lambda_1^\SYMG$.
We identify $\Lambda_1^\SYMG$
 with the column vectors $\Z^3$
 with the Grammian matrix
\begin{equation*}
 \begin{pmatrix}
 4 & 1 & 0 \\
 1 & 4 & 0 \\
 0 & 0 & 20
 \end{pmatrix}.
\end{equation*}
The nontrivial elements
 in $\OP{SO}(\Lambda_1^\SYMG)\cong(\Z/2\Z)^2$ are
$$
 \rho_1=
 \begin{pmatrix}
    & -1 &   \\
 -1 &    &   \\
    &    & -1
 \end{pmatrix},~
 \rho_2=
 \begin{pmatrix}
   & 1 &   \\
 1 &   &   \\
   &   & -1
 \end{pmatrix},~
 \rho_3=
 \begin{pmatrix}
 -1 &    &   \\
    & -1 &   \\
    &    & 1
 \end{pmatrix}.
$$
By Lemma~\ref{LEMext_lambda},
 there exists a unique extension
 $\bar{\rho}_j\in\OP{O}(\Lambda_1)$
 of $\rho_j$
 which commutes with the action
 of $\SYMG$.
We consider the case
 $\rho=\alpha^{-1}\circ\rho_1\circ\alpha$.
By the uniqueness of the extension $\bar{\rho}_1$,
 we see that the action of $C_2$ on $H(X)$
 coincides with that of
 $\langle
 \alpha^{-1}\circ\bar{\rho}_1\circ\alpha
 \rangle$.
The invariant lattice
 $\Lambda_1^{\SYMG\times
 \langle \bar{\rho}_1 \rangle}$
 is generated by the following $l_1$:
$$l_1
 =\begin{pmatrix} 1 \\ -1 \\ 0 \end{pmatrix},~
 \BF{l_1}{l_1}=6,~
 (l_1)^\bot_{\Lambda_1^\SYMG} \cong
 \begin{pmatrix} 10 &  0 \\  0 & 20 \end{pmatrix},~
 [\Lambda_1^\SYMG
 :\Z l_1 \oplus (l_1)^\bot_{\Lambda_1^\SYMG}]=2.$$
Therefore, we find that $X$ corresponds the case
 in the first row in Table~\ref{TBLtranslattice}.
Similarly, we can determine
 the second and third rows
 in Table~\ref{TBLtranslattice}.
We have
$$ (\alpha\otimes\C)(\omega_X)
 =
 \sqrt{2}\begin{pmatrix} 1 \\ 1 \\ 0 \end{pmatrix}
 \pm\sqrt{-1}\begin{pmatrix} 0 \\ 0 \\ 1 \end{pmatrix},~
 \alpha(l)=\pm l_1$$
 for some $\omega_X\in H^{2,0}(X)$.
Moreover, we may assume that
$$(\alpha\otimes\C)(\omega_X)
 =\sqrt{2}\begin{pmatrix} 1 \\ 1 \\ 0 \end{pmatrix}
 +\sqrt{-1}\begin{pmatrix} 0 \\ 0 \\ 1 \end{pmatrix},~
 \alpha(l)=l_1,$$
 after replacing $\alpha$ with $\pm\bar{\rho}_j\circ\alpha$
 if necessary.
Since the $\bar{\rho}_j$ commute
 with each other,
 the isomorphism $\alpha$ is compatible
 with the action of
 $\SYMG\times\langle \bar{\rho}_1 \rangle$.
The uniqueness in the case
 corresponding to the first row
 in Table~\ref{TBLtranslattice}
 follows by Theorem~\ref{THMtorelli}.
On the other hand,
 since $(\Lambda_i)_\SYMG$ contains no vector
 $x$ with $\BF{x}{x}=-2$,
 the existence in this case follows
 by Theorem~\ref{THMtorelli} and \ref{THMsurjperiod}.
The remaining five cases in Table~\ref{TBLtranslattice}
 are similar.
\end{proof}

\section{A certain family of $K3$ surfaces %
 with an action of $\MF{S}_5$}
\label{SECTfamily}
Let $s_k\in\C[x_1,\ldots,x_5]$
 be the $k$-th power sum, i.e.,
 $s_k=x_1^k+\cdots+x_5^k$.
We set
\begin{equation*} 
 \X=\{
 ((x_1:\cdots:x_5),(t_0:t_1))\in\P^4\times\P^1\bigm|
 s_1(x)=t_0 s_4(x) + t_1 s_2(x)^2 = 0
 \}.
\end{equation*}
In the rest of the paper,
 we study the family
\begin{equation} \label{defp}
 p:\X\longrightarrow\P^1;\quad
 ((x_1:\cdots:x_5),(t_0:t_1))\longmapsto (t_0:t_1).
\end{equation}
For $\T\in\P^1$,
 the fiber $\X_\T:=p^{-1}(\T)$ is
 considered as a quartic surface in $\P^3$
 (perhaps singular).
If $\X_\T$ is nonsingular,
 it is a $K3$ surface.
We see that $\mathfrak{S}_5$ acts
 on the family
 by the permutation of
 the variables $x_1,\ldots,x_5$.
In this section,
 we determine singular fibers
 of the family
 and construct a local modification
 of the family
 around each singular fiber.
We use these results
 to determine the generic transcendental
 lattice of the family.

In Subsection~\ref{SUBSECTfamilyA},
 we determine singular fibers
 of the family.
We will find that
 each singular fiber corresponds
 to each conjugacy class
 of elliptic points
 of an arithmetic
 Fuchsian group
 in Section~\ref{SECTquaternion}.
In Subsection~\ref{SUBSECTfamilyB},
 we construct a local modification
 $\bar{\X}^{(\nu)}\rightarrow\Delta$
 of the family
 around each singular fiber,
 which is a smooth family of
 $K3$ surfaces
 with an action of
 $\MF{A}_5$ ($\subset\MF{S}_5$).

\subsection{Singular fibers of the family (\ref{defp})}
\label{SUBSECTfamilyA}
In this subsection,
 we determine singular fibers
 of the family defined by (\ref{defp}).
There are five singular fibers.
Four of them are quartic surfaces
 with $A_1$-singularities,
 so the minimal desingularization of them
 are $K3$ surfaces.

\begin{prop} \label{PROPsingpts}
We have the following.
\begin{itemize2}
\item[(i)]
There exists exactly five singular fibers in the family
$\{\X_\T\}_{\T\in\P^1}$.
We have the following list of singular fibers.
\begin{equation} \label{table_sing}
 \begin{array}{c|c|c|c|c}
   & \T
   & \begin{matrix} \rm{coordinate~of}
        \\ \raise0.5ex\hbox{\rm{one~singular~point}}\end{matrix}
   & \rm{type}
   & \begin{matrix}\rm{number~of}
        \\ \raise0.5ex\hbox{\rm{singular~points}}\end{matrix}
        \\ \hline
  \T_1 & (4:-1) & (1: -1: 1: -1: 0) & A_1 & 15 \\
  \T_2 & (2:-1) & (1: -1: 0: 0: 0) & A_1 & 10 \\
  \T_3 & (30:-7) & (2: 2: 2: -3: -3) & A_1 & 10 \\
  \T_4 & (20:-13) & (1: 1: 1: 1: -4) & A_1 & 5 \\
  \T_5 & (0:1) &
     & \text{\rm double quadric} &
 \end{array}
\end{equation}
In each case, except for $\X_{\T_5}$,
 the group $\mathfrak{S}_5$
 (or even $\mathfrak{A}_5$) acts on
 the set of singular points transitively,
 and we list the coordinate
 of one singular point on each fiber
 in the table.
\item[(ii)]
The set $\operatorname{Sing}\X$ is equal to
$\{\T=(0,1),\, s_1=s_2=s_4=0\}$.
\end{itemize2}
\end{prop}
\begin{proof}
 [Proof (cf.~\rm{\cite{vandergeer82}}\textit{).}]
In the case $\T=\T_5=(0:1)$, the assertion (i) is clear.
We assume that $\T=(t_0:t_1)\neq\T_5$ and put $t=t_1/t_0$.
Let $\X_\T$ be a singular fiber and
$\xi=(x_1:\cdots:x_5)$ a singular point on $\X_\T$.
Since $\X_\T$ is defined by $s_1=s_4+ts_2^2=0$, we have
$$x_i^3+ts_{2}x_i=u,\qquad i=1, \ldots, 5,$$
for some $u\in\C$ by the Jacobian criterion.
Hence $x_i$'s take at most three distinct values.
If $x_i$'s take exactly three distinct values,
then the sum of the three values is zero.
By considering the action of $\mathfrak{S}_5$,
we may assume that $\xi$ is of the form
$(\alpha: \beta: -(\alpha +\beta ): \alpha: \alpha)$ or
$(\alpha: \beta: -(\alpha +\beta ): \alpha: \beta)$
 without loss of generality.
In the former case, $s_1=0$ implies $\alpha=0$ and $t=-1/2$.
In the latter case, we have $\alpha+\beta=0$ and $t=-1/4$.
If $x_i$'s take exactly two distinct values, then we have
$\xi=(1: 1: 1: 1: -4)$ and $t=-13/20$, or
$\xi=(2: 2: 2: -3: -3)$ and $t=-7/30$.
We see that these singular points are of type $A_1$
 by computing the Hessian matrix.
The proof of the assertion (ii) is similar.
\end{proof}

\subsection{Local modifications by base change}
\label{SUBSECTfamilyB}
By Proposition~\ref{PROPsingpts},
 we see that
 the family $\{ \X_\T \}$
 is a smooth family
 of $K3$ surfaces on
 $\P^1\smallsetminus \{ \T_1,\ldots,\T_5 \}$.
On the other hand,
 the minimal desingularization
 of any quartic surface
 with only $A$-$D$-$E$-singularities
 is a $K3$ surface.
Hence, the minimal desingularization
 of $\X_{\T_\nu}$
 is a $K3$ surface for $\nu=1,\ldots,4$.
In this subsection,
 using the result in Atiyah~\cite{atiyah58},
 we construct a local modification
 $\bar{\X}^{(\nu)}\rightarrow\Delta$
 of the family around $\T_\nu$
 ($\nu=1,\ldots,4$).
This modification is
 a smooth deformation
 of the minimal desingularization
 of $\X_{\T_\nu}$ to a generic member
 of the family
 preserving the action
 of $\MF{A}_5$ ($\subset\MF{S}_5$).
We show a similar result
 in the case $\nu=5$,
 where $\X_{\T_5}$ is a double quadric
 surface.
The result is summarized in
 Proposition~\ref{PROPdeform}.

We recall a modification of
 a deformation of an $A_1$-singularity \cite{atiyah58}.
We set
\begin{equation} \label{defy_tau}
 \begin{matrix}
 \MC{Y}_\tau=\{(y_1,y_2,y_3)\in\C^3\bigm|
 y_1^2+y_2^2+y_3^2+\tau=0\},\hspace{4ex}\vspace{1ex} \\
 \MC{Y}_{~}=
 \{ ((y_1,y_2,y_3),\bar{\tau})\in\C^3\times\C\bigm|
 (y_1,y_2,y_3)\in\MC{Y}_{\bar{\tau}} \}.
 \end{matrix}
\end{equation}
We see that $\MC{Y}_\tau$ is nonsingular
 for $\tau\neq 0$ and
 $\MC{Y}_0$ has an $A_1$-singularity at the origin.
Let $\{\MC{Y}'_{\bar{\tau}}\}_{\bar{\tau}\in\C}$
 be the pull-back of
 $\{\MC{Y}_{\tau}\}_{\tau\in\C}$ by
 the base change $\tau=\bar{\tau}^2$
 and $\MC{Y}'\subset\C^3\times\C$ the total space of
 $\{\MC{Y}'_{\bar{\tau}}\}$, i.e.,
\begin{align*}
 \MC{Y}'_{\bar{\tau}}&=\{(y_1,y_2,y_3)\in\C^3\bigm|
 y_1^2+y_2^2+y_3^2+\bar{\tau}^2=0\}, \\
 \MC{Y}'_{~}&=
 \{ ((y_1,y_2,y_3),\bar{\tau})\in\C^3\times\C\bigm|
 (y_1,y_2,y_3)\in\MC{Y}'_{\bar{\tau}} \}.
\end{align*}
Let $\MC{Y}''$ be the blow-up of $\MC{Y}'$ at the origin and
$E\cong\P^1\times\P^1$ the exceptional divisor.
Finally, we contract $E$ to $\bar{E}\cong\P^1$
in $\MC{Y}'$ along a ruling of $E$ and
obtain a smooth threefold $\bar{\MC{Y}}$:
$$\begin{CD}
 \MC{Y}'' @= \MC{Y}'' \\
 @VVV @VVV\\
 \bar{\MC{Y}} @>>> \MC{Y}' @>{2:1}>> \MC{Y} \\
 @VVV @VVV @VVV \\
 \C @= \C @>{2:1}>> ~~\C~.
\end{CD}$$
Under a natural map
$\bar{\MC{Y}}\rightarrow\MC{Y}'\rightarrow\C$,
we obtain a smooth family
$\{\bar{\MC{Y}}_{\bar{\tau}}\}_{\bar{\tau}\in\C}$
with the total space $\bar{\MC{Y}}$.
We see that $\bar{\MC{Y}}_{\bar{\tau}} =
\MC{Y}_{\bar{\tau}^2}$
for $\bar{\tau}\neq 0$.
On the other hand,
 a natural map $\bar{\MC{Y}}_0\rightarrow\MC{Y}_0$ is
 the minimal resolution
 of the $A_1$-singularity of $\MC{Y}_0$
 with the exceptional divisor $\bar{E}$.
Note that there are two possibilities
of the choice of a contraction $E\rightarrow\bar{E}$.
Let $\varphi'$ be an automorphism of
 the family $\MC{Y}'\rightarrow\C$.
Then, $\varphi'$ lifts
 to an automorphism of
 the family $\bar{\MC{Y}}\rightarrow\C$
 if and only if
 the lift of $\varphi'$
 to $\MC{Y}''$ preserves
 the contraction $E\rightarrow\bar{E}$.

\begin{lem} \label{LEMquadric}
Let $Q$ be a nonsingular quadric surface in $\P^3$
 defined by a polynomial $q$
 and $g$ an element in $\OP{GL}(4,\C)$
 with $g^* q=aq$ for some $a\in\C^\times$.
Then, the following conditions are equivalent
 to each other:
\begin{itemize2}
\item[(i)]
$g$ preserves each ruling $Q\rightarrow\P^1$ of $Q$;
\item[(ii)]
the action of $g$ on $H^2(Q,\Z)$ is trivial;
\item[(iii)]
$\det(g)=a^2$.
\end{itemize2}
\end{lem}
\begin{proof}
The assertion
 (i)$\Leftrightarrow$(ii) is trivial.
It is sufficient to show
 the assertion (ii)$\Leftrightarrow$(iii)
 in the case $a=1$.
The group of elements in $\OP{GL}(4,\C)$
 whose actions fix $q$ is $\OP{O}(4,\C)$.
We consider $H^2(Q,\Z)\cong\Z\oplus\Z$
 as a $\Z$-module.
The image of a natural homomorphism
 $\OP{O}(4,\C)\rightarrow \OP{GL}(H^2(Q,\Z))$ is
 an abelian group of order $2$.
Since $\OP{O}(4,\C)$ has
 exactly two connected components,
 the kernel of this homomorphism
 is the identity component of $\OP{O}(4,\C)$,
 i.e.\ $\OP{SO}(4,\C)$.
\end{proof}

In the rest of the subsection,
 we construct local modifications of
 the family $p:\X\rightarrow\P^1$ around singular fibers,
 which we determine in Proposition~\ref{PROPsingpts}.
Let $\Delta=\{z\in\C\bigm| |z|<\varepsilon\}$
 be a disc for sufficiently small $\varepsilon>0$
 and $U_\nu\cong\Delta$ a neighborhood of $\T_\nu$
 such that $\X_\T$ is nonsingular
 for all $\T\in U_\nu\smallsetminus \{ \T_\nu \}$
 ($\nu=1,\ldots,5$).
We fix an identification
 $(U_\nu,\T_\nu)=(\Delta,0)$
 and define the family
 $p^{(\nu)}:\X^{(\nu)}\rightarrow\Delta$
 as the restriction of the family
 $p:\X\rightarrow\P^1$ to $U_\nu=\Delta$:
\begin{equation}\begin{CD}
 \quad\X^{(\nu)} @>{\subset}>> \X \\
 @V{p^{(\nu)}}VV            @VV{p}V \\
 (\Delta,0)\, @>{\subset}>> \,(\P^1,\T_\nu).
\end{CD}\end{equation}

\begin{prop} \label{PROPdeform}
For $\nu=1,\ldots,5$,
 there exists a modification
 $\bar{p}^{(\nu)}:\bar{\X}^{(\nu)}\rightarrow\Delta$
 of the family
 $p^{(\nu)}:\X^{(\nu)}\rightarrow\Delta$
 around $0\in\Delta$
 by the base change $\delta:z\mapsto z^2$:
\begin{equation} \label{cd_modi}
 \begin{CD}
 \quad\bar{\X}^{(\nu)} @>>> \hspace{2ex}\X^{(\nu)} \\
 @V{\bar{p}^{(\nu)}}VV     @VV{p^{(\nu)}}V \\
 \Delta @>>{\delta:z\mapsto z^2}> \Delta
 \end{CD}
\end{equation}
which is a smooth family of $K3$ surfaces
 and satisfies the following conditions:
\begin{itemize2}
\item[(i)]
 $\bar{\X}^{(\nu)}|_{\Delta^*}\rightarrow\Delta^*$
 is the pull-back of
 $\X^{(\nu)}|_{\Delta^*}\rightarrow\Delta^*$
 by the base change $\delta$,
 where $\Delta^*=\Delta\smallsetminus\{0\}$.
\item[(ii)]
 $\bar{\X}^{(\nu)}_0 \rightarrow \X^{(\nu)}_0$
 is the minimal resolution
 (resp.\ a double covering) of $\X^{(\nu)}_0$
 if $\nu=1,\ldots,4$ (resp.\ $\nu=5$);
\item[(iii)]
 the action of $\mathfrak{A}_5$
 (resp.\ $\mathfrak{S}_5$) on the family
 $p^{(\nu)}:\X^{(\nu)}\rightarrow\Delta$ lifts to
 that on the family
 $\bar{p}^{(\nu)}:\bar{\X}^{(\nu)}\rightarrow\Delta$
 if $\nu=1,\ldots,4$ (resp.\ $\nu=5$).
\end{itemize2}
\end{prop}

First, we consider the cases $\nu=1,\ldots,4$.
We fix a singular point $\xi=(x_1:\cdots:x_5)$
 on $\X^{(\nu)}_0=\X_{\T_\nu}$.
Since we have $s_2(x)\neq 0$
 by Proposition~\ref{PROPsingpts},
 the family $\X^{(\nu)}\rightarrow\Delta$
 is locally isomorphic to
 the family $\MC{Y}\rightarrow\C$ given
 in (\ref{defy_tau}) near $\xi$.
Therefore, we can modify
 the family $\X^{(\nu)}\rightarrow\Delta$
 near $\xi$
 in the same way as
 the family $\MC{Y}\rightarrow\C$.
Let $\MC{Z}'\rightarrow\Delta$ be the pull-back
 of $\X^{(\nu)}\rightarrow\Delta$
 by the base change $\delta:z\mapsto z^2$.
Let $\MC{Z}''$ be the blow-up of $\MC{Z}'$
 at the point in $\MC{Z}'$
 corresponding to $\xi$
 with the exceptional divisor
 $E\cong\P^1\times\P^1$.
We fix a ruling
 $\pi:E\rightarrow\bar{E}\cong\P^1$.
We contract $E$ to $\bar{E}$
 along the ruling $\pi$
 and obtain a family
 $\bar{\MC{Z}}\rightarrow\Delta$,
 which is a smooth family near $\bar{E}$.
Let $\sigma$ be an element in $\mathfrak{S}_5$
 which fixes $\xi$.
The action of $\sigma$
 on $\X^{(\nu)}\rightarrow\Delta$
 lifts to that on $\MC{Z}''\rightarrow\Delta$.
This action of $\sigma$ induces
 an action on $\bar{\MC{Z}}\rightarrow\Delta$
 if and only if
 the action of $\sigma$ on $E$
 preserves the ruling $\pi$.
We set $V=\{x\in\C^5\bigm|s_1(x)=0\}$
 and consider $\MC{Z}$
 as a hypersurface of $\P(V)\times\Delta$.
Let $W_1$ (resp.\ $W_2$) be
 the tangent space of $\P(V)$ (resp.\ $\Delta$)
 at $\xi$ (resp.\ the origin).
Then, the tangent space of $\P(V)\times\Delta$
 at $(\xi,0)$ is $W_1\times W_2\cong\C^3\times\C$
 and the exceptional divisor $E$ is considered
 as a quadric surface
 in $\P(W_1\times W_2)\cong\P^3$.
Since $\sigma$ fixes $\xi=(x_1:\cdots:x_5)$,
 we have $\sigma\cdot x=\lambda x$
 for some $\lambda\in\C^\times$.
Hence, we find that
 $\sigma$ is represented by
 $\begin{pmatrix}
 A    & 0 \\
 * & \lambda
 \end{pmatrix}$
 for some $A\in \OP{GL}(3,\C)$
 as an element in $\OP{GL}(V)\cong\OP{GL}(4,\C)$
 and the action of $\sigma$ on $W_1\cong\C^3$
 is represented by the matrix $\lambda^{-1}A$.
We have $\OP{sgn}(\sigma)=\lambda\det(A)$.
We denote the action of
 $\sigma$ on $W_1\times W_2\cong\C^4$
 by $\bar{\sigma}$.
We can see that $\bar{\sigma}$
 is represented by the matrix
 $\begin{pmatrix}
 \lambda^{-1}A &   \\
               & 1
 \end{pmatrix}.$
Let $q$ be the defining equation of
 $E$ in $\P(W_1\times W_2)$.
Since the defining equation of
 $\X^{(\nu)}_0=\X_{\T_\nu}$
 is invariant under $\sigma$,
 we have $\bar{\sigma}^\ast q=\lambda^{-2} q$
 by $\sigma\cdot x=\lambda x$.
On the other hand, we have
$$\det(\bar{\sigma})=\lambda^{-3}\det(A)
 =\lambda^{-3}\cdot(\lambda^{-1}\OP{sgn}(\sigma))
 =\lambda^{-4}\OP{sgn}(\sigma).$$
Let $H= \{ \sigma\in\MF{A}_5 \bigm|
 \sigma\cdot\xi=\xi \}$
 be the stabilizer of $\xi$ in $\mathfrak{A}_5$.
By Lemma~\ref{LEMquadric},
 the action of $H$
 on $\X^{(\nu)}\rightarrow\Delta$ lifts to
 that on
 $\bar{\MC{Z}}\rightarrow\Delta$.
We apply the same modification,
 i.e.\ blowing-up and contraction, to
 other singular points on $\bar{\MC{Z}}$
 and obtain a smooth family
 $\bar{p}^{(\nu)}:\bar{\X}^{(\nu)}\rightarrow\Delta$.
Here we chose contractions of exceptional divisors
 such that the action of $\MF{A}_5$
 preserves them.
Since the minimal resolution of
 a quartic surface in $\P^3$ with
 only $A$-$D$-$E$-singularities is a $K3$ surface,
 the family $\bar{\X}^{(\nu)}\rightarrow\Delta$ is
 a smooth family of $K3$ surfaces.

\begin{rem}
For $\nu=1,\ldots,4$,
 a local modification of the family
 $p^{(\nu)}:\X^{(\nu)}\rightarrow\Delta$
 around $0\in\Delta$ is not unique.
It depends on the choice of contractions
 $\bar{E}\rightarrow E$
 of exceptional divisors
 at singular points
 on $\X^{(\nu)}_0$.
\end{rem}

Next, we consider the case $\nu=5$.
We may assume that the family
 $p_5:\X^{(5)}\rightarrow\Delta$ is given by
\begin{gather*}
 \X^{(5)}
 = \{ ((x_1:\cdots:x_5),z)\in\P^5\times\Delta\bigm|
 s_1(x)=z s_4(x)+s_2(x)^2=0 \}, \\
 p_5:\X^{(5)}\rightarrow\Delta;\quad
 ((x_1:\cdots:x_5),z)\longmapsto z
\end{gather*}
 after shrinking $U_5=\Delta$ and
 changing the coordinate $z\in\Delta$
 if necessary.
We define the family
 $\bar{p}_5:\bar{\X}^{(5)}\rightarrow\Delta$ by
\begin{gather*} \label{defeqX5}
 \bar{\X}^{(5)}
 = \{ ((x_1:\cdots:x_5; \mu), z)
 \in\P(1^5,2)\times\Delta\bigm|
 s_1(x)=s_2(x)-z\mu=s_4(x)-\mu^2=0 \}, \\
 \bar{p}_5:\bar{\X}^{(5)}\longrightarrow\Delta;\quad
 ((x_1:\cdots:x_5;\mu),z)\longmapsto z.
\end{gather*}
Here $\P(1^5,2)$ is the weighted projective space
 of weight $(1,1,1,1,1,2)$
 and
 $(x_1:\cdots:x_5;\mu)
 =(x'_1:\cdots:x'_5;\mu')\in\P(1^5,2)$
 if and only if
 $(x_1,\ldots,x_5,\mu)
 =(\lambda x'_1,\ldots,\lambda x'_5,\lambda^2\mu')$
 for some $\lambda\in\C^\times$.
We have the commutative diagram
$$\begin{CD}
\quad\bar{\X}^{(5)} @>{\varphi}>> \quad\X^{(5)} \\
@V{\bar{p}_5}VV       @VV{p_5}V \\
\Delta @>>{\delta:z\mapsto z^2}> ~~\Delta~,
\end{CD}$$
where we define
 $\varphi:((x_1:\cdots:x_5;\mu),z)
 \mapsto ((x_1:\cdots:x_5),z^2)$.
The action of $\mathfrak{S}_5$
 on the family $p_5:\X^{(5)}\rightarrow\Delta$
 lifts to that on the family
 $\bar{p}_5:\bar{\X}^{(5)}\rightarrow\Delta$.
The assertions (i)--(iii) follows immediately.
It is sufficient to show the following lemma.

\begin{lem} \label{LEMinfty}
$\bar{p}_5:\bar{\X}^{(5)}\rightarrow\Delta$
 is a smooth family of $K3$ surfaces.
\end{lem}
\begin{proof}
We shall show that $\bar{\X}^{(5)}_0$ is nonsingular.
Let $\xi=(x_1:\cdots:x_5;\mu)$ is
 a singular point on $\bar{\X}^{(5)}_0$.
We may assume that $x_1\neq 0$
 without loss of generality.
We take the local coordinate system
 $(\bar{x}_2,\ldots,\bar{x}_5,\bar{\mu})
 =(x_2/x_1,\ldots,x_5/x_1,\mu/x_1^2)$
 of $\P(1^5,2)$ at $\xi$.
Since $\bar{\X}^{(5)}_0$ is defined by
 $s_1(x)=s_2(x)=s_4(x)-\mu^2=0$,
 the Jacobian matrix at $\xi$ is
$$\begin{pmatrix}
    1  & \cdots & 1      & 0 \\
 2\bar{x}_2  & \cdots & 2\bar{x}_5   & 0 \\
 4\bar{x}_2^3 & \cdots & 4\bar{x}_5^3 & -2\bar{\mu}
\end{pmatrix}.$$
Since we have
$$
 x_1\begin{pmatrix} 1 \\ 2
     \\ 4 \end{pmatrix}
 +
 x_2\begin{pmatrix} 1 \\ 2\bar{x}_2
     \\ 4\bar{x}_2^3 \end{pmatrix}
 +\cdots+
 x_5\begin{pmatrix} 1 \\ 2\bar{x}_5
     \\ 4\bar{x}_5^3 \end{pmatrix}
 +
 \frac{2\mu}{x_1}\begin{pmatrix} 0 \\ 0
     \\ -2\bar{\mu} \end{pmatrix}=0,$$
the following matrix is singular:
$$\begin{pmatrix}
    1  & \cdots & 1      & 0 \\
 2x_1  & \cdots & 2x_5   & 0 \\
4x_1^3 & \cdots & 4x_5^3 & -2\mu
\end{pmatrix}.$$
Hence we have
$$\det \begin{pmatrix}
1 & 1 & 1 \\
x_i & x_j & x_k \\
x_i^3 & x_j^3 & x_k^3
\end{pmatrix}
= (x_i-x_j)(x_j-x_k)(x_k-x_i)(x_i+x_j+x_k)= 0
$$
for all $i,j,k=1,\ldots,5$.
We see that $x_i$'s take
 three distinct values the sum of which is zero,
 or two distinct values.
By the same argument
 as in the proof of Proposition~\ref{PROPsingpts},
 we find that such $\xi$ does not exist.
Therefore $\bar{\X}^{(5)}_0$ is nonsingular.
Since $\bar{\X}^{(5)}_0\rightarrow\X^{(5)}_0$ is
 the double covering branched along
 $\{s_1=s_2=s_4=0\}$,
 which is a divisor of
 $\X^{(5)}_0\cong\P^1\times\P^1$
 of bi-degree $(4,4)$,
 it follows that
 $\bar{\X}^{(5)}_0$ is a $K3$ surface.
\end{proof}

\section{Transcendental lattices of singular fibers and %
 a generic fiber of the family (\ref{defp})}
\label{SECTtrans}
In this section,
 we study the period map
 of the family $\{ \X_\T \}$.
To do this,
 we determine the transcendental lattice $T$
 of generic $\X_\T$.
Moreover,
 we determine
 the transcendental lattice of
 ${\bar{\X}^{(\nu)}_0}$
 for $\nu=1,\ldots,5\,$
 (see Proposition~\ref{PROPdeform}).

\subsection{The period map of the family %
$\{ \X_\T \}$}
\label{SUBSECTtransA}
In this subsection,
 we define the period map of the family
 $\{ \X_\T \}$.
For generic
 $\T\in\P^1\smallsetminus\{ \T_\nu \}$,
 the structure of
 the transcendental lattice
 $T(\X_\T)$ of $\X_\T$
 is unique up to isomorphism.
We denote it by $T$.
We fix a generic point
 $\mathsf{s}\in\P^1\smallsetminus\{ \T_\nu \}$
 and identifications
 $\Lambda=H(\X_{\mathsf{s}})$,
 $T=T(\X_{\mathsf{s}})$.
The group $\MF{S}_5$ acts on
 $\Lambda$ through the identification
 $\Lambda=H(\X_{\mathsf{s}})$.
Let $l_\T\in H(\X_\T)$ be
 the natural polarization class of
 $\X_\T$.
We set $l=l_{\mathsf{s}}\in\Lambda$,
 $\Lambda_\PRIM=(l)^\bot_\Lambda$
 and $H_\PRIM(\X_\T)=(l_\T)^\bot_{H(\X_\T)}$.

\begin{defn}
For $\T\in\P^1\smallsetminus\{ \T_\nu \}$,
 a polarized marking $\alpha$
 of $\X_\T$ is
 an isomorphism
 $\alpha:H(\X_\T)
 {\displaystyle\mathop{\rightarrow}^{\sim}}
 \Lambda$
 such that $\alpha(l_\T)=l$.
In addition,
 if $\alpha$ is compatible with
 the action of $\MF{S}_5$,
 we call $\alpha$
 a polarized $\MF{S}_5$-marking
 of $\X_\T$.
\end{defn}

Since the family $\{ \X_\T \}$ is
 a smooth family with an action of $\MF{S}_5$
 for $\T\in\P^1\smallsetminus\{\T_\nu\}$,
 there exists
 a polarized $\MF{S}_5$-marking of $\X_\T$
 for any $\T\in\P^1\smallsetminus\{\T_\nu\}$.
We set
\begin{equation}
 \MC{M}=\{
 (\T,\alpha) \bigm|
 \T\in\P^1\smallsetminus\{\T_\nu\},~
 \alpha:
 \text{a polarized $\MF{S}_5$-marking of $\X_\T$}
 \}
\end{equation}
 and
\begin{equation}
 \OP{Aut}_{\MF{S}_5}(\Lambda,l)
 =\{ \varphi\in\OP{O}(\Lambda) \bigm|
 \varphi(l)=l,~
 \varphi\sigma=\sigma\varphi\,
 (\forall\sigma\in\MF{S}_5) \}.
\end{equation}
Then, $\MC{M}$ becomes
 a complex manifold.
Moreover, the natural map
 $\MC{M}\rightarrow
 \P^1\smallsetminus\{\T_\nu\}$
 is a principal
 $\OP{Aut}_{\MF{S}_5}(\Lambda,l)$-bundle.
For $(\T,\alpha)\in\MC{M}$,
 the period of the marked $K3$ surface
 $(\X_\T,\alpha)$
 is defined as
 $(\alpha\otimes\C)(H^{2,0}(\X_\T))
 \in\P(\Lambda_\PRIM\otimes\C)$.
By Riemann bilinear relations,
 the period of $(\X_\T,\alpha)$ is
 contained in the period domain $\Omega_{\Lambda_\PRIM}$
 defined by
\begin{equation}
 \Omega_{\Lambda_\PRIM}=
 \{ \C\omega\in\P(\Lambda_\PRIM\otimes\C) \bigm|
 \BF{\omega}{\omega}=0,~
 \BF{\omega}{\bar{\omega}}>0 \}.
\end{equation}
Let $\PER'_0$ be the map defined by
\begin{equation}
 \PER'_0:\MC{M}\longrightarrow\Omega_{\Lambda_\PRIM};
 \quad
 (\T,\alpha)\longmapsto(\alpha\otimes\C)(H^{2,0}(\X_\T)).
\end{equation}
Then, the map $\PER'_0$ is holomorphic.
For $(\T,\alpha)\in\MC{M}$,
 the action of $\MF{A}_5$ on
 $\X_\T$ is symplectic
 because $\MF{A}_5$ is
 a non-commutative simple group.
Hence, we have
 $\PER'_0(\T,\alpha)
 \in\P(\Lambda_\PRIM^{\MF{A}_5}\otimes\C)$.
We will determine
 the generic transcendental lattice $T$
 in the following subsection
 (see Lemma~\ref{LEMtransgen}).
In particular,
 we have $T=\Lambda_\PRIM^{\MF{A}_5}$
 and $T$ is of signature $(2,1)$.
As a consequence,
 we have $\PER'_0(\T,\alpha)\in
 \Omega_T=\P(T\otimes\C)\cap\Omega_{\Lambda_\PRIM}$
 for any $(\T,\alpha)\in\MC{M}$.
The period domain
 $\Omega_T$ is isomorphic to the union of
 two copies of
 the upper half-plane $\MF{H}$
 (i.e.\ $\P^1$ minus
 a real projective line)
 and
 the natural action of $\OP{O}(T)$
 on $\Omega_T$ is
 properly discontinuous.
Also, $\OP{Aut}_{\MF{S}_5}(\Lambda,l)$
 acts on $\Omega_T$.
By forgetting markings,
 the map $\PER'_0$ descends to
 the following map $\PER'$:
\begin{equation}\begin{CD}
 \MC{M} @>{\PER'_0}>> \Omega_T \\
 @VVV @VVV \\
 \P^1\smallsetminus\{\T_\nu\} @>{\PER'}>>
 \OP{Aut}_{\MF{S}_5}(\Lambda,l)\backslash\Omega_T.
\end{CD}\end{equation}

\begin{defn}
We call $\PER'$ the period map of
 the family
 $\{ \X_\T \}_{\T\in\P^1\smallsetminus\{\T_\nu\}}$.
\end{defn}

\begin{lem} \label{LEMpi_injective}
The period map
$\PER'$ is injective.
\end{lem}
\begin{proof}
Let $\T,\T'$ be elements
 in $\P^1\smallsetminus\{\T_\nu\}$
 such that $\PER'(\T)=\PER'(\T')$.
By the definition of $\PER'$,
 there exists an isomorphism
 $\varphi:H(\X_\T)\rightarrow H(\X_{\T'})$
 which preserves the Hodge structure,
 maps $l_\T$ to $l_{\T'}$,
 and is compatible with the action
 of $\MF{S}_5$.
Set $V=\{ s_1=0 \}\subset\C^5$.
Recall that $\X_\T,\X_{\T'}$ are
 quartic surfaces in $\P(V)\cong\P^3$
 (see (\ref{defp})).
By Theorem~\ref{THMtorelli},
 there exists a projective transformation
 $f$ of $\P(V)$ such that
 $f(\X_{\T'})=\X_{\T}$,
 (the restriction of) $f$ induces $\varphi$, and
 $f$ is compatible with
 the action of $\MF{S}_5$.
Let $\bar{f}$ be an element in $\OP{GL}(V)$
 which induces $f$.
We have
 $\sigma \bar{f} \sigma^{-1}
 =\lambda(\sigma) \bar{f}$
 with $\lambda(\sigma)\in\C^\times$
 for any $\sigma\in\MF{S}_5$.
Since $\MF{A}_5$ is a non-commutative simple group,
 we have $\lambda(\sigma)=1$ for any $\sigma\in\MF{A}_5$.
Since $V$ as a representation space of $\MF{A}_5$ is
 still irreducible,
 $f$ is trivial by Schur's lemma.
As a consequence, we have $\X_\T=\X_{\T'}\subset V$,
 i.e.\ $\T=\T'$.
\end{proof}

\begin{rem} \label{REMnofixedpoint}
By Lemma~\ref{LEMpi_injective},
 the composite map
 $\MC{M}\rightarrow
 \P^1\smallsetminus\{\T_\nu\}
 \rightarrow\OP{Im}(\PER')$
 is an unramified covering.
Hence,
 there exists no fixed point
 in $\OP{Im}(\PER'_0)$
 under the action of
 $\OP{Aut}_{\MF{S}_5}(\Lambda,l)$.
\end{rem}

We extend the period map $\PER'$
 to $\P^1$ using the modifications
 $\bar{p}^{(\nu)}:\bar{\X}^{(\nu)}
 \rightarrow\Delta$
 in Proposition~\ref{PROPdeform},
 as follows.
For $\nu=1,\ldots,5$,
 we chose a trivialization
 $\alpha^{(\nu)}:R^2(\bar{p}^{(\nu)})_*\Z
 {\displaystyle\mathop{\rightarrow}^{\sim}}
 \Lambda$
 on $\Delta$
 such that
 $\alpha^{(\nu)}_z:
 H(\bar{\X}^{(\nu)}_z)\rightarrow\Lambda$
 is a polarized $\MF{S}_5$-marking
 for all $z\neq 0$.
Let $\overline{\PER}'$ be the extension of $\PER'$
 defined by
\begin{equation} \label{def_barpi}
 \overline{\PER}'(z^2)
 =[(\alpha^{(\nu)}_z\otimes\C)
 (H^{2,0}(\bar{\X}^{(\nu)}_z))]
 \in\OP{Aut}_{\MF{S}_5}(\Lambda,l)
 \backslash\Omega_T,\quad
 z\in\Delta,
\end{equation}
for $\nu=1,\ldots,5$.
This definition of $\overline{\PER}'$ does not
 depend on the choice of $\alpha^{(\nu)}$,
 and $\overline{\PER}'$ is holomorphic.

\begin{defn}
We call $\overline{\PER}'$
 the period map of the family
 $\{ \X_\T \}_{\T\in\P^1}$.
\end{defn}

We will see later that
 the natural map
 $\OP{Aut}_{\MF{S}_5}(\Lambda,l)
 \rightarrow\OP{O}(T)/\pm 1$
 is surjective
 (see Lemma~\ref{LEMtransgen}).
Since the action of $-1\in\OP{O}(T)$
 on $\Omega_T$ is trivial,
 we have 
 $\OP{Aut}_{\MF{S}_5}(\Lambda,l)
 \backslash\Omega_T
 = \OP{O}(T)\backslash\Omega_T$.
Set
 $\OP{O}(T)^\circ
 =\OP{O}(T)\cap \OP{O}(T\otimes\R)^\circ$,
 where $\OP{O}(T\otimes\R)^\circ$
 is the identity component of
 $\OP{O}(T\otimes\R)$.
We can check that
 $[\OP{O}(T):\OP{O}(T)^\circ]=4$
 by the Grammian matrix
 (\ref{gramtransgen}) of $T$,
 so $\OP{O}(T)\backslash\Omega_T$
 is a connected nonsingular curve.
Hence, Lemma~\ref{LEMpi_injective}
 implies the following.

\begin{prop} \label{PROPpisbiholomorphic}
The period map $\overline{\PER}'$ is biholomorphic.
\end{prop}

We define the domain
 $\MC{D}_T$ isomorphic to $\Omega_T$.
In the rest of the paper,
 we fix an orientation of
 $T$ (or $T\otimes\R$).
Set
\begin{equation} \label{domain_T}
 \MC{D}_T=\{
 v\in T\otimes\R \bigm| \BF{v}{v}=-1 \}.
\end{equation}
We identify $\Omega_T$
 with $\MC{D}_T$ by the correspondence
\begin{equation} \label{identification_omega}
 \Omega_T\longleftrightarrow\MC{D}_T;\quad
 \C\omega\longleftrightarrow v,
\end{equation}
where $\omega$ and $v$ satisfy
 $\omega\bot v$ and
 $\OP{Re}\omega\wedge\OP{Im}\omega
 \wedge v>0$.
Let $\MC{D}_T^\circ$ be a connected component
 of $\MC{D}_T$.
Then, we have
 $\MC{D}_T=\MC{D}_T^\circ\sqcup-\MC{D}_T^\circ$.

The following composite map $\PER$ is also called
 the period map of $\{ \X_\T \}$:
\begin{equation} \label{def_pi}
 \PER:\P^1
 \mathop{\longrightarrow}^{\sim}_{\displaystyle\overline{\PER}'}
 \OP{Aut}_{\MF{S}_5}(\Lambda,l)\backslash\Omega_T
 =\OP{O}(T)\backslash\Omega_T
 \cong
 \OP{O}(T)^\circ\backslash\MC{D}_T^\circ.
\end{equation}

\subsection{Determination of transcendental lattices}
\label{SUBSECTtransB}
In this subsection,
 we determine the transcendental lattices
 of generic $\X_\T$ and
 $\bar{\X}^{(\nu)}_0$ ($\nu=1,\ldots,5$)
 defined in Proposition~\ref{PROPdeform}.

\begin{lem} \label{LEMtransinf}
$T({\bar{\X}^{(5)}_0})\cong
{\displaystyle\begin{pmatrix}
20 & 10 \\
10 & 20
\end{pmatrix}}$.
\end{lem}
\begin{proof}
Let $\iota$ be
 the nontrivial covering transformation of
 ${\bar{\X}^{(5)}_0}\rightarrow\X^{(5)}_0$.
The group
 $\widetilde{G}:=G\times \langle \iota \rangle$
 acts on $\bar{\X}^{(5)}_0$,
 where
 $G:=\MF{A}_5\cup \{ \iota\sigma \bigm|
 \sigma\in\MF{S}_5\smallsetminus\MF{A}_5 \}
 \cong\MF{S}_5$.
The action of $G$ is symplectic.
Now we apply Theorem~\ref{THMs5x2}.
Since the action
 of $\widetilde{G}$ preserves
 the natural polarization
 of ${\bar{\X}^{(5)}_0}$,
 which is of degree $4$,
 the assertion follows
 from Table~\ref{TBLtranslattice}.
\end{proof}

\begin{lem} \label{LEMtrans0}
$T({\bar{\X}^{(2)}_0})\cong
{\displaystyle\begin{pmatrix}
4 & 1 \\
1 & 4
\end{pmatrix}}$.
\end{lem}
\begin{proof}
The defining equation of $\X_{\T_2}$ is
$$\sigma_1(x)=\sigma_4(x)=0,$$
where $\sigma_k(x)\in\C[x_1,\ldots,x_5]$ is
the $k$-th elementary symmetric polynomial.
Let $\iota$ be a birational automorphism of $\X_{\T_2}$
given by $(x_1:\ldots:x_5)\mapsto(1/x_1:\ldots:1/x_5)$.
We also denote by $\iota$
 the (biholomorphic) automorphism of ${\bar{\X}^{(2)}_0}$
induced by $\iota$.
Since there is no fixed point
 on ${\bar{\X}^{(2)}_0}$
 under the action of $\iota$,
 the action of $\iota$ is
 a non-symplectic automorphism
 of ${\bar{\X}^{(2)}_0}$
 by Theorem~\ref{THMfixedpts}.
The group
 $\widetilde{G}:=G\times \langle \iota \rangle$
 acts on $\bar{\X}^{(5)}_0$,
 where
 $G:=\MF{A}_5\cup \{ \iota\sigma \bigm|
 \sigma\in\MF{S}_5\smallsetminus\MF{A}_5 \}
 \cong\MF{S}_5$.
The action of $G$ is symplectic.
To apply Theorem~\ref{THMs5x2},
 we shall determine the invariant lattice
 $H(\bar{\X}^{(2)}_0)^{\widetilde{G}}$
 of rank $1$,
 or a polarization fixed under
 the action of $\widetilde{G}$.

Let $D_{ij}$, $1\leq i<j\leq 5$,
 be the smooth rational curve
 on ${\bar{\X}^{(2)}_0}$
 which maps to
 $\{ x_i=x_j=\sigma_1(x)=0 \}
 \subset \X_{\T_2}$.
Let
 $d_{ij} \in
 H(\bar{\X}^{(2)}_0)$
 be the class of $D_{ij}$,
 $d$ the sum of the $d_{ij}$
 and $\delta$ the sum of the classes
 of the ten smooth rational curves
 on ${\bar{\X}^{(2)}_0}$
 which contract
 to $A_1$-singularities on $\X_{\T_2}$
 (see Proposition~\ref{PROPsingpts}).
We have
\begin{gather*}
d_{ij}d_{kl}=
\begin{cases}
-2 & \text{if~~} d_{ij}=d_{kl}, \\
 0 & \text{if~~} d_{ij}\neq d_{kl},
\end{cases} \\
d_{ij}\delta=3,\quad d\delta=30,\quad d^2=\delta^2=-20, \\
\iota^\ast d=\delta.
\end{gather*}
Since we have
 $\theta:=d+\delta \in
 H(\bar{\X}^{(2)}_0)^{\widetilde{G}}$
 and $\theta^2=20$,
 it follows by Table~\ref{TBLtranslattice}
 that $T({\bar{\X}^{(2)}_0})$
 is isomorphic to either
 $\begin{pmatrix} 4 & 1 \\ 1 & 4 \end{pmatrix}$
 or
 $\begin{pmatrix} 4 & 2 \\ 2 & 16 \end{pmatrix}$.
Assume that the latter case occurs.
Then, by Table~\ref{TBLtranslattice},
 there exists $v\in T({\bar{\X}^{(2)}_0})$
 such that
 $(v+\theta)/2\in H(\bar{\X}^{(2)}_0)$.
However,
 we have $d_{ij}\cdot (v+\theta)/2=1/2\not\in\Z$,
 which is a contradiction.
Hence the former case occurs.
\end{proof}

\begin{rem} \label{REMnontrivial}
Since we have
 $T({\bar{\X}^{(2)}_0})\not\cong
 T({\bar{\X}^{(5)}_0})$
 by Lemma~\ref{LEMtransinf}
 and \ref{LEMtrans0},
 the family $\{\X_\T\}$
 is not locally trivial
 by Proposition~\ref{PROPdeform}.
\end{rem}

Now we determine
 the generic transcendental lattice $T$
 of the family $\{ \X_\T \}$.
Recall that
 we fix generic
 $\mathsf{s}\in\P^1
 \smallsetminus \{ \T_\nu \}$
 and identifications
 $\Lambda=H(\X_{\mathsf{s}})$,
 $T=T(\X_{\mathsf{s}})$
 (see Subsection~\ref{SUBSECTtransA}).

\begin{lem} \label{LEMtransgen}
We have the following.
\begin{itemize2}
\item[(i)]
 $\Lambda_{\MF{A}_5}$ is
 isomorphic to $K$
 in Lemma~\ref{LEMsurjectivityA5}.
\item[(ii)]
 $T$ is equal to
 $\Lambda_\PRIM^{\MF{A}_5}$ and
 isomorphic to
\begin{equation} \label{gramtransgen}
 \begin{pmatrix}
 4 & 1 & 0 \\
 1 & 4 & 0 \\
 0 & 0 & -20
 \end{pmatrix}\text{.}
\end{equation}
\end{itemize2}
\end{lem}
\begin{proof}
By Proposition~\ref{PROPdeform},
 there exists a deformation of
 $\bar{\X}^{(2)}_0$
 to $\X_{\mathsf{s}}$
 preserving the action of $\MF{A}_5$
 and the natural (quasi-)polarization.
This induces an isomorphism
 $\varphi:H(\bar{\X}^{(2)}_0)
 {\displaystyle\mathop{\rightarrow}^{\sim}}
 H(\X_{\mathsf{s}})=\Lambda$
 compatible with the action of $\MF{A}_5$.
Hence we have
 $H(\bar{\X}^{(2)}_0)_{\MF{A}_5}\cong
 \Lambda_{\MF{A}_5}$.
Let $G$ and $\delta$ be
 as in the proof of Lemma~\ref{LEMtrans0}.
Since the action of $G$
 on $\bar{\X}^{(2)}_0$ is symplectic,
 we have
 $H(\bar{\X}^{(2)}_0)_G\cong S_{(\SYMG)}$
 by Proposition~\ref{PROPcharacterization}.
Therefore, we have
 $\Lambda_{\MF{A}_5}
 \cong(S_{(\SYMG)})_{\MF{A}_5}
 =K$,
 and the assertion (i) follows.

By Theorem~\ref{THMsymplectic},
 we have $\Lambda_{\MF{A}_5}\oplus\Z l
 \subset S(\X_{\mathsf{s}})$
 and
 $\Lambda_\PRIM^{\MF{A}_5}\supset T$.
Since we have $\rank K=18$,
 we have
 $\rank \Lambda_\PRIM^{\MF{A}_5}=3$.
On the other hand,
 by Remark~\ref{REMnontrivial},
 we have $\rank T\geq 3$.
Hence we have $\rank T=3$ and
 $T=\Lambda_\PRIM^{\MF{A}_5}$.
Let $l'\in H(\bar{\X}^{(2)}_0)$ be
 the natural quasi-polarization class
 of $\bar{\X}^{(2)}_0$.
Since we have
 $\varphi(l')=l$,
 we see that
 $\varphi(T({\bar{\X}^{(2)}_0})\oplus\Z\delta)$
 is a sublattice of $T$ of finite index.
By a direct computation,
 we can check that
 there exists
 no proper even overlattice of
 $T({\bar{\X}^{(2)}_0})\oplus\Z\delta$
 (cf.\ Proposition~\ref{PROPoverlattice}).
Therefore, we have
 $T({\bar{\X}^{(2)}_0})\oplus\Z\delta
 \cong T$.
The Grammian matrix of
 $T(\bar{\X}^{(2)}_0)\oplus\Z\delta$
 is as in the assertion (ii).
\end{proof}

\begin{rem} \label{REMmaximality}
Let $X$ be an algebraic $K3$ surface
 with a faithful (symplectic) action
 of $\MF{A}_5$.
By Theorem~\ref{THMfixedpts},
 we have
$$\rank H(X)^{\MF{A}_5}+2
 =\frac{1}{\lvert\MF{A}_5\rvert}
 \sum_{\sigma\in\MF{A}_5} \chi(\sigma)
 =6,$$
 where $\chi(\sigma)=24,8,6,4\,$
 if $\,\ord(\sigma)=1,2,3,5$.
By the same argument as in the proof
 of Lemma~\ref{LEMtransgen},
 we find that $X$ has
 Picard number $\geq 19$.
This implies that
 the moduli of algebraic $K3$ surfaces
 with faithful (symplectic) actions
 of $\mathfrak{A}_5$ is $1$-dimensional.
Hence our family $\{ \X_{\T} \}$
 is a maximal family
 of algebraic $K3$ surfaces
 with an action of $\MF{S}_5$.
\end{rem}

Let $\OP{Aut}(\X_\mathsf{s},l)$ denote
 the automorphism group of the polarized
 $K3$ surface $(\X_\mathsf{s},l)$,
 i.e., $\OP{Aut}(\X_\mathsf{s},l)$ consists of
 the automorphisms of $\X_\mathsf{s}$
 preserving the polarization $l$.
In general,
 the automorphism group of
 a polarized $K3$ surface is
 a finite group
 \cite{pss71}.

\begin{lem} \label{LEMautopol}
$\OP{Aut}(\X_\mathsf{s},l)=\MF{S}_5$.
\end{lem}
\begin{proof}
Set $\widetilde{G}
 =\OP{Aut}(\X_\mathsf{s},l)$
 and $G=\{ g\in\widetilde{G} \bigm|
 g^*\omega_{\X_\mathsf{s}}
 =\omega_{\X_\mathsf{s}} \}$.
We have $\MF{S}_5\subset\widetilde{G}$
 and $\MF{A}_5\subset G$.
By Theorem~\ref{THMsymplectic}
 and Lemma~\ref{LEMtransgen}(ii),
 we have $\Lambda^G=\Lambda^{\MF{A}_5}$
 and
 the action of $G$ on
 $A_{\Lambda_{\MF{A}_5}}$
 is trivial.
By Lemma~\ref{LEMsurjectivityA5}
 and \ref{LEMtransgen}(i),
 we have $G=\MF{A}_5$.
Let $m$ be the order
 of the cyclic group
 $Q:=\widetilde{G}/G$.
By Theorem~3.1(c) in \cite{nikulin79fin},
 all eigenvalues of the action of
 a generator of $Q$ on $T$ are
 primitive $m$-th roots of unity.
In particular,
 we have
 $\varphi(m) | \rank T $,
 where $\varphi$ is the Euler function.
By Lemma~\ref{LEMtransgen}(ii),
 we have $\rank T=3$,
 and we have
 $m=1$ or $2$.
On the other hand,
 we have $[\MF{S}_5:G]=2$.
Therefore, we have $m=2$ and
 $\widetilde{G}=\MF{S}_5$.
\end{proof}

\begin{lem} \label{LEMsurjectivity_T}
The natural map
 $\OP{Aut}_{\MF{S}_5}(\Lambda,l)
 \rightarrow\OP{O}(T)/\pm1$
 is surjective.
\end{lem}
\begin{proof}
Let $\varphi$ be an element in $\OP{O}(T)$.
By Lemma~\ref{LEMsurjectivityA5}
 and \ref{LEMtransgen},
 we have
 $\disc(\Lambda^{\MF{A}_5})=\disc(T)=-300$
 and $\BF{l}{l}=4$.
Hence we have
 $[\Lambda^{\MF{A}_5}:T\oplus\Z l]=2$.
By a direct computation,
 we can check that
 $\Lambda^{\MF{A}_5}/(T\oplus\Z l)$ is
 generated by the image of
 $(\delta+l)/2$
 (cf.\ Proposition~\ref{PROPoverlattice}).
We shall show that
 $\varphi$ can be extended to
 an element
 $\varphi_1\in\OP{O}(\Lambda^{\MF{A}_5})$
 with $\varphi_1(l)=l$.
It is sufficient to show that
 the action of $\varphi$ fixes
 $\delta/2 \bmod T\in A_T$.
In fact,
 since $\delta/2 \bmod T$ is
 a unique element of order $2$ in $A_T$,
 the action of $\varphi$ fixes
 $\delta/2 \bmod T$.
By Lemma~\ref{LEMsurjectivityA5},
 there exists an extension
 $\Phi\in\OP{O}(\Lambda)$ of $\varphi_1$
 (cf.\ Proposition~\ref{PROPembedding}).
We shall show that
 $\MF{S}'_5:=
 \Phi\cdot\MF{S}_5\cdot\Phi^{-1}=\MF{S}_5$
 in $\OP{O}(\Lambda)$.
By the proof of Lemma~\ref{LEMautopol},
 we have
 $\sigma|_T=\OP{sgn}(\sigma)\cdot\OP{id}_T$
 for any $\sigma\in\MF{S}_5$.
Hence we have
 $\sigma'|_T=\OP{sgn}(\sigma')\cdot\OP{id}_T$
 for any $\sigma'\in\MF{S}'_5$.
By Theorem~\ref{THMtorelli},
 the action of $\MF{S}'_5$ on $\Lambda$
 induces that on $\X_\mathsf{s}$
 preserving the polarization $l$.
By Lemma~\ref{LEMautopol},
 we have $\MF{S}'_5=\MF{S}_5$.
Since the outer automorphism group
 of $\MF{S}_5$ is trivial,
 there exists $\sigma\in\MF{S}_5$
 such that
 $\Phi':=\sigma\cdot\Phi$ commutes
 with $\MF{S}_5$.
We have
 $\Phi'\in\OP{Aut}_{\MF{S}_5}(\Lambda,l)$
 and $\Phi' |_T=\pm \varphi$.
\end{proof}

We determine the transcendental lattice of
 $\bar{\X}^{(\nu)}_0$ for $\nu=1,3,4$.
We prepare the following lemma.

\begin{lem} \label{LEMindex}
Let $L$ be a lattice
 and $\rho$ an element in $\OP{O}(L)$
 of order $2$.
We set
$$L_+=\left\{v\in L\bigm| \rho (v)=v \right\},\quad
 L_-=\left\{v\in L\bigm| \rho (v)=-v \right\}.$$
Then, $L/(L_+ \oplus L_-)$ is
 a $2$-elementary abelian group,
 i.e., we have
 $L/(L_+ \oplus L_-)\cong(\Z/2\Z)^n$.
Moreover, we have
 $n\leq\min \{ \rank L_+,\rank L_- \}$.
Conversely,
 for a non-degenerate primitive sublattice
 $M\subset L$
 such that $L/(M \oplus (M)^\bot_L)$ is
 $2$-elementary,
 there exists an element $\rho\in\OP{O}(L)$
 with $L_+=M$.
\end{lem}
\begin{proof}
For $v\in L$, we have a decomposition
 $v=v_+ + v_-$
 with
 $v_+\in L_+\otimes\Q$
 and $v_-\in L_-\otimes\Q$.
We have $2v_+=v+\rho v\in L$.
We define $\varphi(v \bmod{L_+ \oplus L_-})=2v_+ \bmod 2L_+$.
We can easily see that
$\varphi:L/(L_+ \oplus L_-) \rightarrow L_+/2L_+$
 is a well-defined injection.
Hence we have $L/(L_+ \oplus L_-)\cong(\Z/2\Z)^n$
 with $n\leq\rank L_+$.
Similarly, we have $n\leq\rank L_-$.

Conversely, for $M$ as in the assertion,
 any element $w\in L$ is of the form
 $w_+ + w_-$
 with $2 w_+\in M$ and $2 w_-\in(M)^\bot_L$.
Let $\rho$ be an element in $\OP{O}(L\otimes\Q)$
 with $\rho|_M=\OP{id}$
 and $\rho|_{(M)^\bot_L}=-\OP{id}$.
Then, we have
 $\rho(w)
 =w_+ - w_-
 =w-2w_-\in L$.
Hence, we have $\rho(L)=L$ and $L_+=M$.
\end{proof}

\begin{thm} \label{THMtranslist}
For $\nu=1,\ldots,5$,
 the transcendental lattice
 of $\bar{\X}^{(\nu)}_0$
 is as follows:
\begin{center}
\begin{tabular}{c|ccccc}
$\T$ & $\T_1$ & $\T_2$ & $\T_3$ & $\T_4$ & $\T_5$ \\
\hline
$X$ &${\bar{\X}^{(1)}_0}$ & ${\bar{\X}^{(2)}_0}$ & ${\bar{\X}^{(3)}_0}$ & ${\bar{\X}^{(4)}_0}$ &${\bar{\X}^{(5)}_0}$ \\
$T(X)$ &
$\begin{pmatrix} 4 & 0 \\ 0 & 10 \end{pmatrix}$ &
$\begin{pmatrix} 4 & 1 \\ 1 & 4 \end{pmatrix}$ & 
$\begin{pmatrix} 4 & 2 \\ 2 & 16 \end{pmatrix}$ &
$\begin{pmatrix} 6 & 0 \\ 0 & 20 \end{pmatrix}$ &
$\begin{pmatrix} 20 & 10 \\ 10 & 20 \end{pmatrix}$
\\
$(T(X))^\bot_T$ & $\langle -30\rangle$ & $\langle -20\rangle$ &
$\langle -20 \rangle$ & $\langle -10 \rangle$
& $\langle -4 \rangle$\rule{0ex}{2.8ex} \\
$\begin{matrix}[T:T(X) \\ \ \ \ \ \oplus (T(X))^\bot_T]\end{matrix}$
& $2$ & $1$ & $2$ & $2$ & $2$\rule{0ex}{2.5ex}
\end{tabular}\lower6.8ex\hbox{.}
\end{center}
Here we identify
 the generic transcendental lattice $T$
 with $H_\PRIM(X)^{\MF{A}_5}$ in each case
 (cf.\ Proposition~\ref{PROPdeform}
 and Lemma~\ref{LEMtransgen}(ii)).
\end{thm}
\begin{proof}
Let $l'\in H(X)$
 be the natural (quasi-)polarization class
 of $X:={\bar{\X}^{(\nu)}_0}$.
Set $m=[T:T(X)\oplus({T(X)})^\bot_T]$.

First, we consider the case $\nu=5$.
We have already determined
 the transcendental lattice of $X$
 in Lemma~\ref{LEMtransinf}.
Recall that
 $\X^{(5)}_0
 \cong\P^1\times\P^1$
 and $X={\bar{\X}^{(5)}_0}$ is
 a double covering of $\X^{(5)}_0$.
Let
 $\delta',\delta'' \in H(X)$ be
 the pull-backs
 of the classes of bi-degree $(1,0),(0,1)$
 in $H^2(\X^{(5)}_0,\Z)$, respectively.
By Lemma~\ref{LEMquadric},
 we find that $\delta',\delta''$ is
 an element in
 $H(X)^{\MF{A}_5}$.
Set $\delta=\delta'-\delta''$.
We have
 $\BF{\delta}{\delta}=-4$
 and $\delta\bot l'$.
(In particular,
 $\delta$ is a primitive element.)
Hence,
 $T(X)\oplus\Z\delta$
 is a sublattice of $T$ of finite index.
We have
 $\disc(T(X)\oplus\Z\delta)
 =-1200$.
On the other hand,
 we have $\disc(T)=-300$
 by Lemma~\ref{LEMtransgen}.
Therefore we have $m=2$.

In the remaining cases,
 we define
 $\delta\in H(X)^{\MF{A}_5}$
 as the sum of the classes
 of the (smooth) rational curves on
 $X={\bar{\X}^{(\nu)}_0}$
 which contract to
 $A_1$-singularities on $\X^{(\nu)}_0$
 (see Proposition~\ref{PROPsingpts}).
We have
 $\delta^2=-2k$
 and
 $\delta\bot l'$,
 where $k$ is the number
 of $A_1$-singularities on $\X^{(\nu)}_0$.
By the same argument as above,
 we find that
 $\delta$ is a primitive element.
Since $T(X)\oplus\Z\delta$ is a sublattice
 of $T$ of finite index,
 we have
 $\lvert \disc(T) \rvert
 =2k\cdot\disc(T(X))/m^2=300$.
By Picard-Lefschetz formula
 (cf.\ \cite{bhpv}),
 the local monodromy around $\T_\nu$ is
 the identity on $T(X)$
 and maps $\delta$ to $-\delta$.
Applying Lemma~\ref{LEMindex},
 we find that $m=1$ or $2$.

In the case $\nu=1$, we have $k=15$,
 and we have
 $\disc(T(X))=10$ or $40$.
Since $T(X)$
 is an even lattice of rank $2$,
 we have $\disc(T(X))
 \equiv 0$ or $3 \bmod 4$.
Hence, we have $m=2$
 and $\disc(T(X))=40$.
We can check that
 $T(X)
 \cong
 \begin{pmatrix} 2 & 0 \\ 0 & 20 \end{pmatrix}$
 or
 $\begin{pmatrix} 4 & 0 \\ 0 & 10 \end{pmatrix}.$
(For the classification of lattices of rank $2$
 of given discriminant,
 we refer to Chapter 15 in \cite{SP}.)
Assume that the former case occurs.
Set $T'=T(X)\oplus\Z\delta$.
Since we have $[T:T']=m=2$,
 we have $T\otimes\Z_5\cong T'\otimes\Z_5$.
However, by a direct computation,
 we have
 $(A_T)_5\cong(A_{T'})_5\cong(\Z/5\Z)^2$,
 $(q_T)_5\cong\LF{1/5}\oplus\LF{2/5}$
 and
 $(q_{T'})_5\cong\LF{1/5}\oplus\LF{1/5}$.
By Theorem~\ref{THMfinquad},
 this is a contradiction.
Hence the latter case occurs.

In the case $\nu=2$,
 the assertion follows
 by the proof of Lemma~\ref{LEMtransgen}.

In the case $\nu=4$,
 we have $k=5$.
By the same argument
 as in the case $\nu=1$,
 we can check that
 $T(X)$ is isomorphic
 to one of the following lattices:
$$\begin{array}{c|c|c|c}
 m & Q & (q_{T'})_3 & (q_{T'})_5
   \\ \hline
 2 & \begin{pmatrix} 2 & 0 \\ 0 & 60 \end{pmatrix}
   & \LF{2/3}
   & \LF{1/5}\oplus\LF{1/5}
   \\
 2 & \begin{pmatrix} 4 & 0 \\ 0 & 30 \end{pmatrix}
   & \LF{1/3}
   & \LF{1/5}\oplus\LF{2/5}
   \\
 2 & \begin{pmatrix} 6 & 0 \\ 0 & 20 \end{pmatrix}
   & \LF{2/3}
   & \LF{1/5}\oplus\LF{2/5}
   \\
 2 & \begin{pmatrix} 10 & 0 \\ 0 &  12 \end{pmatrix}
   & \LF{1/3}
   & \LF{1/5}\oplus\LF{1/5}
\end{array}.$$
Here we set $T'=Q\oplus\Z\delta$.
Since we have
 $(q_T)_3\cong\LF{2/3}$ and
 $(q_T)_5\cong\LF{1/5}\oplus\LF{2/5}$,
 we find that
 $T(X)\cong
 \begin{pmatrix} 6 & 0 \\ 0 & 20 \end{pmatrix}$.

In the case $\nu=3$,
 we have $k=10$.
By the above argument
 and Lemma~\ref{LEM42216} below,
 we have
 $T(X)
 \cong
 \begin{pmatrix} 4 & 2 \\ 2 & 16 \end{pmatrix}$.
Since we have $\disc(T(X))=60$,
 we have $m=2$.
\end{proof}

\begin{lem} \label{LEM42216}
We have
$T(\bar{\X}^{(\nu)}_0)\cong
 \begin{pmatrix} 4 & 2 \\ 2 & 16 \end{pmatrix}$
 for some $\nu\in \{ 1,\ldots,5 \}$.
\end{lem}
\begin{proof}
Let $e_1,e_2,e_3$ be a basis of $T$
 for which the Grammian matrix of $T$
 is equal to (\ref{gramtransgen}).
The sublattice $Q$ of $T$ generated by
 $2e_1-2e_2+e_3,-3e_1-e_3$
 is primitive and isomorphic to
 $\begin{pmatrix} 4 & 2 \\ 2 & 16 \end{pmatrix}$.
By a direct computation, we can check that
 $(Q)^\bot_T$ is generated by
 $\delta:=6e_1-4e_2+3e_3$.
Let $\varphi$ be an element
 in $\OP{O}(T\otimes\Q)$
 such that $\varphi|_Q=\OP{id}_Q$ and
 $\varphi(\delta)=-\delta$.
By a direct computation, we have
 $\varphi(T)=T$, i.e.\ $\varphi\in\OP{O}(T)$.
By Lemma~\ref{LEMsurjectivity_T},
 there exists an extension
 $\Phi\in\OP{Aut}_{\MF{S}_5}(\Lambda,l)$
 of $\varphi$ or $-\varphi$.
Let $\omega\neq 0$ be an element
 in $Q\otimes\C$
 with $\BF{\omega}{\omega}=0$.
We see that $\C\omega\in\Omega_T$
 is fixed under the action of $\Phi$.
By Remark~\ref{REMnofixedpoint}
 and Proposition~\ref{PROPpisbiholomorphic},
 we have
 $\bar{\PER}'{}^{-1}([\C\omega])=\T_\nu$
 for some $\nu\in \{ 1,\ldots,5 \}$.
By the definition of $\bar{\PER}'$,
 we have
 $T(\bar{\X}^{(\nu)}_0)\cong Q$.
\end{proof}

\section{Quaternion algebras}
\label{SECTquaternion}
\subsection{Short recall of quaternion algebras and Shimura curves}
\subsubsection{Properties on quaternion algebras and orders}
In this paragraph, we briefly recall several properties on
quaternion algebras.
For details, see \cite{vigneras80}.

Let $B$ be a quaternion algebra over a field $K$,
i.e., $B$ is a central simple algebra over $K$
with $\dim_{K} B=4$.
We denote by $B^\times$ the group of invertible elements in $B$.
The reduced trace and the reduced norm of $x\in B$ are defined by
$\operatorname{Tr}(x)=x+\bar{x}$,
$\operatorname{Nr}(x)=x\bar{x}\in K$, respectively,
where $\bar{\ }$ denotes the canonical involution.

Let $R$ be $\Z$ or $\Z_p$,
$K$ the quotient field of $R$,
and $B$ a quaternion algebra over $K$.
A subring $\MC{O}$ of $B$
which is a free $R$-module
and generates $B$ over $K$
is called
an order of $B$ over $R$.
We have
$\operatorname{Tr}(x),\operatorname{Nr}(x)\in R$
for all $x\in\MC{O}$.
Any order $\MC{O}$ is contained in a maximal order.
If $\MC{O}$ is the intersection of two maximal orders,
$\MC{O}$ is said to be an Eichler order.
In particular, a maximal order is an Eichler order.
Let $\MC{O}^\vee$ be the dual of $\MC{O}$
with respect to the bilinear form $(x,y)\mapsto \operatorname{Tr}(xy)$,
i.e., $\MC{O}^\vee
=\{ x\in B\bigm| \operatorname{Tr}(x\MC{O})\subset R \}$.
We define the different $\MC{D}_\MC{O}$ of $\MC{O}$
as the inverse of $\MC{O}^\vee$, i.e.,
$\MC{D}_\MC{O}=(\MC{O}^\vee)^{-1}
=\{ x\in B \bigm|
\MC{O}^\vee x \MC{O}^\vee \subset \MC{O}^\vee \}$.
The reduced discriminant $d(\MC{O})$ of $\MC{O}$ is
defined as
the fractional ideal of $R$ generated by the reduced norms of
elements in $\MC{D}_\MC{O}$.
Let $u_1,\ldots,u_4$ be a basis of $\MC{O}$ over $R$.
We have the following:
\begin{equation}
 R\det( (\operatorname{Tr}(u_i u_j))_{ij} )
 =d(\MC{O})^2.
\end{equation}
If $\MC{O}'\subset\MC{O}$ is an order,
we have $[\MC{O}:\MC{O}']=[d(\MC{O}):d(\MC{O}')]$.
In particular, $\MC{O}$ is a maximal order if $d(\MC{O})=R$.
We identify $d(\MC{O})$
with its positive generator in the case $R=\Z$.
We set
\begin{equation}
 \MC{O}^1
 =\{ x\in\MC{O} \bigm| \OP{Nr}(x)=1 \},\quad
 N(\MC{O})
 =\{x\in B^\times\bigm|x\MC{O}x^{-1}=\MC{O}\}.
\end{equation}
We use the following theorem.

\begin{thm}[cf.~\cite{serre77,vigneras80}]
 \label{THMlocalquat}
Let $B$ be a quaternion algebra over $\Q_p$ and
 $\MC{O}$ be an Eichler order of $B$ over $\Z_p$.
Then we have the following.
\begin{itemize2}
\item[(i)]
The quaternion algebra $B$ is isomorphic
to either the quaternion division algebra $H$
or the matrix algebra $M(2,\Q_p)$.
The quaternion division algebra
$H$ is given by
$$H=L\oplus Lu,\quad
u^2=p,\, ux=x^\ast u~(\forall x\in L),$$
where $L/\Q_p$ is the unique unramified quadratic field
and $x^\ast$ is the conjugate of $x$ in $L/\Q_p$.
\item[(ii)]
The quaternion division algebra $H$
has a unique maximal order $\MC{O}=R_L[u]$,
where $R_L$ is the ring of integers of $L$
 and $u$ is as in (i).
Moreover, we have $d(\MC{O})=p\Z_p$.
In particular, the quaternion division algebra $H$
has a unique Eichler order $\MC{O}$.
\item[(iii)]
If $B=M(2,\Q_p)$,
 then $\MC{O}$ is conjugate to
 $\MC{O}^{(p)}(n)$ in $B$ for some $n\geq 0$,
 where 
$$\MC{O}^{(p)}(n)=
 \begin{pmatrix}
  \Z_p & \Z_p \\
  \Z_p & \Z_p
 \end{pmatrix}
 \cap
 \begin{pmatrix}
  \Z_p & p^{-n}\Z_p \\
  p^n\Z_p & \Z_p
 \end{pmatrix}
 =
 \begin{pmatrix}
 \Z_p & \Z_p \\
 p^n\Z_p & \Z_p
 \end{pmatrix}.$$
Moreover, we have $d(\MC{O}^{(p)}(n))=p^n\Z_p$.
\item[(iv)]
The quotient group
$N(\MC{O})/\Q_p^\times\MC{O}^\times$
is
trivial if and only if
$B\cong M(2,\Q_p)$
and $\MC{O}$ is maximal.
Otherwise, $N(\MC{O})/\Q_p^\times\MC{O}^\times$ is
isomorphic to $\Z/2\Z$ and its generator is as follows:
$$\begin{array}{c|cc}
 \MC{O} & R_L[u] \text{ in (ii)  }
        & \MC{O}^{(p)}(n) \text{ in (iii)} \\
 \hline
 \text{generator is represented by} & u
                  & \begin{pmatrix}0&1\\p^n&0\end{pmatrix}
\end{array}.$$
\item[(v)]
 The normalizer of $\MC{O}^1$ in $B^\times$
 coincides with $N(\MC{O})$.
\end{itemize2}
\end{thm}

The following is a conclusion of Theorem~\ref{THMlocalquat}
and the fact that any subalgebra of $M(2,\mathbb{F}_p)$ of index $p$ is
conjugate to
$\begin{pmatrix} \mathbb{F}_p & \mathbb{F}_p \\ 0 & \mathbb{F}_p \end{pmatrix}$
in $M(2,\mathbb{F}_p)$.

\begin{thm}[cf.~\cite{alsinabayer}] \label{THMindexp}
Any suborder $\MC{O}\subset M(2,\Z_p)$ (over $\Z_p$)
of index $p$ is an Eichler order.
\end{thm}

Let $B$ be a quaternion algebra over $\Q$ and
$\MC{O}$ an order of $B$ over $\Z$.
The reduced discriminant $d(B)$ of $B$ is defined
as the product of the prime numbers $p$
such that $B\otimes\Q_p$ is a division algebra.
Quaternion algebras over $\Q$
of the same reduced discriminant
are isomorphic to each other
 (cf.\ \cite{vigneras80}).
By Theorem~\ref{THMlocalquat},
we have $d(\MC{O})=d(B)N$
 for some positive integer $N$.
(Note that we have
$d(\MC{O}\otimes\Z_p)=d(\MC{O})\Z_p$.)
If $\MC{O}$ is an Eichler order, $d(B)$ and $N$ are coprime.
The following fact is useful to determine $d(B)$.

\begin{thm}[cf.~\cite{vigneras80}] \label{THMramification}
Let ${\displaystyle \left( \frac{a,b}{\Q} \right)}$
 be a quaternion algebra generated by $1,i,j,k$
 as a $\Q$-vector space
 with relations $i^2=a,j^2=b,ij=-ji=k$
 for $a,b\in\Q^\times$.
Let $v$ be a place of $\Q$,
i.e., $v$ is a prime number or $\infty$.
The following conditions are equivalent:
\begin{itemize2}
\item[(i)]
${\displaystyle \left( \frac{a,b}{\Q} \right)} \otimes \Q_v
\cong M(2,\Q_v)$;
\item[(ii)]
$(a,b)_v=1$.
\end{itemize2}
Here $\Q_\infty=\R$
and $(\ ,\ )_v$ is the Hilbert symbol.
\end{thm}

\subsubsection{Fuchsian groups %
$\Gamma^1(\MC{O})$ and $\Gamma^+(\MC{O})$}
\label{PARAfuchsian}
In this paragraph,
 we recall Fuchsian groups induced by
 Eichler orders.
Let $B$ be a quaternion algebra
 over $\Q$
 and $\MC{O}$ an order of $B$
 over $\Z$.
Assume that $B$
 is indefinite, i.e.,
 $B\otimes\R\cong M(2,\R)$.
We set
\begin{equation}
\begin{array}{c}
 \Gamma^1(\MC{O})
 =\MC{O}^1/\pm 1,
 \vspace{1ex}\\
 N^+(\MC{O})=\{x\in N(\MC{O}) \bigm|
 \OP{Nr}(x)>0 \} ,\quad
 \Gamma^+(\MC{O})=N^+(\MC{O})/\Q^\times.
\end{array}
\end{equation}
The groups $\Gamma^1(\MC{O})$ and $\Gamma^+(\MC{O})$
 become Fuchsian groups
 by an isomorphism $B\otimes\R\cong M(2,\R)$.
Also, these groups
 act on
the following bounded symmetric domain
 $\MC{D}(\MC{O},I)$ of type IV,
 which is isomorphic to
 the upper half-plane $\MF{H}$.
We set
\begin{gather}
\MC{I}(\MC{O})=\{ x\in\MC{O} \bigm| \OP{Tr}(x)=0 \}.
\end{gather}
Let $Q_{\MC{I}(\MC{O})}$ be the quadratic form on
$\MC{I}(\MC{O})$
defined by
\begin{equation}
 Q_{\MC{I}(\MC{O})}:\MC{I}(\MC{O})\longrightarrow\Z;\quad
 x\longmapsto x^2.
\end{equation}
We also denote the bilinear form on $\MC{I}(\MC{O})$
 corresponding to this quadratic form
 by $Q_{\MC{I}(\MC{O})}$:
\begin{equation}
 Q_{\MC{I}(\MC{O})}(x,y)
 =\frac{1}{2}( (x+y)^2-x^2-y^2 )
 =\frac{1}{2}\OP{Tr}(xy),\quad
 x,y\in\MC{I}(\MC{O}).
\end{equation}
Then, we obtain a lattice
 $(\MC{I}(\MC{O}),Q_{\MC{I}(\MC{O})})$
of signature $(2,1)$.
The quadratic form $Q_{\MC{I}(\MC{O})\otimes\R}$
 on $\MC{I}(\MC{O})\otimes\R$
is defined as $Q_{\MC{I}(\MC{O})}\otimes\R$.
The group
 $\Gamma^1(\MC{O}\otimes\R)$
 acts on
 $(\MC{I}(\MC{O})\otimes\R,Q_{\MC{I}(\MC{O})\otimes\R})$
 by
\begin{equation}
 \pm\gamma\cdot x:=\gamma x \gamma^{-1},\quad
 \pm\gamma\in\Gamma^1(\MC{O}\otimes\R),~
 x\in\MC{I}(\MC{O})\otimes\R.
\end{equation}
Similarly, the groups $\Gamma^1(\MC{O})$
 and $\Gamma^+(\MC{O})$
 acts on the lattice
 $(\MC{I}(\MC{O}),Q_{\MC{I}(\MC{O})})$.
Let $I$ be an element in $\MC{I}(\MC{O})\otimes\R$ with $I^2<0$.
We set
\begin{equation} \label{def_dom}
 \MC{D}(\MC{O},I)
 =\{ x\in\MC{I}(\MC{O})\otimes\R \bigm|
 x^2=-1,~Q_{\MC{I}(\MC{O})\otimes\R}(x,I)<0 \}.
\end{equation}
This domain is isomorphic to $\MC{D}_T^\circ$ defined
 in Subsection~\ref{SUBSECTtransA}.
Since $\Gamma^1(\MC{O}\otimes\R)\cong\OP{PSL}(2,\R)$
 is connected,
 the action of $\Gamma^1(\MC{O}\otimes\R)$
 preserves $\MC{D}(\MC{O},I)$.
Hence, the groups
 $\Gamma^1(\MC{O})$ and $\Gamma^+(\MC{O})$
 act on $\MC{D}(\MC{O},I)$.
We have $\operatorname{vol}(
 \Gamma^1(\MC{O})\backslash\MC{D}(\MC{O},I))<\infty$.
For an Eichler order $\MC{O}$,
 the compactification of 
 $\Gamma^1(\MC{O})\backslash\MC{D}(\MC{O},I)$
 is the $\C$-valued points of
 the Shimura curve associated to $\MC{O}$,
 which is defined over $\Q$.

\subsection{A dual Clifford algebra for the generic %
 transcendental lattice}
\subsubsection{Clifford algebras and quadratic forms on
 pure quaternions}
We recall some properties on Clifford algebras.
Let $V$ be a non-degenerate lattice
equipped with a bilinear form
$\langle\ ,\ \rangle$.
Let $S$ be the tensor algebra of $V$ over $\Z$,
$J$ the two-sided ideal of $S$ generated by elements of the form
$v\otimes v- \langle v,v \rangle$ ($v\in V$).
The Clifford algebra $\operatorname{Cl}(V)$ of $V$
 is defined
 as the quotient ring $S/J$.
Then $V$ is considered as a submodule of $\OP{Cl}(V)$.
The even part $\operatorname{Cl}^+(V)$
of $\operatorname{Cl}(V)$
 is defined as the subring of $\OP{Cl}(V)$
generated by elements of the form $vw$ ($v,w\in V$) over $\Z$.
We have $\rank  \operatorname{Cl}(V)=2^n$
and $\rank \OP{Cl}^+(V)=2^{n-1}$,
where $n=\OP{rank}V$.
We assume that $\rank V=3$.
Then $\OP{Cl}^+(V)\otimes\Q$
 is a quaternion algebra
over $\Q$
and $\OP{Cl}^+(V)$ is an order of
 $\operatorname{Cl}^+(V)\otimes\Q$
over $\Z$.
We have the following canonical isomorphism
$\MC{I}(\OP{Cl}^+(V))
{\displaystyle\mathop{\rightarrow}^{\sim}}
V^\vee[\OP{disc}(V)]$ up to sign.
Let $\OP{Cl}^{(2)}(V)$ denote the $\Z$-submodule of $\OP{Cl}(V)$
consisting of elements of degree at most 2.
Then we have $\OP{Cl}^{(2)}(V)=\OP{Cl}^+(V)\oplus V$.
Recall the canonical isomorphism
 $\OP{Cl}(V)/\OP{Cl}^{(2)}(V)\cong
 \bigwedge^3 V$.
We fix an isomorphism
\begin{equation}
 \alpha:\OP{Cl}(V)/\OP{Cl}^{(2)}(V)
 \CONGL\Z.
\end{equation}
We define a pairing $(~,~)_V$ by
\begin{equation} \label{defpairing}
(~,~)_V:V\times\MC{I}(\OP{Cl}^+(V))\longrightarrow\Z;\quad
(v,x)\longmapsto \alpha(vx \bmod \OP{Cl}^{(2)}(V)).
\end{equation}
Under the pairing $(~,~)_V$,
$V$ and $\MC{I}(\OP{Cl}^+(V))$ are dual to each other over $\Z$.
The pairing $(~,~)_V$ induces an isomorphism $\eta$
between $\Z$-modules:
\begin{equation} \label{defeta}
\eta:\MC{I}(\OP{Cl}^+(V))
\CONGL V^\vee;\quad
x\longmapsto v,~(w,x)_V=\langle w,v \rangle~(\forall w\in V).
\end{equation}
Here we consider $V^\vee$ as a submodule of $V\otimes\Q$, as usual.
We have $\eta(x)^2 =\OP{Nr}(x)/\OP{disc}(V)$.
Hence $\eta$ is an isomorphism
between $( \MC{I}(\OP{Cl}^+(V)), Q_{\MC{I}(\OP{Cl}^+(V))})$
and $V^\vee[\OP{disc}(V)]$.

\subsubsection{Quaternion orders $\ORDa$ %
and $\ORDb$}
We apply the above results to study the automorphism group
 $\OP{O}(T)$ of
 the generic transcendental lattice $T$
 of the family $\{\X_\T \}$.
By Lemma~\ref{LEMtransgen}, we have
$$T\cong
\begin{pmatrix}
4 & 1 & 0 \\
1 & 4 & 0 \\
0 & 0 & -20
\end{pmatrix}.$$
We set
\begin{equation}
T^*=T^\vee[-60]\cong
\begin{pmatrix}
-16 & 4 & 0 \\
4 & -16 & 0 \\
0 & 0 & 3
\end{pmatrix}.\end{equation}
We have the canonical identification
 ${T^*}=T[-1/60]^{\vee}$.
In the case $V=T^*$,
 the map $\eta$ defined by (\ref{defeta}) is as follows:
\begin{equation} \label{isomT12}
 \eta:(\MC{I}(\ORDa),Q_{\MC{I}(\ORDa)})
 \mathop{\longrightarrow}^\sim
 {T^*}^\vee[\disc(T^*)]=T[12].
\end{equation}
We set
\begin{equation} \label{defbasis}
e_1=v_2v_3,~e_2=v_3v_1,~e_3=v_1v_2-4\in\ORDa.
\end{equation}
Then $1,e_1,e_2,e_3$ form an integral basis of $\ORDa$
and $e_1,e_2,e_3$ form an integral basis of $\MC{I}(\ORDa)$.
By a direct computation, we have
\begin{equation} \label{structureconst}
\begin{array}{c}
e_1^2=48,\ e_2^2=48,\ e_3^2=-240,\rule[-1.5ex]{0ex}{0ex} \\
e_1 e_2=-3e_3+12,\ e_2 e_3=16e_1-4e_2,\ e_3 e_1=-4e_1+16e_2.
\end{array}
\end{equation}
We have an isomorphism between quaternion algebras over $\Q$:
$$\ORDa\otimes\Q \cong \left( \frac{3,5}{\Q} \right);\quad
\frac{e_1}{4}\leftrightarrow i,\ \frac{e_1-4e_2}{12}
\leftrightarrow j.$$
By Theorem~\ref{THMramification},
we have
\begin{equation} \label{reduced_disc15}
 d(\ORDa\otimes\Q)=3\cdot5=15
\end{equation}
 and $\ORDa\otimes\Q$ is indefinite.
(For formulae to compute the Hilbert symbol,
 see e.g.\ \cite{serre70}.)
The Grammian matrix of $\MC{I}$ under $e_1,e_2,e_3$ is
$$-12\cdot\begin{pmatrix}
4 & 1 & 0 \\
1 & 4 & 0 \\
0 & 0 & -20
\end{pmatrix}.$$
We define $w_1,w_2,w_3\in\ORDa\otimes\Q$ by
\begin{equation} \label{defw}
w_1=\frac{e_1-e_2}{6},~w_2=\frac{e_2}{2},~
w_3=\frac{1}{2}\left(1+\frac{e_3}{4}\right).
\end{equation}
We consider an extension $\ORDb$ of $\ORDa$:
\begin{equation} \label{defordb}
\ORDb=
\Z\oplus\Z w_1 \oplus \Z w_2 \oplus \Z w_3.
\end{equation}
In fact, we can check that
 $\ORDb$ is an order by a direct computation.

\begin{lem} \label{LEMordext}
We have the following:
\begin{itemize2}
\item[(i)]
$\ORDb$ is an Eichler order
 of reduced discriminant $30$
 and of level $2$;
\item[(ii)]
$N^+(\ORDa)=N^+(\ORDb)$.
\end{itemize2}
\end{lem}
\begin{proof}
By a direct computation,
 we have $d(\ORDa)=2^6 \cdot 3^2 \cdot 5$ and
$[\ORDb:\ORDa]=6\cdot 2\cdot 8=2^5\cdot 3$.
Hence we have
 $d(\ORDb) 
  =30$.
By Theorem~\ref{THMindexp}
 and (\ref{reduced_disc15}),
$\ORDb\otimes\Z_2$ is an Eichler order.
For $p=3,5$,
by Theorem~\ref{THMlocalquat},
$\ORDb\otimes\Z_p$ is the maximal order
because ${\ORDa}\otimes\Q_p$ is a division algebra
 by (\ref{reduced_disc15}).
Hence $\ORDb$ is an Eichler order by the Hasse principal for
quaternion orders.
The assertion (i) follows.

We consider the sequence
$\ORDa\subset\MC{O}_1\subset\MC{O}_2\subset\MC{O}_3
\subset\ORDb$ of orders, where $\MC{O}_i$ is defined as follows:
$$
\MC{O}_1=\Z \{ 1,e_1,e_2,\frac{e_3}{2} \},~
\MC{O}_2=\Z \{ 1,\frac{e_1}{2},\frac{e_2}{2},\frac{e_3}{4} \},~
\MC{O}_3=\MC{O}_2+\Z\frac{e_1-e_2}{6}.
$$
In fact, the $\MC{O}_i$ are orders by a direct computation.
Let $\MC{I}_i=\MC{I}(\MC{O}_i)$ be the pure quaternions of $\MC{O}_i$.
It is sufficient to show that
 $\varphi(\ORDa)=\ORDa
 \Leftrightarrow \varphi(\ORDb)=\ORDb$
 for all $\varphi\in\OP{Aut}(\ORDa\otimes\Q)$.

First, we shall show that $\varphi(\ORDa)=\ORDa\Rightarrow
\varphi(\MC{O}_1)=\MC{O}_1$.
Let $\MC{O}'\subset {\ORDa_\Q}$ be an order
 which contains $\ORDa$ with index $2$.
The image of $2\MC{O}'$ in the ring $\ORDa/2\ORDa$
is a two-sided $\MC{O}'$-module of order $2$
(and is also an ideal of $\ORDa/2\ORDa$)
generated by an element, say,
$\omega=\xi_0+\xi_1 \bar{e}_1+\xi_2 \bar{e}_2+\xi_3 \bar{e}_3$.
By (\ref{structureconst}), we have the following relations:
$$\bar{e}_1^2=\bar{e}_2^2=\bar{e}_3^2=\bar{e}_2\bar{e}_3
=\bar{e}_3\bar{e}_1=\bar{0},~\bar{e}_1\bar{e}_2=\bar{e}_3.$$
(Note that $\ORDa/2\ORDa$ is a commutative ring.)
Since $(\bar{1}+\bar{e}_1)^2=(\bar{1}+\bar{e}_2)^2=\bar{1}$,
 we have
$\omega=(\bar{1}+\bar{e}_1)\omega=(\bar{1}+\bar{e}_2)\omega$,
 i.e.,
$$\bar{e}_1\omega=\xi_0\bar{e}_1+\xi_2 \bar{e}_3=0,~
\bar{e}_2\omega=\xi_0\bar{e}_2+\xi_1 \bar{e}_3=0.$$
Hence we have $\xi_0=\xi_1=\xi_2=0$ and $\omega=\bar{e}_3$,
i.e., $\MC{O}'=\MC{O}_1$.
This implies that $\MC{O}_1$ is a unique order which contains
$\ORDa$ with index $2$.
Therefore $\varphi(\ORDa)=\ORDa\Rightarrow
\varphi(\MC{O}_1)=\MC{O}_1$.

Next, we shall show that $\varphi(\MC{O}_1)=\MC{O}_1\Rightarrow
\varphi(\ORDa)=\ORDa$.
We have decompositions $\ORDa=\Z\oplus\MC{I}(\ORDa)$
and $\MC{O}_1=\Z\oplus\MC{I}_1$.
We set $L=\MC{I}_1[-1/12]$, which is an (integral) lattice.
Let $q:L/2L\rightarrow\Z/4\Z$ be a natural finite
quadratic form.
The Grammian matrices of $L$ and $q$ are
$$\begin{pmatrix}
4 & 1 & 0 \\
1 & 4 & 0 \\
0 & 0 & -5
\end{pmatrix},~
\begin{pmatrix}
0 & 1 & 0 \\
1 & 0 & 0 \\
0 & 0 & -1
\end{pmatrix},$$
respectively.
Since $e_1,e_2,e_3$ generate $\MC{I}(\ORDa)$,
 the image of $\MC{I}(\ORDa)$ in $L/2L$ is generated by
 $\{ \bar{x}\in L/2L\bigm|q(\bar{x})=0 \}$.
This property characterize $\MC{I}(\ORDa)$ uniquely.
(Note that $2\MC{I}_1\subset\MC{I}(\ORDa)$.)
Hence $\varphi(\MC{O}_1)=\MC{O}_1\Rightarrow
\varphi(\ORDa)=\ORDa$.

By the above argument, we have $\varphi(\ORDa)=\ORDa\Leftrightarrow
\varphi(\MC{O}_1)=\MC{O}_1$.
Since $\MC{I}_1=2\MC{I}_2$, it is clear that
$\varphi(\MC{O}_1)=\MC{O}_1\Leftrightarrow
\varphi(\MC{O}_2)=\MC{O}_2$.
Similarly, we have
$\varphi(\MC{O}_2)=\MC{O}_2\Leftrightarrow
\varphi(\MC{O}_3)=\MC{O}_3\Leftrightarrow
\varphi(\ORDb)=\ORDb$.
%
\end{proof}

\subsubsection{A principal symplectic form $\Psi_I$ %
 on $\ORDb$}
In the rest of the paper,
 we fix an element
 $I\in\MC{I}(\ORDb)\otimes\Q$
 with $I^2<0$
 as follows:
\begin{equation} \label{defI}
I=\frac{1}{60}(e_1+e_2+e_3).
\end{equation}
We define the non-degenerate symplectic form $\Psi_I$
on $\ORDa\otimes\R$ by
\begin{equation}
 \Psi_I(v,w)=\OP{Tr}(vI\bar{w}),\quad
 v,w\in\ORDa\otimes\R.
\end{equation}
In fact, we have
$$\Psi_I(w,v)=\OP{Tr}(wI\bar{v})
\setbox0=\hbox{$I$}
=\OP{Tr}(
\bar{\raisebox{0ex}[\ht0][0ex]{$w$}}
\hspace{-0.1ex}
\bar{\raisebox{0ex}[\ht0][0ex]{}}
\hspace{ 0.1ex}
\bar{\raisebox{0ex}[\ht0][0ex]{}}
\bar{\raisebox{0ex}[\ht0][0ex]{$I$}}
\bar{\raisebox{0ex}[\ht0][0ex]{}}
\bar{\raisebox{0ex}[\ht0][0ex]{$\bar{v}$}}
)
=\OP{Tr}(v\bar{I}\bar{w})=-\OP{Tr}(vI\bar{w})=-\Psi_I(v,w).$$
The symplectic form $\Psi_I$ is
integral and principal on $\ORDb$.
In fact, the following integral basis of $\ORDb$ is a
symplectic basis of $(\ORDb,\Psi_I)$:
\begin{equation} \label{alphabeta}
 \alpha_1=1,~\alpha_2=1+w_1,~\beta_1=w_1+w_3,~\beta_2=w_1+w_2+w_3,
\end{equation}
i.e., we have
$\Psi_I(\alpha_i,\alpha_j)=\Psi_I(\beta_i,\beta_j)=0$,
$\Psi_I(\alpha_i,\beta_j)=\delta_{ij}$.
This basis is used in Section~\ref{SECTmodular} to
construct a modular embedding.

\begin{rem}
For any indefinite Eichler order $\MC{O}$ over $\Z$,
 there exists an element $I\in\MC{O}\otimes\Q$
 such that $\Psi_I$ is integral and principal on $\MC{O}$
 (cf.\ \cite{hashimoto95}).
\end{rem}

\subsection{Several Fuchsian groups and their relations}
In this section, we define some Fuchsian groups
and show their relations.
They will appear in the argument on the modular embedding
 in Section~\ref{SECTmodular}.
In Subsection~\ref{SUBSECTtransA},
 we defined the Fuchsian group $\OP{O}(T)^\circ$,
 which acts on $\MC{D}_T^\circ\cong\MF{H}$.
Recall the isomorphism
$$\eta:(\MC{I}(\ORDa),Q_{\MC{I}(\ORDa)})
 \CONGL T[12]$$
 defined by (\ref{isomT12}).
For $g\in\GAMd$,
 we identify $g$ with its action
 on $\MC{I}(\ORDa)$ (see Paragraph~\ref{PARAfuchsian}).
By Lemma~\ref{LEMgammaisom} below,
 under the following correspondences,
 we identify $\MC{D}_T^\circ$ and $\OP{O}(T)^\circ$
 with $\DOMb$ and $\GAMd$,
 respectively:
\begin{equation} \label{identification_gam}
 \begin{array}{c}
 \MC{D}_T^\circ \longleftrightarrow \DOMb;\quad
 \sqrt{12}\cdot(\eta\otimes\R)(v) \longleftrightarrow v,
 \vspace{1ex} \\
 \OP{O}(T)^\circ=\OP{O}(T{[}12{]})^\circ
 \longleftrightarrow \GAMd;\quad
 \eta\circ g\circ\eta^{-1} \longleftrightarrow g.
 \end{array}
\end{equation}
Here we chose
 an appropriate connected component $\MC{D}_T^\circ$
 of $\MC{D}_T$.
By Lemma~\ref{LEMelliptic_ele} below,
 there exist
 exactly five conjugacy classes of elliptic elements
 in $\GAMd$, each of which
 corresponds to a singular fiber
 of the family $\{ \X_\T \}$.
They are represented
 by $g_1,\ldots,g_5\in\GAMd$,
 which are given in Lemma~\ref{LEMelliptic_ele}.

\begin{defn} \label{DEFgams}
We define two normal subgroups
 of $\GAMd$ by
\begin{gather} \label{def_gams}
 \begin{array}{c}
 \GAMa
 :=\{ \gamma\in\ORDb^1 \bigm|
 \gamma \equiv 1 \bmod 2\,\ORDb~
 \} /\pm1, \vspace{1ex} \\
 \GAMc:=\langle \GAMb,g_2,g_4 \rangle.
 \end{array}
\end{gather}
\end{defn}
We summarize relations of
 $\GAMd$, $\GAMc$, $\GAMb$ and $\GAMa$
 in the following table (see
 Lemmas~\ref{LEMtable1}--\ref{LEMtable2}).
The genus of each curve $G\backslash\DOMb$ is
 determined by the Riemann--Hurwitz formula.
{\scriptsize
\begin{equation} \label{table_gamma}
\begin{array}{c|cc|ccccc|c}
G & G/\GAMa & G/\GAMb
    & [\xi_1] & [\xi_2] & [\xi_3] & [\xi_4] & [\xi_5]
    & \text{genus}
 \\ \hline
\GAMd & D_8\times C_2\times C_2
    & C_2\times C_2\times C_2 &1&1&1&1&1& 0 \\
\GAMc & D_8\times C_2 & C_2\times C_2
    &1&2&2&2&1& 0 \\
\GAMb & C_2\times C_2 & \{ 1 \}
    &4&4&4&4&4& 3 \\
\GAMa & \{1\} &
    &16&16&16&16&16& 9
\end{array}
\end{equation}
}
Here 
 $\xi_\nu\in\DOMb$ is the elliptic point
 corresponding to $g_\nu$,
 and
 the right side is the table of
 the cardinality $n(G,[\xi_\nu])$
 of the inverse image of
 each
 $[\xi_\nu]:=\GAMd\cdot\xi_\nu$
 under the covering
 $G\backslash\DOMb
 \rightarrow\GAMd\backslash\DOMb$,
 which is given by
$$n(G,[\xi_\nu])=
 \begin{cases}
  [\GAMd:G]   & \text{if~~} g_\nu\not\in G \\
  [\GAMd:G]/2 & \text{if~~} g_\nu\in G.
 \end{cases}$$

\begin{lem} \label{LEMgammaisom}
The map $\GAMd\ni g\mapsto
 \eta\circ g\circ\eta^{-1}\in\OP{O}(T[12])^\circ$
 is an isomorphism.
\end{lem}
\begin{proof}
By Lemma~\ref{LEMordext}(ii),
 it is sufficient to show that the map
$$\Phi:N(\ORDa)/\Q^\times\ni g\mapsto
 \eta\circ g\circ\eta^{-1}\in\OP{SO}(T[12])
$$
 is an isomorphism.
The injectivity of $\Phi$ is trivial.
We shall show the surjectivity of $\Phi$.
Let $f$ be an element in $\OP{SO}(T[12])$.
Recall that we have the canonical identification
 ${T^*}^\vee[\disc(T^*)]=T[12]$.
We consider $f$
 as an element in $\OP{SO}({T^*}^\vee)$.
The map $f|_{T^*}$ induces
 an automorphism $F$ of $\OP{Cl}({T^*})$.
We have a commutative diagram
$$\begin{CD}
\MC{I}(\ORDa) @>{F}>> \MC{I}(\ORDa) \\
@V{\eta}VV             @VV{\eta}V \\
{T^*}^\vee @>{f}>> {T^*}^\vee.
\end{CD}$$
In fact,
 if $x\in\MC{I}(\ORDa)$ and $v=\eta(x)$,
 we have
$$\langle f(w),f(v)\rangle
 =\BF{w}{v}=(w,x)_{T^*}
 =\det(f)^{-1}(f(w),F(x))_{T^*}
 =(f(w),F(x))_{T^*}$$
 for all $w\in {T^*}$.
Hence, we have
 $\eta(F(x))=f(v)=f(\eta(x))$.
By the Skolem--Noether theorem,
 the restriction of $F$ to
 $\ORDa$ is of the form
 $x\mapsto \gamma x \gamma^{-1}$ for some
 $\gamma\in N(\ORDa)$,
 i.e.,
 we have $F|_{\ORDa}\in N(\ORDa)/\Q^\times$.
As a consequence, we have
 $\Theta(F|_{\ORDa})=f$.
The surjectivity of $\Phi$ follows.
\end{proof}

\begin{lem} \label{LEMelliptic_ele}
There exist exactly five conjugacy classes of
 elliptic points on $\DOMb$
 under the action of $\GAMd$,
 each of which is of order $2$.
They are represented by the following points:
\begin{align}
\xi_1:=&\frac{1}{\sqrt{12\cdot 30}}(-5e_1+5e_2-3e_3), \notag\\
\xi_2:=&\frac{1}{\sqrt{12\cdot 20}}(-e_3), \notag\\
\xi_3:=&\frac{1}{\sqrt{12\cdot 20}}(-6e_1+4e_2-3e_3), \\
\xi_4:=&\frac{1}{\sqrt{12\cdot 10}}(-e_1-e_2-e_3), \notag\\
\xi_5:=&\frac{1}{\sqrt{12\cdot  4}}(-2e_1-e_3). \notag
\end{align}
Moreover, the period map
$$\PER:\P^1\rightarrow
 \OP{O}(T)^\circ\backslash\MC{D}_T^\circ
 =\GAMd\backslash\DOMb$$
 defined in
 Subsection~\ref{SUBSECTtransA} maps $\T_\nu$ to
 $[\xi_\nu]:=\GAMd\cdot\xi_\nu$
 for $\nu=1,\ldots,5$.
The elliptic element corresponding to $\xi_\nu$ is
\begin{equation}
 g_\nu:=\xi_\nu\bmod\Q^\times\in\GAMd.
\end{equation}
\end{lem}
\begin{proof}
For $\nu=1,\ldots,5$,
 by Lemma~\ref{LEMindex}
 and Theorem~\ref{THMtranslist},
 there exists an element
 $g'_\nu\in
 \GAMd=\OP{O}(T[12])^\circ=\OP{O}(T)^\circ$
 such that
\begin{equation} \label{fixed_elliptic_element}
 \{ v\in T \bigm| g'_\nu(v)=v \}
 \cong (T(\bar{\X}^{(\nu)}_0))^\bot_T,\quad
 \{ v\in T \bigm| g'_\nu(v)=-v \}
 \cong T(\bar{\X}^{(\nu)}_0).
\end{equation}
In fact, by a direct computation,
 we can check that
 $g_\nu$ as in the assertion is
 an element in $\GAMd$ and
 satisfies the condition
 (\ref{fixed_elliptic_element})
 (see the table in Theorem~\ref{THMtranslist}).
Since $T(\bar{\X}^{(\nu)}_0)$ for $\nu=1,\ldots,5$
 are not isomorphic to each other,
 $g_\nu$ for $\nu=1,\ldots,5$ are not conjugate
 to each other in $\GAMd$.
On the other hand,
 by Remark~\ref{REMnofixedpoint} and
 Lemma~\ref{LEMtransgen},
 there exist at most five conjugacy classes
 of elliptic elements in $\GAMd$.
By (\ref{def_barpi}),
 we have $\PER(\T_\nu)=[\xi_\nu]$.
\end{proof}

\begin{rem}
The local monodromy around $\T_\nu$ on $T$
 is equal to $g_\nu$
 defined in Lemma~\ref{LEMelliptic_ele}
 up to conjugation and sign.
\end{rem}

\begin{lem} \label{LEMtable1}
We have the following:
\begin{itemize2}
\item[(i)]
$\GAMd/\GAMa\cong D_8\times C_2\times C_2$;
\item[(ii)]
$\GAMd/\GAMb \cong C_2\times C_2\times C_2$;
\item[(iii)]
$\GAMb/\GAMa \cong C_2\times C_2$.
\end{itemize2}
\end{lem}
\begin{proof}
By strong approximation theorem
 for quaternion algebras (cf.~\cite{borel81}),
we have a natural isomorphism
\begin{equation} \label{G235}
 \GAMd/\GAMa\cong G_2\times G_3\times G_5,
\end{equation}
where
\begin{align*}
G_2&=\frac{N(\ORDb\otimes\Z_2)}
 {\Q_2^\times \cdot
 \{ \gamma\in(\ORDb\otimes\Z_2)^1 \bigm|
 \gamma \equiv 1 \bmod 2\,\ORDb\otimes\Z_2 \} }, \\
G_p&=\frac{N(\ORDb\otimes\Z_p)}
 {\Q^\times_p \cdot
 (\ORDb\otimes\Z_p)^1},\quad p=3,5.
\end{align*}
By Theorem~\ref{THMlocalquat},
we have $G_3\cong G_5\cong C_2$.
We shall determine $G_2$.
We fix an identification
 $\ORDb\otimes\Z_2
 ={\displaystyle
\begin{pmatrix}
\Z_2 & \Z_2 \\
2\Z_2 & \Z_2
\end{pmatrix}
}$
 (see Theorem~\ref{THMlocalquat}
 and Lemma~\ref{LEMordext}).
Then,
 the image of $(\ORDb\otimes\Z_2)^1$ in $G_2$ is
 considered as the group in
${\displaystyle
\begin{pmatrix}
\Z/2\Z & \Z/2\Z \\
2\Z/4\Z & \Z/2\Z
\end{pmatrix}
}$ generated by
$\theta_1:={\displaystyle
\begin{pmatrix}
\bar{1} & \bar{1} \\
\bar{0} & \bar{1}
\end{pmatrix}
}$ and
$\theta_2:={\displaystyle
\begin{pmatrix}
\bar{1} & \bar{0} \\
\bar{2} & \bar{1}
\end{pmatrix}
}$,
which is isomorphic to $(\Z/2\Z)^2$.
Again by Theorem~\ref{THMlocalquat},
 $N(\ORDb\otimes\Z_2)
 /\Q^\times_2 \cdot (\ORDb\otimes\Z_2)^1$
 is of order $2$ and
generated by the image of
$z:={\displaystyle
\begin{pmatrix}
0 & 1 \\
2 & 0
\end{pmatrix}
}$.
Since the conjugation 
by $z$ interchanges $\theta_1$ and $\theta_2$,
 we have $G_2\cong D_8$.
This implies the assertion (i).
Since we have
$$\GAMd/\GAMb\cong (G_2/\langle \theta_1,\theta_2 \rangle)
\times G_3\times G_5\cong C_2\times C_2\times C_2,$$
the assertion (ii) holds.
Finally, we have
$$\GAMb/\GAMa = \langle \theta_1,\theta_2 \rangle \cong
C_2\times C_2,$$
and the assertion (iii) follows.
\end{proof}

Now, we shall show that $g_1,\ldots,g_5$
 generate $\GAMd$
 and determine their fundamental relations.
We use the following fact.
Let $\OP{vol}(~)$ denote
 the usual hyperbolic measure on $\MF{H}$.
If $S$ is a geodesic triangle
 with the interior angles $\alpha,\beta$ and $\gamma$,
 we have $\OP{vol}(S)=\pi-\alpha-\beta-\gamma$.

\begin{lem}[cf.\ \cite{vigneras80}]
 \label{LEMgeneral_fuchs}
Let $G$ be a Fuchsian group
 such that $\OP{vol}(G\backslash\MF{H})<\infty$.
Let $e_k,e_\infty$ and $g$
 be the number of elliptic points of order $k$,
 that of cusps and the genus,
 of the compactification of
 $G\backslash\MF{H}$, respectively.
Then, we have
$$\frac{1}{2\pi}\operatorname{vol}
(G\backslash\MF{H})
=2g-2+\sum_k \left( 1-\frac{1}{k} \right) e_k + e_\infty.$$
\end{lem}

\begin{lem} \label{LEMgenerators}
We have the following:
\begin{itemize2}
\item[(i)]
$\GAMd$ is generated by $g_1,\ldots,g_5$
 with the following fundamental relations:
\begin{equation} \label{gammarel}
g_1^2=\cdots=g_5^2=
g_1 \cdot g_3 \cdot g_5
\cdot g_4 \cdot g_2 = 1\text{;}
\end{equation}
\item[(ii)]
$\GAMb$ is
 the minimal normal subgroup of
 $\GAMd$ which contains
 $(g_\nu g_\mu)^2$ for $1 \leq \nu < \mu \leq 4$
 and $g_2 g_3$;
\item[(iii)]
$\GAMc$ is generated
 by $g_\nu$ and $g_1 g_\nu g_1$
 for $2\leq\nu\leq 4$;
\item[(iv)]
$g_\nu\not\in\GAMb$ for $\nu=1,\ldots,5$;
\item[(v)]
$g_\nu\in\GAMc$ if and only if $\nu=2,3,4$.
\end{itemize2}
\end{lem}
\begin{proof}[Proof (cf.~\rm{\cite{vigneras80}}\textit{).}]
We construct the fundamental domain
 of the action of $\GAMd$ explicitly.
Recall that
 the geodesic line segment $\overline{Q_1 Q_2}$
 in $\DOMb$,
 with the end points $Q_1$ and $Q_2$, is given by
$$\overline{Q_1 Q_2}=\left\{
 \frac{sQ_1+(1-s)Q_2}{\| sQ_1+(1-s)Q_2 \|} \in \DOMb
 \Biggm|
 0 \leq s \leq 1 \right\},$$
where $\| x \|=\sqrt{-x^2}$.
Let $P$ be the geodesic pentagon whose edges are
 $\overline{\xi_1 \xi_3}$,
 $\overline{\xi_3 \xi_5}$,
 $\overline{\xi_5 \xi_4}$,
 $\overline{\xi_4 \xi_2}$
 and $\overline{\xi_2 \xi_1}$.
By a direct computation, we find that
 $P$ is not self-intersecting and
 each interior angle of $P$ is $\pi/2$.
Hence, the $g \cdot P$ and $g \cdot P'$,
$g \in\GAMd$, cover $\DOMb$,
where $P'$ is the image of the reflection
 with respect to an edge of $P$,
 say, $\overline{\xi_1 \xi_3}$.
We have $\OP{vol}(P\cup P')=\pi$.
Now we apply Lemma~\ref{LEMgeneral_fuchs}
 to $\GAMd$.
By (\ref{def_pi}) and Lemma~\ref{LEMelliptic_ele},
 we have $e_2=5$, $e_k=0$ for $k>2$, $e_\infty=0$
 and $g=0$.
Hence, we have
 $\OP{vol}(\GAMd\backslash\DOMb)=\pi$
 by Lemma~\ref{LEMgeneral_fuchs}.
Therefore,
 $P\cup P'$ is the fundamental domain of $\GAMd$,
 and the assertion (i) follows.

By the proof of Lemma~\ref{LEMtable1},
 the norm map
 $\OP{Nr}:N^+({\ORDb})\rightarrow\Q^\times$
 induces an injection
 $\overline{\OP{Nr}}:\GAMd/\GAMb \rightarrow
 \Q^\times/\Q^{\times2}$.
We consider $\GAMd/\GAMb$
 as the subgroup of $\Q^\times/\Q^{\times2}$
 generated by
 $\bar{2},\bar{3},\bar{5}$
 through the map $\overline{\OP{Nr}}$.
The images of $g_1,g_2,g_3,g_4,g_5$ are
\setbox0=\hbox{$0123456789$}
$\bar{1}\bar{\rule{0ex}{\ht0}}\bar{0},
\bar{1}\bar{\rule{0ex}{\ht0}}\bar{5},
\bar{1}\bar{\rule{0ex}{\ht0}}\bar{5},
\bar{3}\bar{\rule{0ex}{\ht0}}\bar{0},
\bar{3}$, respectively.
Hence, $\GAMb$ contains
 $(g_\nu g_\mu)^2$ for $1\leq\nu<\mu\leq4$
 and $g_2 g_3$.
Let $G$ be the minimal normal subgroup
 of $\GAMd$
 which contains these elements.
We have $G\subset\GAMb$.
On the other hand,
 we can easily check that
 $[\GAMd:G]=8$.
By Lemma~\ref{LEMtable1}(ii),
 we have $G=\GAMb$.

The assertions (iii)--(v) follow immediately
 from the above argument.
\end{proof}

\begin{lem} \label{LEMtable2}
We have the following:
\begin{itemize2}
\item[(i)]
$\GAMc/\GAMa\cong D_8\times C_2$;
\item[(ii)]
$\GAMc/\GAMb \cong C_2\times C_2$.
\end{itemize2}
\end{lem}
\begin{proof}
We use the same notation
 as in the proof of Lemma~\ref{LEMtable1}.
Let $h_3$ and $h_5$ be
 the generators of $G_3$ and $G_5$, respectively.
By the proof of Lemma~\ref{LEMgenerators},
 the map $\overline{\OP{Nr}}$ gives
 an isomorphism between
 $\GAMc/\GAMb$ and
 \setbox0=\hbox{$0123456789$}
 $\langle \bar{2},
 \bar{1}\bar{\rule{0ex}{\ht0}}\bar{5} \rangle
 \cong C_2\times C_2$.
Hence, by the proof of Lemma~\ref{LEMtable1},
 the isomorphism (\ref{G235}) gives
 an isomorphism between
 $\GAMc/\GAMa$ and
 $G_2\times\langle h_3 h_5 \rangle
 \cong D_8\times C_2$.
\end{proof}

\begin{lem} \label{LEMnormalext}
Any subgroup of $\Gamma^1(\ORDb\otimes\R)$
 which contains $\GAMb$
 as a normal subgroup is contained in $\GAMd$.
\end{lem}
\begin{proof}
Let $G$ be a group as in the assertion.
By the Skolem--Noether theorem,
 $G$ is contained in $\Gamma^+(\ORDb\otimes\Q)$.
Let $g$ be an element in $\ORDb\otimes\Q$
 which represents an element in $G$.
For any prime number $p$,
 strong approximation theorem (cf.\ \cite{borel81}) implies that
 $g$ normalizes $(\ORDb\otimes\Z_p)^1$.
Hence, by Theorem~\ref{THMlocalquat}(v),
 we have $g\in N(\ORDb\otimes\Z_p)$.
As a consequence, we have $g\in N^+(\ORDb)$.
\end{proof}

\section{Quaternion modular embedding} \label{SECTmodular}
In this section, we construct
 a modular embedding
 $(\PHID,\PHIG)$
 of the domain $\DOMb$ defined by (\ref{def_dom})
 into the Siegel upper half-space $\MF{H}_2$
 of degree $2$:
\begin{equation}
 \begin{array}{ccc}
 \DOMb
  & {\displaystyle\mathop{\longrightarrow}^{\PHID}}
  & \MF{H}_2 \\
 \circlearrowleft & & \circlearrowleft \\
 \Gamma^1(\ORDb\otimes\R)
  & {\displaystyle\mathop{\longrightarrow}^{\PHIG}}
  & \OP{PSp}(4,\R).
 \end{array}
\end{equation}
Here $\PHID$ (resp.\ $\PHIG$) is an injection
 of domains (resp.\ groups)
 and the diagram is compatible.
In fact, the modular embedding
 $(\PHID,\PHIG)$ induces
 the following map,
 which is generically one-to-one
 (see Subsection~\ref{SUBSECTmodularB}):
\begin{equation*}
 \GAMc\backslash\DOMb\longrightarrow
 \OP{PSp}(4,\Z)\backslash\MF{H}_2.
\end{equation*}

\subsection{Construction of $(\PHID,\PHIG)$}
We construct a modular embedding of $\DOMb$
into $\MF{H}_2$
as a special case of Kuga-Satake construction.
Recall that we fix $I$ as in (\ref{defI}).
We define a domain $\MC{D}'$,
 which is isomorphic to the Siegel
 upper half-space $\mathfrak{H}_2$ of degree $2$
 (see Lemma~\ref{LEMphi2dash}), by
\begin{equation}
 \MC{D}'
 =\{\rho\in \OP{Sp}(\ORDb\otimes\R,\Psi_I)\bigm|
 \rho^2=-1,~\Psi_I(v,\rho(v))>0~(\forall v\neq 0) \}.
\end{equation}
Each element in $\MC{D}'$ defines a complex structure
 on $\ORDb\otimes\R$.
The symplectic group $\OP{Sp}(\ORDb\otimes\R,\Psi_I)$
acts on $\MC{D}'$ as follows:
\begin{equation}
\gamma\cdot\rho=\gamma\rho\gamma^{-1},\quad
\gamma\in\OP{Sp}(\ORDb\otimes\R,\Psi_I),~\rho\in\MC{D}'.
\end{equation}
We identify $(\MC{D}',\OP{PSp}(\ORDb\otimes\R,\Psi_I))$
 with $(\MF{H}_2,\OP{PSp}(4,\R))$
 by the following lemma.

\begin{lem} \label{LEMphi2dash}
The following correspondence (compatible with the action)
 is an isomorphism between
 $(\MC{D}',\OP{PSp}(\ORDb\otimes\R,\Psi_I))$
 and $(\MF{H}_2,\OP{PSp}(4,\R))$:
\begin{equation} \label{dom_ident}
\begin{array}{c}
 \MC{D}'\longleftrightarrow\mathfrak{H}_2;\quad
 \rho\longleftrightarrow\tau, \vspace{1ex} \\
 \begin{pmatrix} \rho(\alpha_1) \\ \rho(\alpha_2) \end{pmatrix}
 =\OP{Re}\tau \cdot \begin{pmatrix} \alpha_1 \\ \alpha_2 \end{pmatrix}
 +\OP{Im}\tau \cdot \begin{pmatrix} \beta_1 \\ \beta_2 \end{pmatrix},
\end{array}
\end{equation}
\begin{equation}
\begin{array}{c}
 \OP{PSp}(\ORDb\otimes\R,\Psi_I)
 \longleftrightarrow \OP{PSp}(4,\R);\quad
 \pm\rho\longleftrightarrow \pm M, \vspace{1ex} \\
 (\rho(\alpha_1),\rho(\alpha_2),\rho(\beta_1),\rho(\beta_2))
 =(\alpha_1,\alpha_2,\beta_1,\beta_2)M.
\end{array}
\end{equation}
Here $(\alpha_1,\alpha_2,\beta_1,\beta_2)$ is a symplectic basis
 of $(\ORDb,\Psi_I)$ defined by (\ref{alphabeta}).
\end{lem}
\begin{proof}
See e.g.\ \cite{satake}.
\end{proof}

The following is a special case of Kuga--Satake construction.

\begin{prop} \label{PROPembdomain}
Let $\Phi$ be the map defined by
\begin{equation} \label{defphi}
 \Phi:{\ORDb\otimes\R}\longrightarrow
 \OP{End}_\R(\ORDb\otimes\R);\quad
 \gamma\longmapsto (v\longmapsto \gamma v)\text{.}
\end{equation}
Then, $\Phi$ induces the following maps
 $\Phi_{\rm dom}$ and $\Phi_{\rm grp}$
 between domains and groups, respectively:
\begin{gather}
 \Phi_{\rm dom}:
  \DOMb\longrightarrow\MC{D}'
  \mathop{\longrightarrow}^{(\ref{dom_ident})}
  \MF{H}_2;\quad
  x\longmapsto\Phi(x) \\
 \Phi_{\rm grp}:
  \Gamma^1(\ORDb\otimes\R)
  \longrightarrow\OP{PSp}(\ORDb\otimes\R,\Psi_I);\quad
  \pm\gamma\longmapsto\pm\Phi(\gamma).
\end{gather}
Moreover,
 $\Phi_{\rm dom}$ and $\Phi_{\rm grp}$ are
 compatible with the actions.
In other words,
 we have
 $\Phi_{\rm dom}(\pm\gamma\cdot x)
 =\Phi_{\rm grp}(\pm\gamma)\cdot\Phi_{\rm dom}(x)$
 for all $\pm\gamma\in\Gamma^1(\ORDb\otimes\R)$
 and $x\in\DOMb$.
\end{prop}
\begin{proof}
We have
$\Phi((\ORDb\otimes\R)^1)\subset
 \OP{Sp}(\ORDb\otimes\R,\Psi_I)$.
In fact,
 for $\gamma\in(\ORDb\otimes\R)^1$
 and $v,w\in\ORDb\otimes\R$, we have
$$\Psi_I(\Phi(\gamma)(v),\Phi(\gamma)(w))
=\OP{Tr}(\gamma v I \bar{w} \bar{\gamma})
=\OP{Tr}(\gamma v I \bar{w} \gamma^{-1})
=\OP{Tr}(vI\bar{w})=\Psi_I(v,w).$$
For $x=I/\sqrt{-I^2}\in\DOMb$,
 we shall show that $\Phi(x)\in\MC{D}'$.
Using the degenerate symmetric bilinear form
 $b(v,w):=\Psi_I(v,xw)$,
 we have the orthogonal decomposition
$${\ORDb\otimes\R}
 =\R 1\oplus \R x\oplus (x)^\bot_{\MC{I}_\R},$$
 where $\MC{I}_\R=\MC{I}(\ORDb\otimes\R)$.
Since we have
$(x)^\bot_{\MC{I}_\R}
 =\{ y\in\MC{I}_\R \bigm|
 \OP{Tr}(x\bar{y})=0 \}$,
 we find that $b(~,~)$ is positive definite
 on each component in the above decomposition.
Hence we have $\Phi(x)\in\MC{D}'$.
Since $\DOMb$ is connected and $\Phi(x)\in\MC{D}$,
 we have $\Phi(\DOMb)\subset\MC{D}'$.
The compatibility of
 $\Phi_{\rm dom}$ and $\Phi_{\rm grp}$
 follows immediately from the definition.
\end{proof}

Let $\OP{PSp}(4,\Z)[2]= \{ \pm M\in\OP{PSp}(4,\Z) \bigm|
 M\equiv 1_4 \bmod 2 \} $
 be the congruence subgroup of $\OP{PSp}(4,\Z)$ of level $2$.
By the above construction, we have
\begin{equation} \label{pull-back_of_groups}
 \begin{array}{c}
 \PHIG^{-1}(\OP{PSp}(4,\Z))=\GAMb, \vspace{1ex} \\
 \PHIG^{-1}(\OP{PSp}(4,\Z)[2])=\GAMa.
 \end{array}
\end{equation}

\subsection{Properties of $(\PHID,\PHIG)$}
\label{SUBSECTmodularB}
By the previous subsection,
 the modular embedding $(\PHID,\PHIG)$ induces a map
\begin{equation} \label{modular_map_1}
 \GAMb\backslash\DOMb\longrightarrow
 \OP{PSp}(4,\Z)\backslash\MF{H}_2.
\end{equation}
However, this map is not (generically) one-to-one.
In fact, we have the following.

\begin{prop} \label{PROPmodular_inj}
The modular embedding $(\PHID,\PHIG)$ induces
 well-defined and generically one-to-one maps
\begin{gather}
 \GAMc\backslash\DOMb\longrightarrow
 \OP{PSp}(4,\Z)\backslash\MF{H}_2, \\
 \GAMa\backslash\DOMb\longrightarrow
 \OP{PSp}(4,\Z)[2]\backslash\MF{H}_2,
\end{gather}
 where $\GAMc$ and $\GAMa$ are defined
 in Definition~\ref{DEFgams}.
\end{prop}

Let $G_\R$ be the stabilizer of $\PHID(\DOMb)$ in $\OP{PSp}(4,\R)$
 and $G_\Z$ that in $\OP{PSp}(4,\Z)$.
We set
\begin{equation}
Z=\{
z\in\OP{PSp}(4,\R) \bigm| z\cdot \xi=\xi~
(\forall \xi\in\PHID(\DOMb)) \}.
\end{equation}
The map (\ref{modular_map_1}) is decomposed as follows:
\begin{equation} \label{modular_decomposition}
 \GAMb\backslash\DOMb\longrightarrow
 G_\Z\backslash\PHID(\DOMb)\longrightarrow
 \OP{PSp}(4,\Z)\backslash\MF{H}_2.
\end{equation}
Since the second map in (\ref{modular_decomposition}) is
 generically one-to-one,
 it is sufficient to show that $G_\Z/(G_\Z\cap Z)\cong\Gamma'$
 in order to prove the former part of
 Proposition~\ref{PROPmodular_inj}.
First, we study the structure of $G_\R$.

\begin{lem} \label{LEM_G_R}
We have the following:
\begin{itemize2}
\item[(i)]
 $G_\R=\PHIG(\Gamma^1(\ORDb\otimes\R))\times Z$;
\item[(ii)]
 $Z= \{
 \pm(x\mapsto x\eta)\bigm| \eta\in\ORDb\otimes\R,~
 \eta I\bar{\eta}=I \}$;
\item[(iii)]
 $Z\cap\OP{PSp}(4,\Z)= \{ \pm 1_4 \}$.
\end{itemize2}
\end{lem}
\begin{proof}
We have $\ORDb\otimes\R=\R+\R\cdot\DOMb$.
Hence, for any element $\pm\varphi\in Z$, we have
$$\varphi\circ\Phi(x)=\Phi(x)\circ\varphi
 \in\OP{End}_\R(\ORDb\otimes\R) \quad
 (\forall x\in \ORDb\otimes\R).$$
In particular, $Z$ centralizes $\PHIG(\Gamma^1(\ORDb\otimes\R))$.
Since the automorphism group of $\DOMb\cong\MF{H}$ is
 $\Gamma^1(\ORDb\otimes\R)\cong\OP{PSL}(2,\R)$,
 we see that $\Gamma^1(\ORDb\otimes\R)$ and $Z$ generate $G_\R$.
Thus the assertion (i) follows.
We set $\eta=\varphi(1)$.
Since we have $\Phi(x)(1)=x$, we have
 $\varphi(x)=x\eta$ for all $x\in \ORDb\otimes\R$.
The condition $(x\mapsto x\eta)\in\OP{Sp}(\ORDb\otimes\R,\Psi_I)$
 is equivalent to the condition
 $\Psi_I(x\eta,y\eta)=\Psi_I(x,y)$,
 i.e., $\OP{Tr}(x\eta I\bar{\eta}\bar{y})=\OP{Tr}(xI\bar{y})$
 for all $x,y\in\ORDb\otimes\R$.
Since the bilinear form $(x,y)\mapsto\OP{Tr}(xy)$ is non-degenerate,
 this condition is equivalent to
 the condition $\eta I\bar{\eta}=I$.
This implies the assertion (ii).
We set
\begin{align} \label{def_ortho_basis}
 y_1&={\displaystyle\frac{e_1-e_2}{6}=w_1}, \notag \\
 y_2&={\displaystyle\frac{2e_1+2e_2+e_3}{4}=-1+3w_1+2w_2+2w_3}, \\
 y_3&={\displaystyle-\frac{e_1+e_2+e_3}{2}
  =2-3w_1-2w_2-4w_3=-15I}. \notag
\end{align}
We have
 $\OP{Tr}(y_i,I)=y_i I+I y_i=0$ ($i=1,2$).
Assume that $(x\mapsto x\eta)\in\OP{Sp}(\ORDb,\Psi_I)$.
Then, we have $\eta\in\ORDb$, $\eta I \bar{\eta}=I$
 and $\OP{Nr}(\eta)=\pm 1$.
If $\OP{Nr}(\eta)=1$, we have $\eta=a+by_3$ with $a,b\in\Z$.
Since we have $\OP{Nr}(\eta)=a^2+30b^2$, we have $\eta=\pm 1$.
If $\OP{Nr}(\eta)=-1$, we have $\eta=ay_1+by_2$ with $a,b\in\Z$.
Then, we have $\OP{Nr}(\eta)=-2 a^2-15 b^2=-1$,
 which is a contradiction.
Therefore, the assertion (iii) follows.
\end{proof}

\begin{rem} \label{REMpo2}
By the proof of Lemma~\ref{LEM_G_R}, $Z$ is anti-isomorphic to
\begin{equation*}
 \left( \{ \eta=a+by_3 \bigm| \OP{Nr}(\eta)=1 \}
 \sqcup
 \{ \eta=ay_1+by_2 \bigm| \OP{Nr}(\eta)=-1 \} \right)/\pm1
 \cong\OP{PO}(2,\R).
\end{equation*}
\end{rem}

\begin{proof}
[Proof of Proposition~\ref{PROPmodular_inj}]
Let $f$ denote the natural map
\begin{equation}
 f:G_\R\longrightarrow\ G_\R/Z\cong\Gamma^1(\ORDb\otimes\R).
\end{equation}
By Lemma~\ref{LEM_G_R}(i), $f(G_\Z)$ is a normal extension of $\GAMb$.
Hence, we have $f(G_\Z)\subset\GAMd$ by Lemma~\ref{LEMnormalext}.
By Lemma~\ref{LEM_G_R}(ii),
 for $\xi\in\DOMb$
 and $\eta\in\ORDb\otimes\R$ with $\eta I\bar{\eta}=I$, we have
 $(x\mapsto\xi x\eta)\in\OP{Sp}(\ORDb\otimes\R,\Psi_I)$.
By Lemma~\ref{LEMelliptic_ele}, for $\nu=1,\ldots,5$, we have
\begin{align*}
 & \,(x\mapsto\xi_\nu x \eta)\in\OP{Sp}(\ORDb,\Psi_I) \\
\Longleftrightarrow & \,\xi_\nu\ORDb\eta=\ORDb\xi_\nu\eta=\ORDb \\
\Longleftrightarrow & \,\xi_\nu \eta\in\ORDb^\times 
\end{align*}
We set
\begin{equation}
 \eta_2=\frac{2e_1+2e_2+e_3}{\sqrt{240}},\quad
 \eta_4=\xi_4=\frac{-e_1-e_2-e_3}{\sqrt{120}}.
\end{equation}
We have $\eta_\nu I \bar{\eta}_\nu=I$
 and $\xi_\nu \eta_\nu\in\ORDb^\times$ for $\nu=2,4$.
In fact, by a direct computation, we have
 $\xi_2\eta_2=1+(e_1-e_2)/6=1+w_1$, $\OP{Nr}(\xi_2\eta_2)=-1$
 and $\xi_4\eta_4=-1$.
Hence, we have $g_\nu=\xi_\nu\bmod\Q^\times\in f(G_\Z)$
 for $\nu=2,4$.
As a consequence, we have
 $\GAMc\subset f(G_\Z)\subset\GAMd$
 (see (\ref{def_gams}) and (\ref{table_gamma})).
Assume that $f(G_\Z)=\GAMd$.
Then, by Lemma~\ref{LEM_G_R} and Remark~\ref{REMpo2},
 we have a injection
$$C_2^3\cong f(G_\Z)/\GAMb\cong G_\Z/\PHIG(\GAMb)
 \longrightarrow Z\cong\OP{PO}(2,\R).$$
Since $\OP{PO}(2,\R)$ does not contain $C_2^3$,
 this is a contradiction.
By $[\GAMd:\GAMc]=2$, we have $f(G_\Z)=\GAMc$.
This implies the former part of the proposition.

We have a natural sequence
\begin{equation*}
 \GAMd/\GAMa \mathop{\longrightarrow}^{\alpha}
 G_\Z/(G_\Z\cap\OP{PSp}(4,\Z)[2]) \mathop{\longrightarrow}^{\beta}
 \OP{Sp}(4,\mathbb{F}_2),
\end{equation*}
 where $\alpha$ is surjective and $\beta$ is injective.
We shall show that $\left|\OP{Im}(\beta)\right|=2^4$.
Let $M_\nu$ denote the element
 $(x\mapsto x\eta_\nu)
 \in\OP{Sp}(\ORDb\otimes\R,\Psi_I)=\OP{Sp}(4,\R)$ for $\nu=2,4$.
By a direct computation, in $\OP{Sp}(4,\mathbb{F}_2)$, we have
\begin{equation} \label{generators_level_2}
\begin{array}{r}
\PHIG(g_2)M_2\equiv\PHIG(g_1 g_2 g_1)M_2\equiv
\begin{pmatrix}
0 &  1 &  0 & 0 \\
1 &  0 &  0 & 0 \\
0 &  0 &  0 & 1 \\
0 &  0 &  1 & 0
\end{pmatrix}, \vspace{1ex} \\
\PHIG(g_3) M_2\equiv
\begin{pmatrix}
 1 &  0 &  1 & 1 \\
 0 &  1 &  1 & 1 \\
 1 &  1 &  1 & 0 \\
 1 &  1 &  0 & 1
\end{pmatrix}, \vspace{1ex} \\
\PHIG(g_4)M_4\equiv\PHIG(g_1 g_4 g_1)M_4\equiv
\begin{pmatrix}
 1 &  0 &  0 & 1 \\
 0 &  1 &  1 & 0 \\
 0 &  0 &  1 & 0 \\
 0 &  0 &  0 & 1
\end{pmatrix}, \vspace{1ex} \\
\PHIG(g_1 g_3 g_1)M_2\equiv
\begin{pmatrix}
 0 & 1 & 0 &  0 \\
 1 & 0 & 0 &  0 \\
 1 & 1 & 0 &  1 \\
 1 & 1 & 1 &  0
\end{pmatrix}.
\end{array}
\end{equation}
By Lemma~\ref{LEMgenerators}(iii) and \ref{LEM_G_R},
 $\OP{Im}(\beta)$ is generated by these elements.
We can check that $\left|\OP{Im}(\beta)\right|=2^4$.
This implies $\alpha$ is bijective.
Therefore, we have $f(G_\Z\cap\OP{PSp}(4,\Z)[2])=\GAMa$.
This implies the latter part of the proposition.
\end{proof}

\section{The inverse period map and automorphic forms} \label{SECTperiod}
In this section,
we construct the inverse period map of the family $\{\X_\T\}$
by using automorphic forms on the period domain.
To construct automorphic forms on $\DOMb$,
 we pull-back theta constants on $\mathfrak{H}_2$ by
 the modular embedding $(\PHIG,\PHID)$
 defined in Proposition~\ref{PROPembdomain}.
The maps in the following diagram are
 defined in the previous sections:
$$\begin{array}{c@{}c@{}c@{}c@{}c@{}c@{}c}
 \P^1
   & \,\smash{\displaystyle\mathop{\longrightarrow}^{\sim}_{\PER}}\,
   & \OP{O}(T)^\circ\backslash\MC{D}_T^\circ
   & \,\smash{\displaystyle\mathop{\cong}
     _{(\ref{identification_gam})}}
   & \GAMd\backslash\DOMb & & \\
 & & & & {\scriptstyle~\text{gen.}\,2:1}\!\uparrow \\
 & & & & \GAMc\backslash\DOMb
   & {\displaystyle\mathop{\longrightarrow}^{\text{gen.}\,1:1}}
   & \OP{PSp}(4,\Z)\backslash\MF{H}_2 \\
 & & & & {\scriptstyle\text{gen.}\,16:1}\!\uparrow &
   & \uparrow \\
 & & &
   & \GAMa\backslash\DOMb
   & \,{\displaystyle\mathop{\longrightarrow}^{\text{gen.}\,1:1}}\,
   & \OP{PSp}(4,\Z)[2]\backslash\MF{H}_2.
\end{array}$$

\subsection{Theta constants}

We recall some results on theta constants on $\mathfrak{H}_2$
and the Siegel modular variety of genus $2$ and level $2$
 (cf.\ \cite{vandergeer82}).
For a row vector
 $m=(m',m'')\in
 {\displaystyle\frac{1}{2}}\Z^2\times
 {\displaystyle\frac{1}{2}}\Z^2$,
the theta constants $\theta_m(\tau)$ of character $m$
 is defined by
\begin{equation} \label{theta}
\theta_m(\tau)=\sum_{n\in\Z^2}
\exp\left(2\pi i\left(
\frac{1}{2}(n+m')\cdot\tau\cdot
{}^t(n+m')+
(n+m')\cdot
{}^t m''\right)\right),
\end{equation}
 where $\tau\in\MF{H}_2$.
The character $m$ is called even,
 if $4m'\cdot {}^t m''\equiv 1 \bmod 2$.
We have $\theta_m\equiv 0$
 if and only if $m$ is even,
 and we obtain ten non-zero theta constants.
We use the following theorems.

\begin{thm}[\cite{igusa64}] \label{THMconstants}
For
 $M=\begin{pmatrix}A&B\\C&D\end{pmatrix}\in\OP{Sp}(4,\Z)$
 and $m\in{\displaystyle\frac{1}{2}}\Z^4$,
 we set
\begin{equation}
 M\cdot m
 =m M^{-1}
  +\frac{1}{2}
   \left( (C\cdot {}^t\! D)_0, (A\cdot {}^t\! B)_0 \right)
\end{equation}
 and
$$X_0=\begin{pmatrix}x_{11}&x_{22}\end{pmatrix}
 \quad\text{for}\quad
 X
 =\begin{pmatrix}
   x_{11} & x_{12} \\
   x_{21} & x_{22}
  \end{pmatrix}.$$
Then, for $\tau\in\MF{H}_2$,
\begin{equation} \label{theta_trans}
 \theta_{M\cdot m}(M\cdot\tau)^4
 =\pm\det(C\tau+D)^2 \theta_m(\tau)^4,
\end{equation}
 where the sign depends only on $M$ and $m$,
 and positive if $M\in\OP{Sp}(4,\Z)[2]$.
In particular,
 the function $\theta_m^4$ is
 a Siegel modular form of weight $2$ and
 level $2$ if $m$ is even.
(As for the sign of (\ref{theta_trans}),
 see \cite{igusa64}.)
\end{thm}


\begin{thm}[\cite{igusa64,vandergeer82}] \label{THMmodularvar}
The vector space of modular forms generated by
 $\theta_m^4$'s is $5$-dimensional.
The following map becomes a biholomorphic embedding of the Siegel
 modular variety of genus $2$ and level $2$:
\begin{equation}
 \Theta:\OP{Sp}(4,\Z)[2]\backslash\mathfrak{H}_2\longrightarrow
 \P^9;\quad
 [\tau]\longmapsto(\ldots:\theta_m(\tau)^4:\ldots).
\end{equation}
Here $m$ runs through all even characters modulo $\Z$.
\end{thm}

\subsection{Defining equation}

In this subsection, we determine the defining equations
 of the image of
\begin{equation} \label{map_modular_inj}
 F:\MC{C}_4:=\GAMa\backslash\DOMb\longrightarrow
 \OP{PSp}(4,\Z)[2]\backslash\mathfrak{H}_2
 \mathop{\longrightarrow}^{\Theta}\P^4.
\end{equation}
The map $F$ is generically one-to-one
 by Proposition~\ref{PROPmodular_inj}.
We set
\begin{align} \label{def_t}
\Theta_1&=\XXXb+\XXXc, \notag \\
\Theta_2&=\XXXb-\XXXc, \notag \\
\Theta_3&=2\,\XXXa-\XXXb-\XXXc+2\,\XXXd, \\
\Theta_4&=\XXXb-\XXXc+2\,\XXXe, \notag \\
\Theta_5&=2\,\XXXa-\XXXb+\XXXc-2\,\XXXd-2\,\XXXe, \notag
\end{align}
and define an algebra homomorphism
\begin{equation}
Ev:\C[Y_1,\ldots,Y_5]\longrightarrow\C[F^*\Theta_1,\ldots,F^*\Theta_5]:
Y_i\mapsto F^*\Theta_i.
\end{equation}
By linear relations between $\theta_m^4$'s
 (see \cite{vandergeer82}),
 the map $\Theta$ in Theorem~\ref{THMmodularvar}
 factors through
\begin{equation}
\OP{Proj}(\C[F^*\Theta_1,\ldots,F^*\Theta_5])
\longrightarrow \OP{Proj}(\C[Y_1,\ldots,Y_5])\cong\P^4.
\end{equation}
We define actions $\sigma_1,\ldots,\sigma_4$ on $\C[Y_1,\ldots,Y_5]$ by
\begin{equation} \label{linear_action}
\begin{array}{l}
 \sigma_1:(Y_1,Y_2,Y_3,Y_4,Y_5)
  \longmapsto(Y_1,-Y_2,Y_3,Y_4,Y_5), \vspace{1ex} \\
 \sigma_2:(Y_1,Y_2,Y_3,Y_4,Y_5)
  \longmapsto(-Y_1,Y_2,Y_3,-Y_4,Y_5), \vspace{1ex} \\
 \sigma_3:(Y_1,Y_2,Y_3,Y_4,Y_5)
  \longmapsto(Y_1,Y_2,Y_3,Y_5,Y_4), \vspace{1ex} \\
 \sigma_4:(Y_1,Y_2,Y_3,Y_4,Y_5)
  \longmapsto(-Y_1,Y_2,Y_3,Y_4,-Y_5),
\end{array}
\end{equation}
and that of $M=\begin{pmatrix}A&B\\C&D\end{pmatrix}\in\OP{Sp}(4,\Z)$ on
$\Theta_i$'s by
\begin{equation}
(\Theta_i|M)(\tau)=\det(C\tau+D)^{-2}\Theta_i(M\cdot\tau).
\end{equation}
Then using notations in (\ref{generators_level_2}), we have
\begin{align*}
Ev(\sigma_1(Y_k))&=Ev(Y_k)|(\PHIG(g_2)M_2) \\
Ev(\sigma_2(Y_k))&=Ev(Y_k)|(\PHIG(g_3) M_2) \\
Ev(\sigma_3(Y_k))&=Ev(Y_k)|(\PHIG(g_4)M_4) \\
Ev(\sigma_4(Y_k))&=Ev(Y_k)|(\PHIG(g_1 g_3 g_1)M_2).
\end{align*}
The action of each $\sigma_i$ induces an action
 on $\MC{C}_4$ via the map $F$,
 which is also denoted by $\sigma_i$.

Now we come to the key theorem for describing
 the inverse period map of the family $\{ \X_\T \}$.

\begin{thm} \label{THMdef_ideal}
Let $\MF{a}\subset\C[Y_1,\ldots,Y_5]$ be the defining ideal
 of the image of $F$ in (\ref{map_modular_inj}).
Then, the following polynomials are contained in $\MF{a}$:
\begin{equation} \label{def_ideal}
 \begin{array}{l}
 f_2=Y_1^2-Y_2^2-Y_3^2+Y_4^2+Y_5^2, \vspace{1ex} \\
 f_4=Y_1^2 Y_4^2 + Y_1^2 Y_5^2 + Y_4^2 Y_5^2
      - 2 Y_1 Y_3 Y_4 Y_5, \vspace{1ex} \\
 f_5=111375 Y_1^5 - 180630 Y_1^3 Y_2^2
      - 35721 Y_1 Y_2^4 - 107280 Y_1^3 Y_3^2 \\
      + 65808 Y_1 Y_2^2 Y_3^2 - 4096 Y_1 Y_3^4
      + 291600 Y_1^2 Y_3 Y_4 Y_5
      - 84240 Y_2^2 Y_3 Y_4 Y_5 - 61440 Y_3^3 Y_4 Y_5.
 \end{array}
\end{equation}
Moreover, the quotient ring $\C[Y_1,\ldots,Y_5]/(f_2,f_4,f_5)$
is of dimension $2$.
\end{thm}
%
%

For the proof of Theorem~\ref{THMdef_ideal},
 we use computations of some special CM points
 in the image of $F$,
 which are given in Appendix.
We prepare for two lemmas.

\begin{lem} \label{LEMlinear_relation}
The ideal $\MF{a}$ contains no element of degree $1$.
\end{lem}
\begin{proof}
Let $p_j\in\DOMb$ ($j=1,\ldots,8$) be as in Appendix.
From the result in Appendix, we can check that the matrix
 $(\Theta_i\circ\PHID(p_j))_{ij}$ is of maximal rank by a direct computation.
The assertion immediately follows from this.
\end{proof}

\begin{lem} \label{LEMdegree}
$\deg(F^*\MC{O}_{\P^4}(n))=32n$.
\end{lem}
\begin{proof}
Let $f\in\C[Y_1,\ldots,Y_5]$ be a homogeneous polynomial.
Then,
 the pull-back $f(\Theta_1,\ldots,\Theta_5)\circ\PHID$
 is an automorphic form
 of weight $4\deg(f)$ for $\GAMa$.
Hence, by Riemann--Roch formula and (\ref{table_gamma}), we have
\begin{align*}
 \deg(F^*\MC{O}_{\P^4}(n))
 &=\frac{\OP{vol}(\GAMa\backslash\DOMb)}{2\pi}\cdot\frac{4n}{2} \\
 &=\frac{32\OP{vol}(\GAMd\backslash\DOMb)}{2\pi}\cdot 2n \\
 &=32n.
\end{align*}
Here we use $\OP{vol}(\GAMd\backslash\DOMb)=\pi$,
 which we showed in the proof of Lemma~\ref{LEMgenerators}.
\end{proof}

Let $\MF{S}$ be the group generated by
 $\sigma_1,\ldots,\sigma_4$ in (\ref{linear_action}) and
 $\chi$ the character
 of the $1$-dimensional representation space
 $\C Y_1$ of $\MF{S}$.
The group $\OP{Ker}(\chi)$ is generated
 by $\sigma_1,\sigma_2\sigma_4,\sigma_3$.
Note that $\OP{Ker}(\chi)$ is a normal subgroup of $\GAMd/\GAMa$
 (cf.\ Lemma~\ref{LEMgenerators}).
We consider a sequence of curves and natural maps:
\begin{equation}\begin{CD}
 \DOMb   @>{\text{\'etale}}>{\pi_4}>
 \MC{C}_4@>{8:1}>{\pi_3}>
 \MC{C}_3@>{2:1}>{\pi_2}>
 \MC{C}_2@>{2:1}>{\pi_1}>
 \MC{C}_1,
\end{CD}\end{equation}
where the curves are defined as follows:
\begin{equation} \label{table_C_i}
 \begin{array}{c|c|ccccc|c}
  & \text{curve} &[\xi_1]&[\xi_2]&[\xi_3]&[\xi_4]&[\xi_5]
    & \text{genus} \\
  \hline
  \MC{C}_1 & \GAMd\backslash\DOMb              &1&1&1&1&1& 0 \\
  \MC{C}_2 & \MF{S}\backslash \MC{C}_4              &1&2&2&2&1& 0 \\
  \MC{C}_3 & \OP{Ker}(\chi)\backslash \MC{C}_4 &2&4&2&4&2& 0 \\
  \MC{C}_4 & \GAMa\backslash\DOMb         &16&16&16&16&16& 9
 \end{array}.
\end{equation}
Here $[p]$ denotes the image of $p\in\DOMb$ in $\MC{C}_1$,
 i.e., $[p]=\pi_1\circ\pi_2\circ\pi_3\circ\pi_4(p)$.
We already showed 
 this table  in (\ref{table_gamma}),
 except for the row of $\MC{C}_3$.
Let $\Gamma''\subset\GAMd$ be the subgroup corresponding
 to $\OP{Ker}(\chi)$, i.e., $\Gamma''/\GAMa=\OP{Ker}(\chi)$.
The cardinality of $(\pi_1 \circ \pi_2)^{-1}([\xi_\nu])$
 is equal to
 $4$ (resp.\ $2$) if $g_\nu\in\Gamma''$
 (resp.\ $g_\nu\not\in\Gamma''$).
We can easily check that $g_\nu\in\Gamma''$
 if and only if $\nu=2,4$ by (\ref{linear_action}).

\begin{proof}[Proof of Theorem~\ref{THMdef_ideal}]
We only prove that $f_5$ vanishes on the image of $F$.
We can determine $f_2,f_4$ in the same way.
Under the equality $f_2\equiv 0 \pmod{\MF{a}}$,
 the equality $f_4\equiv 0 \pmod{\MF{a}}$ is
 equivalent to the well-known equality (\cite{igusa64,vandergeer82})
\begin{equation} \label{igusa_rel}
 t_8-\frac{1}{4}t_4^2=0,
\end{equation}
where $t_k=\sum\theta_m^{2k}$.
Let $\MC{K}$ be the canonical sheaf of $\DOMb$ and
 $\MC{L}$ an invertible sheaf on $\MC{C}_4$ defined by
\begin{equation}
 \MC{L}=((\pi_4)_*\MC{K})^{\GAMa}.
\end{equation}
We fix an identification
 $F^*(\MC{O}_{\P^4}(1))=\MC{L}^{\otimes 2}$.
For any $\sigma\in\MF{S}$
 and $x\in\MC{C}_4$ fixed by $\sigma$,
 the action of $\sigma$ on the fiber $\MC{L}\otimes\kappa(x)$
 is of order $1$ or $2$,
 and the action of $\sigma$
 on $\MC{L}^{\otimes 2}\otimes\kappa(x)$
 is trivial.
Therefore,
 an invertible sheaf
\begin{equation}
 \MC{M}:=((\pi_2\circ\pi_3)_*\MC{L}^{\otimes 10})^{\OP{Ker}(\chi)}
\end{equation}
 on $\MC{C}_2$
 satisfies
 $(\pi_2\circ\pi_3)^*\MC{M}
 \cong\MC{L}^{\otimes 10}
 =F^*(\MC{O}_{\P^4}(5))$.
By Lemma~\ref{LEMdegree},
 we have
\begin{equation} \label{degree_G}
 \deg(\MC{M})
 =\deg(F^*(\MC{O}_{\P^4}(5)))/\lvert\MF{S}\rvert
 =160/16=10.
\end{equation}
For a module $V$ with an action of $\MF{S}$,
 let $V^\chi$ denote the $\chi$-eigenspace of $V$.
We have a natural isomorphism
\begin{equation}
 \varphi:
 H^0(\MC{C}_3,\pi_2^*\MC{M})^\chi
 \CONGL
 H^0(\MC{C}_2,((\pi_2)_*\MC{O}_{\MC{C}_3})^\chi\otimes\MC{M}).
\end{equation}
Since $\deg(((\pi_2)_*\MC{O}_{\MC{C}_3})^\chi)=-1$,
 we have
\begin{equation} \label{degree_nine}
\deg(((\pi_2)_*\MC{O}_{\MC{C}_3})^\chi\otimes\MC{M})=9.
\end{equation}
We set $R=\C[Y_1,\ldots,Y_5]/(f_2,f_4)=\bigoplus R_d$,
 where $R_d$ is the homogeneous part of $R$ of degree $d$.
By a direct computation, we can check that
 the following nine monomials modulo $(f_2,f_4)$
 form a $\C$-basis of $R_5^\chi$:
\begin{equation} \label{monomials}
\begin{array}{c}
 Y_1^5,~ Y_1^3 Y_2^2,~ Y_1 Y_2^4,~ Y_1^3 Y_3^2,~
 Y_1 Y_2^2 Y_3^2,~ \vspace{1ex} \\
 Y_1 Y_3^4,~ Y_1^2 Y_3 Y_4 Y_5,~
 Y_2^2 Y_3 Y_4 Y_5,~ Y_3^3 Y_4 Y_5.
\end{array}
\end{equation}
Let $M_i=M_i(Y_1,\ldots,Y_5)$ denotes the $i$-th monomial in (\ref{monomials}),
and set
\begin{gather*}
q_j=
(\Theta_1(\PHID(p_j))/\Theta_1(\PHID(p_j)),\ldots,\Theta_5(\PHID(p_j))/\Theta_1(\PHID(p_j))),\\
k_j=\Q(\{
 M_{i}({q_j}) \bigm| i=1,\ldots,9 \} ) .
\end{gather*}
Then, we have
$$k_j=
 \begin{cases}
  \Q(\sqrt{5}) & ~ j=4, \\
  \Q(\sqrt{2}) & ~ j=7, \\
  \Q           & ~ \text{otherwise.}
 \end{cases}
$$
We fix $\alpha_j\in\OP{Gal}(\bar{\Q}/\Q)$ whose restriction
 to $k_j$ is nontrivial for $j=4,7$.

Let
\begin{gather*}
v_i=(M_i(q_1),\ldots,M_i(q_8)), \\
\widetilde{v_i}=(M_i(q_1),\ldots,M_i(q_8),M_i(\alpha_4(q_4)),M_i(\alpha_7(q_7))),
 \\
V=\begin{pmatrix}v_1\\ \vdots \\v_9\end{pmatrix}, \quad
\widetilde{V}=\begin{pmatrix}\widetilde{v_1}\\ \vdots \\ \widetilde{v_9}\end{pmatrix}.
\end{gather*}
We can check that $\bold{c}V=0$ for
 $\bold{c}=(111375,\ldots,-61440)$ using the computation of Appendix.
Since $\bold{c}\in\Q^9$, we have $\bold{c}\widetilde{V}=0$,
 which implies that the function $f_5=111375 M_1+\cdots-61440 M_9$
 vanishes on the ten points
 $q_1,\ldots,q_8,\alpha_4(q_4),\alpha_7(q_7)$.
By the result of \cite{shimura70,deligne71},
 the image $\OP{Im}F$ of $F$ is defined over $\Q$.
Therefore, $[\alpha_j(q_j)]$ ($j=4,7$) is an
 element in $\OP{Im}(F)$.
By the equality (\ref{degree_nine}), we have $f_5\equiv 0$ on the image
 of $F$.



For the latter part of the theorem, it is sufficient to show that
 the height of $\MF{b}=(Y_2,Y_3,f_2,f_4,f_5)$ is $5$.
This can be checked by observing
\begin{equation*}
 \MF{b}
 =(Y_1^5,Y_2,Y_3,
   Y_1^2+Y_4^2+Y_5^2,Y_1^2 Y_4^2+Y_1^2 Y_5^2+Y_4^2 Y_5^2).
\end{equation*}
%
%
\end{proof}


\begin{rem}
By Lemma~\ref{LEMdegree}, the image of $F$ is of degree $32$.
We have $(f_2,f_4,f_5)\subset\MF{p}_i^2$
 and $\MF{a}=\sqrt{ (f_2,f_4,f_5):\MF{p}_1^2\MF{p}_2^2 }$,
 where
$$
 \MF{p}_1=(Y_1,Y_2^2+Y_3^2-Y_4^2,Y_5),~
 \MF{p}_2=(Y_1,Y_2^2+Y_3^2-Y_5^2,Y_4).
$$
\end{rem}

\begin{rem}
The image of $F$ is contained
 in a Humbert surface $H$ of discriminant $8$
 (cf.\ \cite{hashimoto95,hashimotomurabayashi}).
The defining ideal of $H$ is $(f_2,f_4)$.
\end{rem}

\subsection{Explicit presentation of the inverse period map}

In this subsection, we give an explicit presentation
 of the inverse period map
 of the family $\{ \X_\T \}$.
We use the same notations as in the previous subsection.

We write down a generator of
 the function field $K(\MC{C}_2)$ of $\MC{C}_2$.
Let $r_1,r_2,r_3$ be elements in $K(\MC{C}_2)$ defined by
\begin{equation} \label{def_r}
 r_1=\frac{Y_4 Y_5}{Y_1 Y_3}\biggm|_{\MC{C}_2},~
 r_2=\frac{Y_1^2}{Y_3^2}\biggm|_{\MC{C}_2},~
 r_3=\frac{Y_2^2}{Y_3^2}\biggm|_{\MC{C}_2}.
\end{equation}
Here $(f/g)|_{\MC{C}_2}$ denotes the element in $K(\MC{C}_2)$
 corresponding to an invariant rational function
 $f/g\in\C(Y_1,\ldots,Y_5)^\MF{S}$
 with $F^*g\not\equiv 0$.
(By Lemma~\ref{LEMlinear_relation}, we have $F^* Y_i\not\equiv 0$
 for all $i$.)
By the relation $(f_4-Y_1^2 f_2)/Y_1^2 Y_3^2=0$ on $\OP{Im}(F)$, we have
\begin{equation} \label{elm_r3}
 r_1^2-2r_1-r_2+r_3+1=(r_1-1)^2-r_2+r_3=0.
\end{equation}
By $\deg(F^*(\MC{O}_{\P^4}(2)))=4$,
 we have $[K(\MC{C}_2):\C(r_2-r_3)]\leq 4$,
 and $[K(\MC{C}_2):\C(r_1)]\leq 2$ by (\ref{elm_r3}).

By the relation $Y_1 f_5/Y_3^6=0$ on $\OP{Im}(F)$
 and (\ref{elm_r3}), we obtain
\begin{align}
-104976 r_2^2+252072 r_1^2 r_2-296784 r_1 r_2
+210600 r_2-35721 r_1^4 \notag \\
+227124 r_1^3-448614 r_1^2+297300 r_1-105625=0.
\end{align}
Solving for $r_2$, we have
\begin{equation}
 r_2=\frac{1}{324}\left\{
 (389 r_1^2 - 458 r_1 +325) \pm r_1 (340 r_1 - 20)
 \sqrt{1 - \frac{1}{r_1}} \right\}.
\end{equation}
Therefore, for
\begin{equation} \label{defs}
s=\sqrt{1-\frac{1}{r_1}}=
\frac{324 r_2 - 389 r_1^2
+ 458 r_1 - 325}{r_1(340 r_1 - 20)},
\end{equation}
 we have $[K(\MC{C}_2):\C(s)]=1$.

From the result in Appendix,
 we can check that $s$ takes the values $1,-2$
 at the images of $\xi_1,\xi_5$ in $\MC{C}_2$, respectively.
Hence,
 ${\displaystyle\left( \frac{s-1}{s+2} \right)^2\!\biggm|_{\MC{C}_1}}$
 generates $K(\MC{C}_1)$.
This takes the value $5/4$ at $[\xi_2]$ from the result in Appendix.
Therefore, we have
\begin{equation}
 \frac{t_1}{t_0}
 =-\frac{1}{5}\left( \frac{s-1}{s+2} \right)^2-\frac{1}{4},
\end{equation}
 where $\T=(t_0:t_1)$.
In fact, the right hand side takes values $-1/4,-1/2,\infty$
 at $[\xi_1],[\xi_2],[\xi_5]$, respectively
 (see (\ref{table_sing})).
As a consequence, we have the following theorem.

\begin{thm} \label{THMmain}
Let $i$ be a modular embedding defined by
\begin{equation} \label{def_i}
 i:
 \Omega_T^\circ
 \mathop{\longrightarrow}^{\sim}_{(\ref{identification_omega})}
 \MC{D}_T^\circ
 \mathop{\longrightarrow}^{\sim}_{(\ref{identification_gam})}
 \DOMb
 \mathop{\longrightarrow}_{\PHID}
 \MF{H}_2,
\end{equation}
 where $\Omega_T^\circ$ is
 a connected component of $\Omega_T$ corresponding to $\DOMb$
 and $\PHID$ is defined in Proposition~\ref{PROPembdomain}.
For $\T=(t_0:t_1)\in\P^1$,
 let $[\omega_\T]\in\OP{O}(T)^\circ\backslash\Omega_T^\circ$
 ($\omega_\T\in\Omega_T^\circ$) be the period of $\X_\T$
(see Subsection~(\ref{SUBSECTtransA})).
We set $X_1,\ldots,X_5$ by
\begin{equation}
 \begin{array}{c}
 X_1=i^*\theta_{0000}^4,~
 X_2=i^*\theta_{000\HALF}^4,~
 X_3=i^*\theta_{00\HALF0}^4, \vspace{1ex} \\
 X_4=i^*\theta_{00\HALF\HALF}^4,~
 X_5=i^*\theta_{0\HALF00}^4.
 \end{array}
\end{equation}
Then, we have
\begin{equation*}
 \frac{t_1}{t_0}
 =-\frac{1}{5}
 \left(
 \frac
 {s(\omega_\T)-1}
 {s(\omega_\T)+2}
 \right)^2-\frac{1}{4},
\end{equation*}
where
\begin{gather*}
 r_1=\frac
 {(X_2-X_3+2X_5)(2X_1-X_2+X_3-2X_4-2X_5)}
 {(X_2+X_3)(2X_1-X_2-X_3+2X_4)}, \\
 r_2=\frac
 {(X_2+X_3)^2}
 {(2X_1-X_2-X_3+2X_4)^2}, \\
 s=\frac{324 r_2 - 389 r_1^2
 + 458 r_1 - 325}{r_1(340 r_1 - 20)}.
\end{gather*}
\end{thm}

As we mentioned, the curve $\GAMb\backslash\DOMb$ has
 a canonical $\Q$-structure,
 which is called the Shimura curve $X(15,2)$ (cf.\ \cite{alsinabayer}).
By the above argument,
 we determine the defining equation of $X(15,2)$.

\begin{cor}
The defining equation of the Shimura curve $X(15,2)$ is
\begin{equation} \label{minpoly}
 y^4+6(16x^4-19x^2-1)y^2+(48x^4-55x^2-7)^2=0.
\end{equation}
\end{cor}
\begin{proof}
We only give the outline of the computation.
We determine the function field $K(\MC{C})$
 of $\MC{C}:=\GAMb\backslash\DOMb$.
We set
\begin{equation}
   x=\frac{s-1}{s+2}\biggm|_{\MC{C}},~
 r_4=\frac{Y_4^2 - Y_5^2}{Y_3^2}\biggm|_{\MC{C}},~
 r_5=\frac{Y_2}{Y_1}\biggm|_{\MC{C}}.
\end{equation}
Note that $r_4,r_5$ are relative invariants
 under the action of $\MF{S}$.
We have
\begin{equation}
 r_4^2=\frac{(x-1)^4 (x-2)^2 (x^2-2)}{81 x^4 (x+2)^4},~
 r_5^2=-\frac{(2x-1)^2 (4x^2-5)}{(2x+3)^2 (12x^2+1)}.
\end{equation}
We can check that $K(\MC{C})=\C(x,\sqrt{A}+\sqrt{B})$, where
\begin{equation}
 A=x^2-2,~B=-(4x^2-5)(12x^2+1).
\end{equation}
The minimal polynomial of $y:=\sqrt{A}+\sqrt{B}$ over $\C(x)$ is
\begin{equation*}
 y^4-2(A+B)y^2+(A-B)^2=y^4+6(16x^4-19x^2-1)y^2+(48x^4-55x^2-7)^2,
\end{equation*}
and the assertion follows from this.
\end{proof}

\begin{rem}
It is known that $X(15,2)$ is a hyperelliptic curve
 (\cite{ogg83}).
The hyperelliptic involution is induced by $g_2$ in our notation.
By the above argument,
 the defining equation of $X(15,2)_{\Q(\sqrt{2})}$ is
\begin{equation}
 Y^2=-(5X^2+2X+5) (5X^2-2X+5) (3X^4+26X^2+3).
\end{equation}
Note that the quotient of $X(15,2)$
 by the hyperelliptic involution is
 a genus zero curve which is not isomorphic to $\P^1_\Q$. 
\end{rem}

\noindent
\hskip 1in {\bf Address}\\
\hskip 1in Graduate School of Mathematical Science,\\
\hskip 1in The University of Tokyo, 3-8-1 Komaba, Maguro-ku,\\
\hskip 1in Tokyo, 153-8914, Japan\\
\hskip 1in e-mail: hashi@ms.u-tokyo.ac.jp

\newpage

\def\notsubset{\subset \hskip -0.065in /}
\def\Cal{\mathcal} 
\def\subsetne{\underset{\neq}\subset} 
\def\frak{\mathfrak} 
\def\mod{\operatorname{mod}} 
\def\op{\operatorname{op}} 
\def\ot{\leftarrow} 
\def\<<{\langle } 
\def\>>{\rangle }

\numberwithin{equation}{section} 
\newtheorem{theorem}{Theorem}[section] 
\newtheorem{proposition}[theorem]{Proposition} 
\newtheorem{corollary}[theorem]{Corollary} 
\newtheorem{definition}[theorem]{Definition} 
\newtheorem{conjecture}[theorem]{Conjecture} 
\newtheorem{remark}[theorem]{Remark} 
\newtheorem{lemma}[theorem]{Lemma} 
\newtheorem{fact}[theorem]{Fact} 
\newtheorem{example}[theorem]{Example} 
\newtheorem{axiom}[theorem]{Axiom} 
\newtheorem{property}[theorem]{Property} 
\newtheorem{problem}[theorem]{Problem} 
\newtheorem{notation}[theorem]{Notation} 
\newtheorem{observation}[theorem]{Observation} 

\section*{Appendix: Special values of theta constants at CM points} 

\begin{center}
\bf{Tomohide Terasoma} 
\end{center}

\maketitle 
 
\makeatletter 
\renewcommand{\@evenhead}{\tiny \thepage \hfill 
\hfill}

\renewcommand{\@oddhead}{\tiny \hfill 
 \hfill \thepage} 

\setcounter{section}{1}
\setcounter{subsection}{0}

\subsection{Reduction to elliptic theta values 
and $\lambda$-invariants}

In this appendix, we compute the ratio of the theta
constants evaluated at several CM points in
$\DOMb$ (for the definition, see 
(5.8) or (5.18). Let
$e_1,e_2,e_3$ be the base of $\Cal I(\ORDa)$ defined in
(5.14). 
We used MapleV and risa/asir for the computation.

An element $a_1e_1+a_2e_2+a_3e_3$ is denoted as
$[a_1,a_2,a_3]$.
Let $\widetilde{p_1}, \dots, \widetilde{p_8}$ be elements of $\ORDa$ defined by
\begin{align*}
\widetilde{p_1}=[2,0,1],\widetilde{p_2}=[4,-4,3],
\widetilde{p_3}=[2,4,3], 
\widetilde{p_4}=[0,0,1], \\
\widetilde{p_5}=[2,-2,3],
\widetilde{p_6}=[5,1,3],
\widetilde{p_7}=[1,1,1],
\widetilde{p_8}=[5,-5,3].
\end{align*}
Then $\widetilde{p_i}^2\in \bold R_{ <0}$ and
$p_i=\displaystyle\frac{\widetilde{p_i}}
{\sqrt{-\widetilde{p_i}^2}}$ is an element of 
$\DOMb$ defined in (5.8).

Let $\tau_1, \dots, \tau_8$ be the images of $p_1, \dots,p_8$
under the map 
$$
\PHID:\DOMb \to \frak H_2
$$
in (6.7). Then we have 
$$
\tau_1=\begin{pmatrix}\frac{3\,\sqrt{3}\,i+5}{2} & 
-\frac{\sqrt{3}\,i+1}{2 } \\ -\frac{\sqrt{3}\,i+1}{2} 
& \frac{\sqrt{3}\,i-1}{2} \\  \end{pmatrix},
\tau_2=\begin{pmatrix}\frac{\sqrt{7}\,i+1}{2} & 0 \\ 0 & 
\frac{\sqrt{7}\,i -3}{2} \\ \end{pmatrix},
$$
$$
\tau_3=\begin{pmatrix}\frac{3\,\sqrt{7}\,i-3}{2} & 
-\frac{\sqrt{7}\,i-3}{2 } \\ -\frac{\sqrt{7}\,i-3}{2} 
& \frac{\sqrt{7}\,i-1}{2} \\  \end{pmatrix},
\tau_4=\begin{pmatrix}\frac{\sqrt{15}\,i-1}{2} & 1 \\ 1 & \frac{\sqrt{15}
  \,i-1}{2} \\ \end{pmatrix},
$$
$$
\tau_5=\begin{pmatrix}\frac{\sqrt{13}\,i}{2} & \frac{1}{2} \\ 
\frac{1}{2}  & \frac{\sqrt{13}\,i-2}{2} \\ \end{pmatrix},
\tau_6=\begin{pmatrix}\frac{\sqrt{22}\,i+2}{2} & -\frac{1}{2} \\ -\frac{1 }{2}
& \frac{\sqrt{22}\,i-2}{4} \\ \end{pmatrix},
$$
$$
\tau_7=\begin{pmatrix}\frac{\sqrt{30}\,i+2}{2} & \frac{1}{2} \\ \frac{1}{2 } &
\frac{\sqrt{30}\,i-2}{4} \\ \end{pmatrix},
\tau_8=\begin{pmatrix}\frac{\sqrt{10}\,i+3}{4} & -\frac{1}{4} 
\\ -\frac{1 }{4} & \frac{\sqrt{10}\,i-7}{4} \\ \end{pmatrix}.
$$
In this appendix, we compute the ratio of theta constants defined in (7.1).
\begin{equation}
\label{ratio}
(\theta_{v_1}(\tau_i)^4:\theta_{v_2}(\tau_i)^4:
\theta_{v_3}(\tau_i)^4:\theta_{v_4}(\tau_i)^4:
\theta_{v_5}(\tau_i)^4)
\end{equation}
where $v_1=(0,0,0,0)$,$v_2=(0,0,0,1/2)$,$v_3=(0,0,1/2,0)$,
$v_4=(0,0,1/2,1/2)$ and $v_5=(0,1/2,0,0)$.
We set
$$P=\begin{pmatrix}1 & 1 \\ 0 & 1 \\ \end{pmatrix}$$
and define $\overline{\tau_1},\overline{\tau_3}$ by
$$
\overline{\tau_1}=P\tau_1^tP=\begin{pmatrix}\sqrt{3}\,i+1 & -1 \\ -1 & 
\frac{\sqrt{3}\,i-1}{2}\end{pmatrix},  
\overline{\tau_3}=P\tau_3^tP=\begin{pmatrix}\sqrt{7}\,i+1 & 1 \\ 1 & 
\frac{\sqrt{7}\,i-1}{2} \\  \end{pmatrix}.
$$
Since
$\theta_{v_2}(\tau_j)=\theta_{v_4}(\overline{\tau_j})$,
$\theta_{v_4}(\tau_j)=\theta_{v_2}(\overline{\tau_j})$ and
$\theta_{v_i}(\tau_j)=\theta_{v_i}(\overline{\tau_j})$ for
$j=1,3$ and $i=1,3,5$, the computations of $\theta_{v_i}(\tau_j)$
is reduced to that of $\theta_{v_i}(\overline{\tau_j})$.

Let $\tau$ be an element of $\{\overline{\tau_1},
\tau_2,\overline{\tau_3},\tau_6,\tau_5,\tau_7,
\tau_8\,\tau_4\}$. Then the $(1,2)$-component $\tau_{1,2}$ of $\tau$ is
an element in $\frac{1}{4}\bold Z$.
We use the following formula for theta functions.
\begin{equation}
\label{twice formula for g=2theta}
\theta_{a_1,a_2,b_1,b_2}(\tau)^2=
\sum_{\substack{\epsilon_1=0,1/2 \\\epsilon_2=0,1/2  }}
s(\epsilon)\theta_{\epsilon_1,\epsilon_2,b_1,b_2}(2\tau)
\theta_{\epsilon_1+a_1,\epsilon_2+a_2,b_1,b_2}(2\tau),
\end{equation}
where 
$
s(\epsilon_1,\epsilon)=(-1)^{4(\alpha_1+\epsilon_1)\beta_1+
4(\alpha_2+\epsilon_2)\beta_2}
$.
By the above formula, the computations 
of $\theta_{a_1,a_2,b_1,b_2}(\tau_i)$
is reduced to those of 
$\theta_{a_1,a_2,b_1,b_2}(\widetilde{\tau})$,
where $\widetilde{\tau}\in\{  
\overline{\tau_1},
\tau_2,\overline{\tau_3},\tau_6,2\tau_5,2\tau_7,2\tau_8,4\tau_4\}$.
We recall the following classical formulae for theta constants.
We use the elliptic theta constant:
$$
\theta_{m_1,m_2}(\tau)=\sum_{n\in \bold Z}\exp (2\pi
     i(\frac{1}{2}(n+m_1)\tau (n+m_1)+(n+m_1)m_2)),
$$
for $m_1,m_2\in \bold Q,\tau \in \frak H=\{\tau\in\bold C|
\OP{Im}(\tau)>0 \}$.

\begin{proposition}
\begin{enumerate}
\item
Let $\widetilde{\tau}=
\left(\begin{matrix}\widetilde{\tau}_{11},\widetilde{\tau}_{12} \\
\widetilde{\tau}_{21},\widetilde{\tau}_{22}\end{matrix}\right)$
be an element of the Siegel upper half plane such that 
$\widetilde{\tau_{12}}=\widetilde{\tau_{21}}
\in \bold Z$ and $a_1,a_2,b_1,b_2\in \{0,1/2\}$.
Then we have
\begin{equation}
\label{reduction to elliptic case}
\theta_{a_1,a_2,b_1,b_2}(\widetilde{\tau})=
\theta_{a_1,b_1+\widetilde{\tau}_{12}a_2}(\widetilde{\tau}_{11})
\theta_{a_2,b_2+\widetilde{\tau}_{12}a_1}(\widetilde{\tau}_{22}).
\end{equation}
\item
We have
\begin{align}
\label{jacobi}
& \theta_{0,1/2}(\tau)^4+
\theta_{1/2,0}(\tau)^4=\theta_{0,0}(\tau)^4,\quad
\theta_{1/2,1/2}(\tau)=0, \\
\label{twice for lambda}
& \theta_{0,0}(2\tau)^2=\frac{1}{2}(\theta_{0,0}(\tau)^2
+\theta_{0,1/2}(\tau)^2),\quad
 \theta_{0,1/2}(2\tau)^2=\theta_{0,0}(\tau)\theta_{0,1/2}(\tau).
\end{align}
\item
We define the $\lambda$-invariant by
$\displaystyle
\lambda(\tau)=\frac{\theta_{0,1/2}(\tau)^4}{\theta_{0,0}(\tau)^4}$.
Then we have
\begin{equation}
\label{j and lambda}
\displaystyle
j(\tau)=256\frac{(\lambda^2-\lambda+1)^3}{\lambda^2(\lambda-1)^2}
\end{equation}
where $j(\tau)$ is the $j$-invariant of $\tau$.
\end{enumerate}
\end{proposition}
By the above proposition, the computations of the ratio (\ref{ratio})
is reduced to those of
$j(\widetilde{\tau}_{11})$ and
$j(\widetilde{\tau}_{22})$.

\subsection{Computations of $\lambda$-invariants}
Let $D$ be a positive integer, $\Cal O_D$ the integer ring
of $\bold Q_(\sqrt{-d})$, $I_1, \dots ,I_k$ a complete
representative of ideal class of $\Cal O_D$.
Then the elliptic curve $E_i=\bold C/I_i$ depends only on the
ideal class of $I_i$. These elliptic curves are defined over the Hilbert 
class field $F$ of $\bold Q_D$.
We choose a $\bold Z$ base $\alpha,\beta$ of $I_i$ such that 
$\tau_i=Im(\alpha/\beta)>0$.
Then the $j$-invariant $j(E_i)$ of the elliptic curve $E_i$
is equal to $j(\tau_i)$.
By the theory of complex multiplication, 
the set $\{ j(E_1), \dots, j(E_k)\}$
of $j$-invariants is stable under the action of $Gal(F/\bold Q)$,
and $\lambda(E_i)$is an algebraic integer.
Therefore elementary symmetric 
polynomials $s_1, \dots, s_k$ of $\{ j(E_1),\dots, j(E_k)\}$ are
integers. Using this integrality, we can compute $s_1, \dots, s_k$ by 
numerical computation.
Thus, the set $\{j(E_1),\dots, j(E_k)\}$ is the set of solutions
of the equation $x^k-s_1x^{k-1}\cdots +(-1)^ks_k=0$. 

We compute the equation satisfied by the $j$-invariants of
$E_i$'s for $D=3,7,10,13,15,22,30$.

\vskip 0.3in
\hskip 0in
\begin{tabular}{c|c}
\hline
equation & 
the set of solutions \\
\hline\hline
$x=0$ & $\{j(\frac{-1+\sqrt{-3}}{2})\}$  \\
\hline
$x+3375=0$ & $\{j(\frac{1+\sqrt{-7}}{2})\}$ \\ 
\hline
$x^2+191025x-121287375$ &  
$
\begin{matrix}
\{j((\sqrt{-15}+1)/2), \\
j((\sqrt{-15}+1)/4)\}
\end{matrix}
$ 
\\
\hline
$x^2-6896880000x-567663552000000$ &
$
\begin{matrix}
\{j(\sqrt{-13}), \\
j((\sqrt{-13}+1)/2)\}
\end{matrix}
$ 
\\
\hline
$
\begin{matrix}
x^2-6294842640000x \\
+15798135578688000000
\end{matrix}
$ &
$\{j(\sqrt{-22}),j(\sqrt{-22}/2) \}
$ \\
\hline
$\begin{matrix}
x^4-883067971104000x^3 \\
+26329406807264910336000x^2 \\
-2588458316335175909376000000x \\
+4934510722321469030006784000000 \\
\end{matrix}$ & 
$\begin{matrix}
\{j(\sqrt{-30}),
j(\sqrt{-30}/2), \\
j(\sqrt{-30}/3),
j(\sqrt{-30}/5)\}
\end{matrix}$ \\
\hline
$x^2-425692800x+9103145472000$ &
$\{j(\sqrt{-10}/2),j(\sqrt{-10})\}$\\
\hline
\end{tabular}

\vskip 0.1in
By solving the equation, we have
\begin{align*}
j(\frac{\sqrt{-3}-1}{2})&=0,
j(\frac{\sqrt{-7}+1}{2})=
j(\frac{\sqrt{-7}-3}{2})=
j(\frac{\sqrt{-7}-1}{2})=-3375, \\
j(\frac{\sqrt{-15}-1}{2})&=\frac{-85995\sqrt{5}-191025}{2}, \\
j(\sqrt{-13})&=j(\sqrt{-13}-2)=
956448000\sqrt{13}+3448440000, \\
j(\sqrt{-22}+2)&=3147421320000+2225561184000\sqrt{2},\\
j(\frac{\sqrt{-22}-2}{2})&=3147421320000-2225561184000\sqrt{2}, \\
j(\sqrt{-30}+2)&=98729993940480\sqrt{5}+156105837619200\sqrt{2}\\
&+69812648236800\sqrt{10}+220766992776000,\\
j(\frac{\sqrt{-30}-2}{2})&=-98729993940480\sqrt{5}+156105837619200\sqrt{2}\\
&-69812648236800\sqrt{10}+220766992776000, \\
j(\sqrt{-10}+3)&=
j(\sqrt{-10}-7)=95178240\sqrt{5}+212846400.
\end{align*}

\subsection{$j$-invariants to $\lambda$-invariants and theta values}
Using the identity (\ref{j and lambda}), the $\lambda$-invariants
are obtained from the $j$-invariants. Using the $\lambda$-invariants, we
compute the ratio's of theta constants.

(1) The case of $\tau_1$.
\begin{align*}
\lambda(\frac{\sqrt{-3}-1}{2})&=(\sqrt{-3}+1)/2 \text{ and }
\lambda(\sqrt{-3}+1)=8-4\sqrt{3} \text{ by (\ref{twice for lambda})}.
\end{align*}
By using the formula (\ref{reduction to elliptic case}), we have
\begin{align*}
&(\theta_{v_1}(\tau_1)^4:\theta_{v_2}(\tau_1)^4:
\theta_{v_3}(\tau_1)^4:\theta_{v_4}(\tau_1)^4:
\theta_{v_5}(\tau_1)^4) \\
&=
(1:4\sqrt{-3}-6 i-2\sqrt{3}+4:8-4\sqrt{3}:
\sqrt{-3}/2+1/2:-4\sqrt{-3}+6 i-2\sqrt{3}+4
).
\end{align*}

(2) The case of $\tau_2$.
\begin{align*}
\lambda(\frac{\sqrt{-7}+1}{2})&=
\lambda(\frac{\sqrt{-7}-3}{2})=(-3\sqrt{-7}+31)/32.
\end{align*}
By using the formula (\ref{reduction to elliptic case}), we have
\begin{align*}
&(\theta_{v_1}(\tau_2)^4:\theta_{v_2}(\tau_2)^4:
\theta_{v_3}(\tau_2)^4:\theta_{v_4}(\tau_2)^4:
\theta_{v_5}(\tau_2)^4) \\
=&
(1:31/32-3\sqrt{-7}/32:31/32-3\sqrt{-7}/32:\\
&449/512-93\sqrt{-7}/512:3\sqrt{-7}/32+1/32).
\end{align*}

(3) The case of $\tau_3$.
\begin{align*}
\lambda(\frac{\sqrt{-7}-1}{2})&=(3\sqrt{-7}+31)/32,
\text{ and }
\lambda(\sqrt{-7}+1)=-48\sqrt{7}+128 \text{ by (\ref{twice for lambda})}
\end{align*}
By using the formula (\ref{reduction to elliptic case}), we have
\begin{align*}
&(\theta_{v_1}(\tau_3)^4:\theta_{v_2}(\tau_3)^4:
\theta_{v_3}(\tau_3)^4:\theta_{v_4}(\tau_3)^4:
\theta_{v_5}(\tau_3)^4) \\
=&(
32: 384\sqrt{-7}-1008 \sqrt{-1}-1488\sqrt{7}+3968: 4096-1536\sqrt{7}: \\
& 3\sqrt{-7}+31: -384\sqrt{-7}+1008 \sqrt{-1}-48\sqrt{7}+128).
\end{align*}

(4) The case of $\tau_4$.
\begin{align*}
\lambda(\frac{\sqrt{-15}-1}{2})&=17/64+21\sqrt{5}/64+17\sqrt{-3}/64
-7\sqrt{-15}/64, 
\end{align*}
By using the formula (\ref{reduction to elliptic case}), we have
\begin{align*}
&(\theta_{v_1}(\tau_4)^4:\theta_{v_2}(\tau_4)^4:
\theta_{v_3}(\tau_4)^4:\theta_{v_4}(\tau_4)^4:
\theta_{v_5}(\tau_4)^4) \\
=&(
1024: -112\sqrt{-15}+272\sqrt{-3}+336\sqrt{5}+272: 
-112\sqrt{-15}+272\sqrt{-3}+336\sqrt{5}+272: \\
&-223\sqrt{-3}+119\sqrt{-15}+357\sqrt{5}+223: 
112\sqrt{-15}-272\sqrt{-3}-336\sqrt{5}+752
).
\end{align*}

(5) The case of $\tau_5$.
\begin{align*}
\lambda(\sqrt{-13})&=
\lambda(\sqrt{-13}-2)=1/2+3\sqrt{-18+5\sqrt{13}}.
\end{align*}
Therefore, we have
\begin{align*}
\theta(2\tau_5,0,1/2,0,0)^4/\theta(2\tau_5,0,0,0,0)^4=& 
\theta(2\tau_5,1/2,0,0,0)^4/\theta(2\tau_5,0,0,0,0)^4 \\
=& 649/4-45\sqrt{13} \\
=& (5\sqrt{13}-18)^2/4, \\
\theta(2\tau_5,1/2,1/2,0,0)^4/\theta(2\tau_5,0,0,0,0)^4=&0.
\end{align*}
By using the formula (\ref{twice formula for g=2theta}), we have
\begin{align*}
&(\theta_{v_1}(\tau_5)^4:\theta_{v_2}(\tau_5)^4:
\theta_{v_3}(\tau_5)^4:\theta_{v_4}(\tau_5)^4:
\theta_{v_5}(\tau_5)^4) \\
=&(
1:614-170\sqrt{13}:686-190\sqrt{13}:1:36-10\sqrt{13}).
\end{align*}

(6) The case of $\tau_6$.
\begin{align*}
\lambda(\sqrt{-22}+2)&=
2(\sqrt{2}+1)^6(10-3\sqrt{11})(3\sqrt{11}-7\sqrt{2}),
\\
1-\lambda(\sqrt{-22}+2)&=
(10-3\sqrt{11})^2(3\sqrt{11}-7\sqrt{2})^2,\\
\lambda(\frac{\sqrt{-22}-2}{2})&=
(3\sqrt{11}+10)^2(3\sqrt{11}-7\sqrt{2})^2, \\
1-\lambda(\frac{\sqrt{-22}-2}{2})&=
-2(\sqrt{2}-1)^6(3\sqrt{11}+10)(3\sqrt{11}-7\sqrt{2}). 
\end{align*}
Therefore, we have
\begin{align*}
\theta(2\tau_6,0,1/2,0,0)^4/\theta(2\tau_6,0,0,0,0)^4=& 
-4(3\sqrt{11}-7\sqrt{2})^2, \\
\theta(2\tau_6,1/2,0,0,0)^4/\theta(2\tau_6,0,0,0,0)^4  =& 
(3\sqrt{11}-7\sqrt{2})^4,
\\
\theta(2\tau_6,1/2,1/2,0,0)^4/\theta(2\tau_6,0,0,0,0)^4=&0.
\end{align*}
By using the formula (\ref{twice formula for g=2theta}), we have
\begin{align*}
&(\theta_{v_1}(\tau_6)^4:\theta_{v_2}(\tau_6)^4:
\theta_{v_3}(\tau_6)^4:\theta_{v_4}(\tau_6)^4:
\theta_{v_5}(\tau_6)^4) \\
=&(
1176\sqrt{-11}-2758\sqrt{-2}-4074\sqrt{22}+19109: \\
&1182\sqrt{-11}-2772\sqrt{-2}-4116\sqrt{22}+19306: \\
&-1182\sqrt{-11}+2772\sqrt{-2}-4116\sqrt{22}+19306: \\
&-1176\sqrt{-11}+2758\sqrt{-2}-4074\sqrt{22}+19109: \\
&14\sqrt{-2}-6\sqrt{-11}
).
\end{align*}

(7) The case of $\tau_7$.
\begin{align*}
\lambda(\sqrt{-30}+2)&=
2(-\sqrt{5}+\sqrt{6})((1+\sqrt{5})/2)^6
(\sqrt{10}+3)^2(4-\sqrt{15})(2-\sqrt{3})(5-2\sqrt{6}), \\
1-\lambda(\sqrt{-30}+2)&=
(2-\sqrt{3})^2(5-2\sqrt{6})^2
(4-\sqrt{15})^2(-\sqrt{5}+\sqrt{6})^2, \\
\lambda(\frac{\sqrt{-30}-2}{2})&=
(2+\sqrt{3})^2(5+2\sqrt{6})^2(4-\sqrt{15})^2(\sqrt{5}-\sqrt{6})^2, \\
1-\lambda(\frac{\sqrt{-30}-2}{2})&=
2(\sqrt{5}-\sqrt{6})((1-\sqrt{5})/2)^6
(-\sqrt{10}+3)^2(4-\sqrt{15})(2+\sqrt{3})(5+2\sqrt{6}).
\end{align*}
Therefore, we have
\begin{align*}
\theta(2\tau_7,0,1/2,0,0)^4/\theta(2\tau_7,0,0,0,0)^4=& 
-4(\sqrt{6}-\sqrt{5})^2(4-\sqrt{15})^2,
\\
\theta(2\tau_7,1/2,0,0,0)^4/\theta(2\tau_7,0,0,0,0)^4  =& 
(\sqrt{6}-\sqrt{5})^4(4-\sqrt{15})^4,
\\
\theta(2\tau_7,1/2,1/2,0,0)^4/\theta(2\tau_7,0,0,0,0)^4=&0.
\end{align*}
By using the formula (\ref{twice formula for g=2theta}), we have
\begin{align*}
&(\theta_{v_1}(\tau_7)^4:\theta_{v_2}(\tau_7)^4:
\theta_{v_3}(\tau_7)^4:\theta_{v_4}(\tau_7)^4:
\theta_{v_5}(\tau_7)^4). \\
&(-3966\sqrt{-10}-5608\sqrt{-5}+5120\sqrt{-6}
+7240\sqrt{-3} \\
&-21038\sqrt{30}-29752\sqrt{15}
+81480\sqrt{2}+115229: \\
&-3972\sqrt{-10}-5616\sqrt{-5}
+5128\sqrt{-6}+7250\sqrt{-3}\\
&-21100\sqrt{30}
-29840\sqrt{15}+81720\sqrt{2}+115570: \\
&3972\sqrt{-10}
+5616\sqrt{-5}-5128\sqrt{-6}-7250\sqrt{-3}\\
&-21100\sqrt{30}
-29840\sqrt{15}+81720\sqrt{2}+115570: \\
&3966\sqrt{-10}+5608\sqrt{-5}-5120\sqrt{-6}
-7240\sqrt{-3}\\
&-21038\sqrt{30}
-29752\sqrt{15}
+81480\sqrt{2}+115229: \\
&6\sqrt{-10}+8\sqrt{-5}
-8\sqrt{-6}-10\sqrt{-3}).
\end{align*}

(8) The case of $\tau_8$.
\begin{align*}
\lambda(\sqrt{-10}+3)&=
\lambda(\sqrt{-10}-7)=
(\sqrt{10}+3)
(1/2-\sqrt{5}/2)^6(1+\sqrt{2})^2/2, \\
1-\lambda(\sqrt{-10}+3)&=
1-\lambda(\sqrt{-10}-7)=
(-\sqrt{10}+3)(1/2-\sqrt{5}/2)^6(1-\sqrt{2})^2/2. 
\end{align*}
Therefore, we have
\begin{align*}
\theta(4\tau_8,0,1/2,0,0)^4/\theta(4\tau_8,0,0,0,0)^4=& 
\theta(4\tau_8,1/2,0,0,0)^4/\theta(4\tau_8,0,0,0,0)^4 \\
 =& 
-(\sqrt{5}/2-1/2)^{12}/4,
\\
\theta(4\tau_8,1/2,1/2,0,0)^4/\theta(4\tau_8,0,0,0,0)^4=&0.
\end{align*}
By using the formula (\ref{twice formula for g=2theta}), we have
\begin{align*}
\theta(2\tau_8,0,1/2,0,0)^2/\theta(2\tau_8,0,0,0,0)^2=& 
\theta(2\tau_8,1/2,0,0,0)^2/\theta(2\tau_8,0,0,0,0)^2 \\
=& 
(\sqrt{5}/2-1/2)^3(1+\sqrt{-1}), \\
\theta(2\tau_8,1/2,1/2,0,0)^2/\theta(2\tau_8,0,0,0,0)^2=&
-(\sqrt{5}/2-1/2)^6.
\end{align*}
Again by using the formula (\ref{twice formula for g=2theta}), we have
\begin{align*}
&(\theta_{v_1}(\tau_8)^4:\theta_{v_2}(\tau_8)^4:
\theta_{v_3}(\tau_8)^4:\theta_{v_4}(\tau_8)^4:
\theta_{v_5}(\tau_8)^4) \\
=&
(-16\sqrt{-5}+36\sqrt{-1}-12\sqrt{5}+27:94-42\sqrt{5}:
14-6\sqrt{5}:\\
&16\sqrt{-5}-36\sqrt{-1}-12\sqrt{5}+27:
\\
&-16\sqrt{-5}+36\sqrt{-1}+8\sqrt{10}+18\sqrt{5}-18\sqrt{2}-40).
\end{align*}

\end{document}